\documentclass[a4paper,11pt]{amsart}

\usepackage{amssymb,amsthm,amsmath,amscd,amsfonts,bbm,mathrsfs}
\usepackage{stmaryrd}
\usepackage{hyperref}
\usepackage{xcolor}
\hypersetup{
  bookmarksnumbered=true,
  colorlinks=false,
  bookmarksopen=true,
  bookmarksdepth=2,
  pageanchor=true}

\usepackage[cmtip,all]{xy}

\theoremstyle{plain}
\newtheorem{LeABC}{Lemma}
\newtheorem{ThmABC}[LeABC]{Theorem}
\newtheorem{PrABC}[LeABC]{Proposition}
\newtheorem{Lemma}[equation]{Lemma}
\newtheorem{Thm}[equation]{Theorem}
\newtheorem{Prop}[equation]{Proposition}
\newtheorem{Cor}[equation]{Corollary}
\theoremstyle{definition}
\newtheorem{Defn}[equation]{Definition}
\newtheorem{Const}[equation]{Construction}
\newtheorem*{ToDo}{\color{blue}ToDo}
\theoremstyle{remark}
\newtheorem{Remark}[equation]{Remark}
\newtheorem*{Remark*}{Remark}
\newtheorem{Obs}[equation]{Oberservation}
\newtheorem{Notation}[equation]{Notation}
\newtheorem{Terminology}[equation]{Terminology}
\newtheorem*{Notation*}{Notation}
\newtheorem{Ass}[equation]{Assumption}
\newtheorem{Example}[equation]{Example}

\newcommand{\bLe}{\begin{Lemma}}
\newcommand{\eLe}{\end{Lemma}}
\newcommand{\bTh}{\begin{Thm}}
\newcommand{\eTh}{\end{Thm}}
\newcommand{\bPr}{\begin{Prop}}
\newcommand{\ePr}{\end{Prop}}
\newcommand{\bCo}{\begin{Cor}}
\newcommand{\eCo}{\end{Cor}}
\newcommand{\bCon}{\begin{Const}}
\newcommand{\eCon}{\end{Const}}
\newcommand{\bDe}{\begin{Defn}}
\newcommand{\eDe}{\end{Defn}}
\newcommand{\bAs}{\begin{Ass}}
\newcommand{\eAs}{\end{Ass}}
\newcommand{\bReX}{\begin{Remark*}}
\newcommand{\eReX}{\end{Remark*}}
\newcommand{\bRe}{\begin{Remark}}
\newcommand{\eRe}{\end{Remark}}
\newcommand{\bOb}{\begin{Obs}}
\newcommand{\eOb}{\end{Obs}}
\newcommand{\bEx}{\begin{Example}}
\newcommand{\eEx}{\end{Example}}
\newcommand{\bNo}{\begin{Notation}}
\newcommand{\eNo}{\end{Notation}}
\newcommand{\bTe}{\begin{Terminology}}
\newcommand{\eTe}{\end{Terminology}}
\newcommand{\bNoX}{\begin{Notation*}}
\newcommand{\eNoX}{\end{Notation*}}
\newcommand{\bToDo}{\begin{ToDo}\color{blue}}
\newcommand{\eToDo}{\end{ToDo}}

\newcommand{\bproof}{\begin{proof}}
\newcommand{\eproof}{\end{proof}}
\newcommand{\beqn}{\begin{equation}}
\newcommand{\eeqn}{\end{equation}}

\renewcommand{\u}{\underline}
\newcommand{\uM}{M} 
\newcommand{\uN}{N} 

\numberwithin{equation}{section}

\usepackage{relsize}
\usepackage[bbgreekl]{mathbbol}

\DeclareSymbolFontAlphabet{\mathbb}{AMSb}
\DeclareSymbolFontAlphabet{\mathbbl}{bbold}

\newcommand{\AAA}{\mathcal{A}}
\newcommand{\BBB}{\mathcal{B}}

\newcommand{\FFF}{\mathcal{F}}

\newcommand{\MMM}{\mathcal{M}}

\newcommand{\OOO}{\mathcal{O}}

\newcommand{\SSS}{\mathcal{S}}

\newcommand{\XXX}{\mathcal{X}}
\newcommand{\YYY}{\mathcal{Y}}

\newcommand{\Fa}{\mathfrak{a}}
\newcommand{\Fb}{\mathfrak{b}}

\newcommand{\Fm}{\mathfrak{m}}

\newcommand{\FY}{\mathfrak{Y}}
\newcommand{\FZ}{\mathfrak{Z}}

\newcommand{\DD}{{\mathbb{D}}}
\newcommand{\FF}{{\mathbb{F}}}
\newcommand{\GG}{{\mathbb{G}}}

\newcommand{\LL}{{\mathbb{L}}}

\newcommand{\NN}{{\mathbb{N}}}
\newcommand{\QQ}{{\mathbb{Q}}}

\newcommand{\WW}{{\mathbb{W}}}
\newcommand{\ZZ}{{\mathbb{Z}}}

\DeclareMathOperator{\coker}{coker}
\DeclareMathOperator{\colim}{colim}
\DeclareMathOperator{\cone}{cone}
\DeclareMathOperator{\can}{can}
\DeclareMathOperator{\crys}{cris}
\DeclareMathOperator{\csyn}{csyn}
\DeclareMathOperator{\essemiperf}{essp}
\DeclareMathOperator{\fppf}{fppf}
\DeclareMathOperator{\fpqc}{fpqc}
\DeclareMathOperator{\gr}{gr}
\DeclareMathOperator{\height}{ht}
\DeclareMathOperator{\id}{id}
\DeclareMathOperator{\nil}{nil}
\DeclareMathOperator{\nnil}{-nil}
\DeclareMathOperator{\plain}{plain}
\DeclareMathOperator{\perf}{perf}
\DeclareMathOperator{\qrsp}{qrsp}
\DeclareMathOperator{\qsyn}{qsyn}
\DeclareMathOperator{\rel}{rel}
\DeclareMathOperator{\red}{red}
\DeclareMathOperator{\rk}{rk}
\DeclareMathOperator{\syn}{syn}
\DeclareMathOperator{\semiperf}{sp}
\DeclareMathOperator{\uni}{uni}

\DeclareMathOperator{\Ab}{Ab}
\DeclareMathOperator{\Aff}{Aff}
\DeclareMathOperator{\Aut}{Aut}
\DeclareMathOperator{\Bil}{Bil}
\DeclareMathOperator{\Ch}{Ch}
\DeclareMathOperator{\nCrys}{nCris}
\DeclareMathOperator{\BT}{BT}
\DeclareMathOperator{\Disp}{Disp}
\DeclareMathOperator{\sDisp}{{}^sDisp}
\DeclareMathOperator{\End}{End}
\DeclareMathOperator{\Ext}{Ext}
\DeclareMathOperator{\FG}{FG}
\DeclareMathOperator{\GL}{GL}
\DeclareMathOperator{\Hom}{Hom}
\DeclareMathOperator{\Isom}{Isom}
\DeclareMathOperator{\LF}{LF}
\DeclareMathOperator{\Lie}{Lie}
\DeclareMathOperator{\Mod}{Mod}
\DeclareMathOperator{\Nil}{Nil}
\DeclareMathOperator{\Pic}{Pic}
\DeclareMathOperator{\Spf}{Spf}
\DeclareMathOperator{\Spec}{Spec}
\DeclareMathOperator{\Tor}{Tor}
\DeclareMathOperator{\Win}{Win}

\newcommand{\Vnil}{V{\nnil}}
\newcommand{\limperf}{^{\perf}}
\newcommand{\colimperf}{_{\perf}}

\newcommand{\sC}{{}^{s}C}
\newcommand{\sD}{{}^{s\!}D}

\newcommand{\sM}{{}^{s\!}M}

\newcommand{\sMMM}{{}^{s}\!\MMM}
\newcommand{\sFZ}{{}^s\FZ}
\newcommand{\sW}{{}^sW}

\newcommand{\usW}{{}^s\u W}
\newcommand{\bsW}{{}^s\bar W}

\newcommand{\srho}{{}^{s\!}\rho}

\newcommand{\cclearpage}{}

\begin{document}

\title{Sheared displays and $p$-divisible groups}
\author{Manuel Hoff and Eike Lau}
\date{\today}
\address{Fakult\"{a}t f\"{u}r Mathematik,
Universit\"{a}t Bielefeld, D-33501 Bielefeld}

\begin{abstract}
We develop a Dieudonn\'e theory for $p$-divisible groups over general $p$-nilpotent rings using sheared Witt vectors.
\end{abstract}

\maketitle

\setcounter{tocdepth}{1}
\tableofcontents

\cclearpage
\section{Introduction}

Let $p$ be a prime, and let $R$ be a $p$-nilpotent ring.
The theory of displays gives an equivalence between infinitesimal $p$-divisible groups and $V$-nilpotent displays over $R$. In this article, we show that \emph{general} $p$-divisible groups are equivalent to sheared displays, the analogue of displays where the ring of Witt vectors $W(R)$ is replaced by its sheared version $\sW(R)$, introduced recently in \cite{Drinfeld:Ring-stacks-conjecturally, BMVZ}.  
Such an equivalence was conjectured by Drinfeld in \cite{Drinfeld:Talk-Manin}.

A recent result of Gardner--Madapusi \cite{Gardner-Madapusi} shows that $p$-divisible groups over $R$ are equivalent to vector bundle $F$-gauges\footnote{Vector bundle $F$-gauges over $R$ are vector bundles over the syntomification $R^{\syn}$ defined by Bhatt--Lurie and Drinfeld.} with Hodge-Tate weights in $\{0,1\}$. A priori the latter form an $\infty$-category, which ultimately turns out to be classical, while sheared displays are classical by definition. Most of the present work does not use the language of higher categories.

The above result of \cite{Gardner-Madapusi} follows from a similar classification of $n$-truncated Barsotti--Tate groups, and the conjectures of \cite{Drinfeld:Talk-Manin} also include a truncated version, which, however, involves mild animation a priori. For this reason, we will restrict to the case of $p$-divisible groups here and return to the truncated version later.

Before stating the main results in Theorems \ref{Th:Intro-sFZR} and \ref{Th:Intro-sDR}, we need some preparations.

\subsection{Displays and $p$-divisible groups}

Let us recall two classical constructions, whose sheared analogues will be given below. 
First, there is a functor $\FZ_R\colon\Disp(R)\to\FG(R)$ from displays to commutative (smooth) formal groups, defined in \cite{Zink:The-display} by an explicit formula. Moreover, the Dieudonn\'e crystal yields a functor $D_R\colon\BT(R)\to\Disp(R)$ from $p$-divisible groups to displays such that the composition $\FZ_R\circ D_R\colon\BT(R)\to\FG(R)$ is the formal completion, by \cite{Lau:Smoothness}. These functors restrict to mutually inverse equivalences $\FZ_{R,\nil}$ and $D_{R,\nil}$ between infinitesimal $p$-divisible groups and $V$-nilpotent displays.

\subsection{Sheared Witt vectors}

As usual let $\hat W(R)\subseteq W(R)$ denote the ideal of Witt vectors whose entries are nilpotent and eventually zero.
Following \cite{Drinfeld:Ring-stacks-conjecturally, BMVZ} we set $Q=W/\hat W$ as an fpqc sheaf  of rings  on the category of $p$-nilpotent rings. The sheaf of sheared Witt vectors is defined as
\[
\sW=W\times_QQ\limperf,
\] 
where $Q\limperf$ is the limit perfection of $Q$ with respect to the Witt vector Frobenius; the terminology reflects the role of $\sW$ in the theory of sheared prismatisation of \cite{BKMVZ}.
For example, if $R$ is an Artinian local ring with perfect residue field $k$, then $\sW(R)=W(k)\oplus\hat W(\Fm_R)$ as a subring of $W(R)$. 
This ring appears in \cite{Zink:Dieudonne}.
If $R$ is a semiperfect $\FF_p$-algebra, i.e.\ the Frobenius endomorphism is surjective, then $\sW(R)=W(R\limperf)/\hat W(J)$, where $J$ is the kernel of $R\limperf\to R$ and where $\hat W(J)$ is taken in the topological sense. 
This ring appears in \cite{Lau:Gerbe}.
The sheaf $\sW$ is derived $p$-complete, and its classical $p$-completion gives back $W$.
The Frobenius $F$ of $W(R)$ and a modified Verschiebung\footnote{For $p\ge 3$ we take $\tilde V=V$. For $p=2$ we define $\tilde V(x)=V(u_0x)$ where $V(u_0)=p-[p]$.}
$\tilde V$ carry over to $\sW(R)$,
and there is an exact sequence 
\beqn
0\to\sW(R)\xrightarrow{\tilde V}\sW(R)\xrightarrow\pi R.
\eeqn
For a nilpotent ideal $N\subseteq R$ there is another exact sequence
\beqn
\label{Eq:Intro-hWN-sWR-sWRN}
0\to\hat W(N)\to\sW(R)\xrightarrow\pi\sW(R/N).
\eeqn
In both cases $\pi$ is surjective if $R_{\red}$ is perfect.
For details we refer to \S\ref{Se:Sheard-Witt-vect}.
The analogue of \eqref{Eq:Intro-hWN-sWR-sWRN} for $W(R)$ involves $W(N)$ instead of $\hat W(N)$, and this difference allows the theory of sheared displays to work for all $p$-divisible groups, rather than just for infinitesimal $p$-divisible groups.

\subsection{Sheared displays}
\label{Se:Intro-sheared-disp}

We define plain sheared displays by substituting $\sW(R)$ for $W(R)$ and $\tilde V$ for $V$ in the classical definition of displays. Thus a plain sheared display is a collection 
\[
\uM=(M_0,M_1,\phi_0,\phi_1)
\] 
where $M_1\to M_0$ is a homomorphism of $\sW(R)$-modules that is a direct summand of a sum of copies of $\id_{\sW(R)}$ and $\tilde V$,
and where $\phi_i\colon M_i\to M_0$ is an $F$-linear map such that $\phi_1(M_1)$ generates $M_0$, and $\phi_1(\tilde V(a)x)=a\phi_0(x)$ for $a\in\sW(R)$ and $x\in M_0$.
These objects form a prestack, and the category of sheared displays is 
defined to be the associated stack for the topology generated by extractions of $p^\infty$-roots, or equivalently the associated fpqc stack.
See Definition \ref{De:Sh-disp} and Remark \ref{Re:sDisp-fpqc-stack}.
If $R_{\red}$ is perfect one can show that all sheared displays over $R$ are plain, using results of \cite{BMVZ}. 
See Corollary \ref{Co:sheared-disp-over-good-rings-are-plain}.
There is a functor 
\[
\sDisp(R)\to\Disp(R), \qquad 
\uM\mapsto\uM^\wedge
\] 
from sheared displays to displays given by classical $p$-completion. 

\subsection{Sheared functors}

The equivalence between $\BT(R)$ and the category $\sDisp(R)$ of sheared displays admits different descriptions depending on the direction of the functor.

\begin{ThmABC}
\label{Th:Intro-sFZR}
There is a functor\/ $\sFZ_R\colon \sDisp(R)\to\BT(R)$ which is determined by the following property. For $G=\sFZ_R(\uM)$ there is a
natural exact sequence of fpqc sheaves on $\Spec R$
\beqn
\label{Eq:TpG-M1-M0}
0\longrightarrow T_pG\longrightarrow\tilde M_1\xrightarrow{\;\phi_1-{\id}\;}\tilde M_0\longrightarrow 0
\eeqn
where the sheaves $\tilde M_i$ are defined by base change of sheared displays. The functor $\sFZ_R$ is an equivalence of exact categories and compatible with duality. There is a natural isomorphism of formal groups $\sFZ_R(\uM)^\wedge\cong\FZ_R(\uM^\wedge)$, so $\sFZ_R$ can be viewed as a decompletion of $\FZ_R$.
\end{ThmABC}

Here $T_pG=\lim G[p^n]$ is the Tate module scheme of $G\in\BT$, and $G^\wedge$ denotes the formal completion of $G$.
Theorem \ref{Th:Intro-sFZR} is proved in \S\ref{Se:Dieud-equiv:conclusions}.

\begin{ThmABC}
\label{Th:Intro-sDR}
There is an additive functor $\sD\colon\BT(R)\to\sDisp(R)$ compatible with base change in $R$ such that the composition with $\sDisp(R)\to\Disp(R)$ restricts to the equivalence $D_{\nil,R}\colon\BT_{\inf}(R)\to\Disp_{\Vnil}(R)$, and this uniquely determines the functor $\sD_R$. The functor $\sD_R$ is an equivalence of exact categories and compatible with duality.
\end{ThmABC}

Theorem \ref{Th:Intro-sDR} is proved in \S\ref{Se:Dieud-equiv:rigidity}.
The functors
$\sFZ_R$ and $\sD_R$ are mutually inverse equivalences $\BT(R)\cong\sDisp(R)$. For a special class of rings, including Artinian local rings with perfect residue field, sheared displays give back the Dieudonn\'e displays of \cite{Zink:Dieudonne, Lau:Relations}, and we recover the equivalence between $p$-divisible groups and Dieudonn\'e displays of loc.\,cit.

\subsection{Outline}

The proofs of Theorems \ref{Th:Intro-sFZR} and \ref{Th:Intro-sDR} are intertwined. 

\begin{enumerate}
\item
\label{It:construct-sFZ}
One proves directly that the homomorphism $\phi_1-\id$ in \eqref{Eq:TpG-M1-M0} is surjective and that its kernel is a Tate module scheme. Hence the functor $\sFZ_R$ is well-defined. See \S\ref{Se:Tate}.
\item
\label{It:construct-sDR}
The crystalline Dieudonn\'e functor produces the functor $\sD_R$ for $\FF_p$-algebras. See \S\ref{Se:sh-disp-functor}.
\item
\label{It:sDR-equiv-Fp}
The resulting functor $\sD\colon\BT\to\sDisp$ of fibered categories over $\Aff_{\FF_p}$ is formally \'etale with respect to nil-immersions of semiperfect $\FF_p$-algebras dominated by a power of Frobenius, and $\sD_R$ an equivalence over perfect $\FF_p$-algebras. 
These properties allow to deduce that $\sD_R$ is an equivalence for all $\FF_p$-algebras. See \S\S\ref{Se:Dieud-equiv:limits}--\ref{Se:Dieud.equiv:deform}.
\item
\label{It:sFZR-equivalence-char-p}
Using this, the relation $T_pG=\u\Hom(\QQ_p/\ZZ_p,G)$ yields that $\sFZ_R\circ\sD_R$ is the identity for $\FF_p$-algebras, hence $\sFZ_R$ is an equivalence for $\FF_p$-algebras. See \S\ref{Se:Dieud-equiv:conclusions}.
\item
\label{It:Rim-Schlessinger-sFZ}
A Rim--Schlessinger condition implies that $\sFZ_R$ is an equivalence in general; see \S\ref{Se:Dieud-equiv:conclusions}.
Define $\sD_R$ to be the inverse of $\sFZ_R$. 
\item
\label{It:Uniqueness-rigidity}
The functor $\sD_R$ is unique by a rigidity property of the fibered category $\BT$. See \S\ref{Se:Dieud-equiv:rigidity}.
\end{enumerate}
For $p\ge 3$, points \eqref{It:construct-sDR}--\eqref{It:sFZR-equivalence-char-p} extend to all $p$-nilpotent rings, so \eqref{It:Rim-Schlessinger-sFZ} 
is necessary only when $p=2$.
Let us now give more details on \eqref{It:construct-sFZ} and \eqref{It:construct-sDR}, corresponding to \S\ref{Se:Tate} and \S\ref{Se:sh-disp-functor}.

\subsection{Details on the Tate module functor}

For $\uM\in\sDisp(R)$ let $\sC(\uM)=[\tilde M_1\xrightarrow{\phi_1-{\id}}\tilde M_0]$ as a complex of fpqc sheaves of abelian groups over $\Spec R$ in degrees $[0,1]$ as in \eqref{Eq:TpG-M1-M0}. 
The existence of the functor $\sFZ_R$ in Theorem \ref{Th:Intro-sFZR} can be expressed as follows.

\begin{PrABC}
\label{Pr:Intro-sFZ}
For $\uM\in\sDisp(R)$ the complex\/ $\sC(\uM)$ is quasi-isomorphic to a Tate module scheme sitting in degree zero.
\end{PrABC}

See Theorem \ref{Th:sC-gives-BT}.
The proof of Proposition \ref{Pr:Intro-sFZ} proceeds along the following lines.
If $R$ is an $\FF_p$-algebra, there is a definition of $n$-truncated displays and a truncation functor $\uM\mapsto\uM_n$ from sheared displays to $n$-truncated displays, and we have an exact triangle 
\[
\FZ_{n,R}(\uM_n)\to\sC(\uM)/^{\LL}p^n\to\FY_{n,R}(\uM_n)\to \FZ_{n,R}(\uM_n)[1]
\]
where $\FZ_{n,R}(\uM_n)$ is an $n$-smooth group scheme, and $\FY_{n,R}(\uM_n)$ is an $n$-cosmooth group scheme.\footnote{A commutative group scheme $G$ in characteristic $p$ is $n$-smooth if locally $G=H[F^n]$ for a commutative formal group $H$, and $G$ is $n$-cosmooth if $G$ is the Cartier dual of an $n$-smooth group scheme.} The functor $\FZ_{n,R}$ is a truncated version of $\FZ_R$, 
and the functor $\FY_{n,R}$ is defined by a truncated version of $\sFZ_R$. 
The exact triangle gives Proposition \ref{Pr:Intro-sFZ} in characteristic $p$.
The general case follows by a deformation argument:  
Let $N\subseteq R$ be a nilpotent ideal. 
For $G\in\BT(R)$ there is an evident deformation sequence
\[
0\to T_pG(R)\to T_pG(R/N)\to \hat G(N)\to 0,
\]
and for $\uM\in\sDisp(R)$ and $T=H^0(\sC(\uM))$ one finds a similar deformation sequence
\[
0\to T(R)\to T(R/N)\to\FZ_R(\uM^\wedge)(R)\to\ldots
\]
as a consequence of \eqref{Eq:Intro-hWN-sWR-sWRN}.
A comparison of the deformation sequences allows us to extend Proposition \ref{Pr:Intro-sFZ} from $R/N$ to $R$. More precisely, assume that $T_{R/N}=T_p\bar G$ for $\bar G\in\BT(R/N)$. The ind-affine group scheme $V_p\bar G=T_p\bar G[p^{-1}]$ extends canonically to an ind-affine group scheme $V$ over $R$, and $T$ will be representable by a closed subscheme of $V$. 

The second deformation sequence also gives the comparison between $\sFZ_R$ and $\FZ_R$.

\subsection{Details on the sheared display functor}

Over $\FF_p$-algebras, the functor $\sD_R$ of Theorem \ref{Th:Intro-sDR} can be constructed stepwise as follows.
Using the terminology of frames and windows, the category of plain sheared displays can be written as $\Win(\usW(R))$, where $\usW(R)$ denotes a natural frame structure on $\sW(R)$.
If $R$ is a quasiregular semiperfect $\FF_p$-algebra, 
the ring $A_{\crys}(R)$ is $p$-torsion free and gives rise to a frame $\u A_{\crys}(R)$, moreover
there is a natural frame homomorphism $\srho\colon\u A_{\crys}(R)\to\usW(R)$, and the functor $\sD_R$ is the composition
\[
\BT(R)\xrightarrow{D_R^{\crys}}\Win(\u A_{\crys}(R))\xrightarrow{\;\srho\;}\sDisp(R)
\]
where $D_R^{\crys}$ is induced by the crystalline Dieudonn\'e functor.
By descent along the topology generated by extractions of $p^{\infty}$-roots, we get the functor $\sD_R$ for all quasisyntomic $\FF_p$-algebras. 
Since the algebraic stacks of truncated BT groups and their truncation morphisms are smooth, the functors $\sD_R$ for all $\FF_p$-algebras are obtained by Kan extension and ind-\'etale descent.

\subsubsection*{Perfect rings}

Over a perfect $\FF_p$-algebra $R$, sheared displays are equivalent to Dieudonn\'e modules, and the functor $\sD_R$ is known to be an equivalence.

\subsubsection*{Grothendieck--Messing theorem}

The functor $\sD\colon\BT\to\sDisp$ of fibered categories over $\Aff_{\FF_p}$ is formally \'etale because the Grothendieck--Messing theorem for $p$-divisible groups matches with a similar classification of deformations of sheared displays. This is a standard argument based on the following extension of the functors $\sD_R$.
For a nilpotent pd thickening of $\FF_p$-algebras $R\to R/N=\bar R$, one can define a frame $\usW(R/\bar R)$ with underlying ring $\sW(R)$, giving a category $\sDisp(R/\bar R)$  of relative sheared displays. This construction yields a Grothendieck--Messing theorem for sheared displays, which technically relies on the fact that the endomorphism $\sigma_1$ of $\hat W(N)$ is locally nilpotent.\footnote{The corresponding deformation result for displays is restricted to $F$-nilpotent or $V$-nilpotent objects because the endomorphism $\sigma_1$ of $W(N)$ is not locally nilpotent. This goes in hand with the fact that displays classify formal or unipotent $p$-divisible groups, but not both at the same time.}
The classification of deformations on both sides of $\sD$ is related via a functor $\sD_{R/\bar R}\colon\BT(\bar R)\to\sDisp(R/\bar R)$ with a crystalline comparison.

\subsubsection*{Odd primes}

If $p\ge 3$, a similar construction gives functors $\sD_R$ for all $p$-nilpotent rings, including the relative version and the Grothendieck--Messing classification.

\subsubsection*{Prismatic Dieudonn\'e theory}

For $p=2$, a similar construction with the prismatic Dieudonn\'e theory of \cite{Anschuetz-LeBras} as a starting point could be used to construct functors $\sD_R$ for all $p$-nilpotent rings. This will not be pursued here.

\subsubsection*{Comparison with the display functor}

One can expect that $\BT(R)\xrightarrow{\sD_R}\sDisp(R)\to\Disp(R)$ coincides with the known functor $D_R$; this is also part of the conjectures of \cite{Drinfeld:Talk-Manin}. The present methods give this comparison for infinitesimal $p$-divisible groups, by duality also for unipotent $p$-divisible groups, or when $R$ is an $\FF_p$-algebra or $p\ge 3$. The general case will be a consequence of the prismatic construction of the functor $\sD_R$. It would be interesting to have a more direct argument based on the characterisation of $\sD_R$ provided by Theorem \ref{Th:Intro-sDR}.

\subsection{Relation with other versions of Dieudonn\'e theory}

Different versions of Dieudonn\'e theory have aimed to classify $p$-divisible groups by suitable linear algebraic data, with increasing sophistication allowing base schemes of increasing generality. Over $\FF_p$-schemes, the original formulation of crystalline Dieudonn\'e theory 
\cite{Grothendieck:ICM, Grothendieck:Montreal, Messing:Crystals, Mazur-Messing, BBM} 
allows a classification of $p$-divisible groups if the base satisfies certain smoothness or regularity assumptions \cite{Berthelot-Messing, Jong:Crystalline}
and yields a fully faithful functor if the base is excellent l.c.i.\ \cite{Jong-Messing} or regular semiperfect \cite{Scholze-Weinstein:Moduli}.
Adding a divided Frobenius extends the classification to general $\FF_p$-schemes with certain finiteness conditions related to the Frobenius \cite{Lau:Semiperf, Lau:Divided}. All these results extend formally to a $p$-adic base when $p$ is odd. Prismatic Dieudonn\'e theory as developed in \cite{Anschuetz-LeBras} provides a classification of $p$-divisible groups over a quasi-syntomic $p$-adic base ring for all $p$, including $p=2$;
the case of integral perfectoid base rings was treated earlier in \cite{Scholze-Weinstein:Berkeley}. 
The recent classification of $p$-divisible groups over a general $p$-adic base by vector bundle $F$-gauges \cite{Gardner-Madapusi} replaces the prismatic site by the prismatization of \cite{Drinfeld:Prismatization, Bhatt-Lurie:Prismatization} and implicitly incorporates a divided Frobenius by means of the Nygaard filtration. 

The theories of displays of \cite{Zink:The-display} and of sheared displays developed here deviate from this line of thought by replacing crystals with projective modules over the Witt vectors or sheared Witt vectors, which seems to be a lossy operation a priori. The main results mean that displays of formal $p$-divisible groups and sheared displays of general $p$-divisible groups preserve all information. This is proved here using classical methods of Dieudonn\'e theory. The results could also be deduced from the classification of \cite{Gardner-Madapusi} by a comparison of vector bundle $F$-gauges and sheared displays. Such comparison, in a more general `Shimurian' setting related to a pair $(G,\mu)$ of 
a smooth affine group scheme with a cocharacter, is the content of \cite[Conjecture D.8.4]{Drinfeld:The-Lau}.

The above historical sketch is far from complete.
For example, there is a considerable literature on explicit Dieudonn\'e theory over complete discrete valuation rings \cite{Fontaine:Groupes, Breuil:Groupes, Kisin:Crystalline, Kim:Classification, Lau:Relations, Liu:Correspondence}. 
An alternative construction of Dieudonn\'e crystals via cohomology of  classifying stacks appears in \cite{Mondal:DieudonneI, Mondal:DieudonneII}. 
Results for $p$-divisible groups often have a finite counterpart. For example, a
classification of finite locally free group schemes in terms of $F$-gauges is given in \cite{Madapusi-Mondal}, based on \cite{Gardner-Madapusi}.

\subsection{Acknowledgements}

We are grateful to Vladimir Drinfeld and Akhil Mathew for helpful discussions, and to Drinfeld and the authors of \cite{BMVZ} for generously sharing unpublished notes and drafts.
This work was 
funded by the Deutsche Forschungsgemeinschaft (DFG, German Research
Foundation) -- Project-ID 491392403 -- TRR 358.

\subsection{Notation}

Divided powers are compatible with the canonical divided powers of $p$.

\cclearpage
\section{Sheared Witt vectors}
\label{Se:Sheard-Witt-vect}

In this section, we recall some aspects of the theory of sheared Witt vectors, developed in \cite{Drinfeld:Ring-stacks-conjecturally, BMVZ}, with some minor additions for later use. 
We include proofs for completeness.
Let $R$ be a $p$-nilpotent ring.

\subsection{First flat cohomology of $\hat W$}

A reduced ring $S$ is called seminormal if for any $x,y\in S$ with $x^3=y^2$ there is an $s\in S$ with $x=s^2$ and $y=s^3$. Every perfect $\FF_p$-algebra is seminormal.

\bPr
\label{Pr:H^1-hatW-seminormal}
If $R_{\red}$ is seminormal then $H^1_{\fpqc}(R,\hat W)=0$.
\ePr

\bproof
The fpqc sheaf $\hat W$ is a direct summand of the sheaf $\Lambda$ defined by
\[
\Lambda(A)=\ker(A[T]^*\xrightarrow{T\mapsto 0} A^*)=\{1+a_0T+\ldots+a_nT^n\mid a_n\in A\text{ nilpotent}\}.
\]
We have to show that for every faithfully flat ring homomorphism $R\to R'$ the \v Cech cohomology group $\check H^1(R'/R,\Lambda)$ vanishes. This group classifies isomorphism classes of invertible $R[T]$-modules equipped with a trivialisation at $T=0$ which become trivial under base change by $R\to R'$. If $R_{\red}$ is seminormal then $\Pic(R[T])\cong\Pic(R)$ by \cite[Chap.~I, Theorem 3.11]{Weibel:K-Book}. It follows that $\check H^1 (R'/R,\Lambda)$ is trivial.
\eproof

\bCo
\label{Co:H^1-hatW-fpqc-syn}
For every $R$ we have
$
H^1_{\fpqc}(R,\hat W)=H^1_{\fppf}(R,\hat W)=H^1_{\syn}(R,\hat W).
$
In other words, each fpqc $\hat W$-torsor over $\Spec R$ is trivialised over a syntomic covering.
\eCo

\bproof
Let $R'=R[a^{1/p^\infty}|a\in R]$.
Then $R'/p$ is semiperfect, so $R'_{\red}$ is perfect. Proposition \ref{Pr:H^1-hatW-seminormal} implies that every $\hat W$-torsor $Z$ over $\Spec R$ becomes trivial over $\Spec R'$ and thus corresponds to an element of $H^1(R'/R,\hat W)$. Since $\hat W$ preserves filtered colimits of rings, $Z$ becomes trivial over $\Spec R''$ with $R''=R[a_1^{1/p^n}\!\!,\ldots,a_r^{1/p^n}]$ for some $a_1,\ldots,a_r\in R$ and $n\ge 0$.
\eproof

\bRe
Corollary \ref{Co:H^1-hatW-fpqc-syn} also extends to the topology in which coverings are finite successions of extractions of $p$th roots. This will not be used in the following.
\eRe

\subsection{Definitions and basic properties}

\bDe
Let $Q=W/\hat W$ as an fpqc quotient sheaf on the category of $p$-nilpotent rings, or equivalently as an fppf or syntomic quotient sheaf by Corollary \ref{Co:H^1-hatW-fpqc-syn}. 
Let $Q\limperf=\lim(Q,F)$, the limit perfection of $Q$ with respect to $F$, and let $\sW=W\times_QQ\limperf$. Then $Q\limperf$ and $\sW$ are fpqc sheaves of $\delta$-rings. 
\eDe

\bRe
\label{Re:hatW-sW-Qperf}
The definition gives an exact sequences of fpqc sheaves
\beqn
\label{Eq:hatW-sW-Qperf}
0\to\hat W\to\sW\to Q\limperf\to 0
\eeqn
which is exact on the level of syntomic sheaves by Corollary \ref{Co:H^1-hatW-fpqc-syn}, moreover $W(R)\to Q(R)$ and $\sW(R)\to Q\limperf(R)$ are surjective when $R_{\red}$ is seminormal by Proposition \ref{Pr:H^1-hatW-seminormal}.
\eRe

\bRe
\label{Re:TFQ-sW-W}
Since $F\colon W\to W$ is surjective as ind-fppf sheaves, the same holds for $F\colon Q\to Q$ and for the projection $Q\limperf\to Q$. Hence we have an exact sequence of fpqc sheaves
\beqn
\label{Eq:TFQ-sW-W}
0\to T_FQ\to\sW\to W\to 0
\eeqn
which is exact as ind-fppf sheaves, where $T_FQ=\lim_nQ[F^n]$ is the $F$-adic Tate module of $Q$. This is an ideal of square zero by the following lemma. 
\eRe

\bLe
\label{Le:TFQ-Qperf0}
Let $Q^0=W(\hat\GG_a)/\hat W$ as an ideal of $Q$, and $Q^{\perf,0}$ the $F$-perfection of $Q^0$ as an ideal of $Q\limperf$. Then $T_FQ\subseteq Q^{\perf,0}$ and $T_FQ\cdot Q^{\perf,0}=0$.
\eLe

\bproof
It suffices to show that $Q[F^n]\subseteq Q^0$, which is clear, and that $Q[F^n]\cdot Q^0=0$. But $V\colon Q^0\to Q^0$ is bijective, and for $a\in Q[F^n]$ and $b\in Q^0$ we have $aV^n(b)=V^n(F^n(a)b)=0$.
\eproof

Let us call an ideal $N\subseteq R$ uniform locally nilpotent if there is an $m$ with $x^m=0$ for all $x\in N$; this hold iff the image of $N$ in $R/p$ is annihilated by a power of the Frobenius.

\bPr
\label{Pr:WN-sWR-sWRN}
For a uniform locally nilpotent ideal $N\subseteq R$ there is an exact sequence
\beqn
\label{Eq:WN-sWR-sWRN}
0\to\hat W(N)\to\sW(R)\to\sW(R/N).
\eeqn
If $R_{\red}$ is seminormal, the homomorphism $\sW(R)\to\sW(R/N)$ is surjective.
\ePr

\bproof 
We can assume that $R_{\red}$ is seminormal by passing to an fpqc covering if necessary.
Proposition \ref{Pr:H^1-hatW-seminormal} gives $Q(S)=W(S)/\hat W(S)$ for $S\in\{R,R/N\}$. Hence $Q(R)\to Q(R/N)$ is surjective with kernel annihilated by a power of $F$, and thus $Q\limperf(R)\to Q\limperf(R/N)$ is bijective. Since \eqref{Eq:hatW-sW-Qperf} is exact on $R$-point and on $R/N$-points, the lemma follows.
\eproof

\bEx
If $R$ is perfect, then $\sW(R)=W(R)$ since $Q(R)=W(R)$ by Remark \ref{Re:hatW-sW-Qperf}.
\eEx

\bEx
\label{Ex:Artin}
Assume that $\Nil(R)$ is uniform locally nilpotent and $R_{\red}$ is perfect; this holds when $R$ is an Artinian local ring with perfect residue field. Then $\sW(R)\to W(R)$ is injective and identifies $\sW(R)$ with the ring $\hat W(R)$ of \cite{Zink:Dieudonne} (when $R$ is Artinian) and $\WW(R)$ of \cite{Lau:Relations}. Indeed, \eqref{Eq:WN-sWR-sWRN} reads $0\to\hat W(R)\to\sW(R)\to W(R_{\red})\to 0$. This sequence splits uniquely and thus $\sW(R)=\hat W(R)\oplus W(R_{\red})$ as a subring of $W(R)=W(\Nil(R))\oplus W(R_{\red})$.
\eEx

\subsection{Rings with perfect reduction}

\bPr
\label{Pr:Qperf-Rred-surj}
If $R_{\red}$ is perfect then $Q\limperf(R)\to R_{\red}$ is surjective.
\ePr

\bproof
Cf.\ \cite[Proposition 3.42]{BMVZ}.
We can assume that $R$ is an $\FF_p$-algebra since $Q\limperf(R)=Q\limperf(R/p)$ as in the proof of Proposition \ref{Pr:WN-sWR-sWRN}. For a given element $a_0\in R$ choose $a_i\in R$ with $a_{i+1}^{p}\equiv a_i$ mod $\Nil(R)$. Then $F([a_{i+1}])\equiv([a_i])+ y_i$ mod $\hat W(R)$ for some $y_i\in V^{i+1} W(R)$. Let $v_i=[a_i]+\sum_{j\ge 0}F^j(y_{i+j})$. Then $F(v_{i+1})\equiv v_i$ mod $\hat W(R)$, and $(v_i)_i$ defines an element of $Q\limperf(R)$ that maps to the image of $a_0$ in $R_{\red}$.
\eproof

\bCo
\label{Co:sWR-R-surj}
If $R_{\red}$ is perfect then $\sW(R)\to R$ is surjective.
\eCo

\bproof
Use that $0\to\hat W(R)\to\sW(R)\to Q\limperf(R)\to 0$ is exact by Remark \ref{Re:hatW-sW-Qperf}.
\eproof

\subsection{Countably syntomic topology} 

\bDe
A ring homomorphism $R\to S$ will be called countably syntomic if $S$ is the colimit of a sequence of  syntomic ring homomorphisms $R=S_0\to S_1\to S_2\to\ldots$ .
\eDe

\bRe
Every countably syntomic ring homomorphism is countably ind-syntomic in the obvious sense, but not conversely, because for a countably syntomic homomorphism the cotangent complex has projective amplitude in $[-1,0]$,\footnote{Indeed, we have $L_{S/R}=\colim_i L_i$ where $L_i=L_{S_i/R}\otimes^L_{S_i}S$ has projective amplitude in $[-1,0]$ since $R\to S_i$ is syntomic. For an $S$-module $M$ it follows that $\Hom(L_{S/R},M[j])$ is zero for $j\not\in[0,2]$ and equal to $\lim^1_i\Hom(L_i,M[1])$ for $j=2$. Since $S_i\to S_{i+1}$ is syntomic, $\cone(L_i\to L_{i+1})$ has projective amplitude in $[-1,0]$, hence the modules $\Hom(L_i,M[1])$ form a surjective system, and $\lim^1_i$ vanishes.}
which does not hold for countably ind-syntomic homomorphisms in general.\footnote{Let $R\to S$ be $k\to k[t^{1/p^\infty}]/(t^{>1})$ for a field $k$ of characteristic $p$. Then $L_{S/R}=M[1]$ with $M=(t^{>1})/(t^{>2})$. The $S$-module $M$ is not projective because projective implies free by Kaplansky, and $tM=0$ while $t\ne 0$ in $S$.}
\eRe

\bLe
\label{Le:N-syntomic-compos}
The family of countably syntomic ring homomorphisms is stable under base change and under composition.
\eLe

\bproof
Base change is clear. For composition we use that the category of syntomic algebras is limit preserving by \cite[
\href{https://stacks.math.columbia.edu/tag/05N9}{Lemma 05N9}, \href{https://stacks.math.columbia.edu/tag/0C33}{Lemma 0C33}]{Stacks}. 
\eproof

\bDe
\label{De:N-syntomic-topology}
A covering in the countably syntomic topology on $\Aff$ is a finite family of jointly surjective countably syntomic morphisms. This defines a site $\Aff_{\csyn}$ by Lemma \ref{Le:N-syntomic-compos}.
\eDe

\bLe
\label{Le:N-syntomic-extend}
Let $A\to A'$ be a surjective ring homomorphism. For each countably syntomic ring homomorphism $A'\to B'$ there is a countably syntomic ring homomorphism $A\to C$ such that $A'\to C'=C\otimes_AA'$ factors as $A'\to B'\xrightarrow u C'$ where $u$ is faithfully flat and countably ind-\'etale, in particular countably syntomic.
\eLe

\bproof
For syntomic and \'etale in place of countably syntomic and countably ind-\'etale the assertion holds by \cite[\href{https://stacks.math.columbia.edu/tag/00T0}{Lemma 00T0}]{Stacks}. The lemma follows by a recursive procedure.
\eproof

\bLe
\label{Le:N-syntomic-cocont}
Let $R\to R'=R/N$ be a surjective ring homomorphism with locally nilpotent kernel $N$. Then the base change functor $\Aff_R\to\Aff_{R'}$ is continuous and cocontinuous for the countably syntomic topology.
\eLe

\bproof
Clearly the functor is continuous. The functor is cocontinuous by Lemma \ref{Le:N-syntomic-extend}, using that $\Spec R'\to\Spec R$ is bijective.
\eproof

\subsection{Frobenius endomorphism}

\bLe
\label{Le:F-on-hatW}
The homomorphism $F\colon \hat W\to\hat W$ is surjective for the syntomic topology.
\eLe

\bproof
This follows from the proof of \cite[Proposition 2.29]{BMVZ}, using that $F\colon W_{n+1}\to W_n$ is a syntomic cover because $F$ is an fppf cover between smooth schemes. 
Namely, $\hat W(R)$ is generated by elements $x=V^m([a])$ with $a\in\Nil(R)$, and we have $x=[b]V^m(1)$ when $b^{p^n}=a$, which exists syntomic locally. Here $F^n([b])=0$ for some $n$, and syntomic locally we have $b=c^p$ and $V^m(1)=F(y)+V^n(z)$ for some $y,z\in W(R)$. Then $x=F([c]y)$.
\eproof

\bLe
\label{Le:hatWFm-sW-sW}
There is an exact sequence of fpqc sheaves
\beqn
0\to\hat W[F^n]\to\sW\xrightarrow{F^n}\sW\to 0
\eeqn
which is exact on the level of syntomic sheaves.
\eLe

\bproof
Use Lemma \ref{Le:F-on-hatW} and the exact sequence \eqref{Eq:hatW-sW-Qperf} where $F$ is bijective on $Q\limperf$.
\eproof

\bLe
\label{Le:TFhW-Wperf-sW}
There is an exact sequence of fpqc sheaves
\beqn
0\to T_F\hat W\to W\limperf\to\sW\to 0
\eeqn
which is exact on the level of countably syntomic sheaves, 
where 
$W\limperf=\lim(W,F)$ is the limit perfection and $T_F\hat W=\lim_n\hat W[F^n]$ is the $F$-adic Tate module.
\eLe

\bproof
Lemma \ref{Le:F-on-hatW} and the syntomic exact sequence $0\to\hat W\to W\to Q\to 0$ give a countably syntomic exact sequence 
$0\to\hat W\limperf\to W\limperf\to Q\limperf\to 0$,
which maps to \eqref{Eq:hatW-sW-Qperf}. If follows that $W\limperf\to\sW$ is quasi-isomorphic to $\hat W\limperf\to\hat W$, which is surjective with kernel $T_F\hat W$.
\eproof

We have the following refinement of Proposition \ref{Pr:WN-sWR-sWRN}.

\bCo
\label{Co:sWR-sWR/N-Nsynt}
Let $\bar R=R/N$ where $N$ is uniform locally nilpotent and $i\colon\Spec\bar R\to\Spec R$ the canonical morphism. Let $\sW_R$ denote the restriction of the fpqc sheaf $\sW$ to $\Aff_R$. Then the natural homomorphism of countably syntomic sheaves $\sW_R\to i_*(\sW_{\bar R})$ is surjective.
\eCo

\bproof
The homomorphism $W_{R}\limperf\to i_*(W_{\bar R}\limperf)$ is bijective, and $i_*(W_{\bar R}\limperf)\to i_*(\sW_{\bar R})$ is surjective by Lemmas \ref{Le:TFhW-Wperf-sW} and \ref{Le:N-syntomic-cocont}.
\eproof

\subsection{Modified Verschiebung}
\label{Se:modif-V}

Let $\sW(\ZZ_p)=\lim_n\sW(\ZZ/p^n)$ and $\hat W(\ZZ_p)=\lim_n\hat W(\ZZ/p^n)$. Example \ref{Ex:Artin} gives $\sW(\ZZ_p)=\ZZ_p\oplus\hat W(\ZZ_p)$.
There is a unique element $\bar u_0\in W(\ZZ_p)/\hat W(\ZZ_p)$ with $V(\bar u_0)=p$, and $\bar u_0$ is a unit. Indeed, $V$ is injective on this quotient, and the residue class of $V^{-1}(p-[p])$ is one solution. 

\bLe 
\label{Le:u0}
The following are equivalent. 
\begin{enumerate}
\item
\label{It:u0-1}
$\bar u_0=1$, 
\item 
\label{It:u0-2}
the element $V(1)$ of $W(\ZZ_p)$ lies in $\sW(\ZZ_p)$,
\item
\label{It:u0-3}
$\sW(\ZZ_p)$ is stable under $V\colon W(\ZZ_p)\to W(\ZZ_p)$,
\item
\label{It:u0-4}
$p\ge 3$.
\end{enumerate}
\eLe

\bproof
$\eqref{It:u0-1}\Rightarrow\eqref{It:u0-2}\Leftrightarrow\eqref{It:u0-3}$ is straightforward, and $\eqref{It:u0-3}\Leftrightarrow\eqref{It:u0-4}\Rightarrow\eqref{It:u0-1}$ follows from \cite[Lemma 2]{Zink:Dieudonne} and its proof; in the notation of loc.\,cit.\ we have $\delta(^V1)=p$ and $^V\delta(1)=V(1)$.
\eproof

Let $u_0\in W(\ZZ_p)$ an inverse image of $\bar u_0$. The subsequent constructions do not depend on the choice of $u_0$. For definiteness we choose $V(u_0)=p-[p]$ for $p=2$ and $u_0=1$ for $p\ge 3$. Define $\tilde V\colon W\to W$ by $\tilde V(x)=V(u_0x)$. Then $\tilde VF\equiv p\equiv F\tilde V$ where $\equiv$ means equal mod $\hat W$, so the induced endomorphism $\tilde V$ of $Q$ satisfies $\tilde VF=F\tilde V=p$, and $\tilde V$ extends to $Q\limperf$ and to $\sW$. 
Let $\tilde p=F\tilde V(1)=pu_0$ in $W(\ZZ_p)$, thus $\tilde p=p-[p^2]$ for $p=2$ and $\tilde p=p$ for $p\ge 3$. The elements $V(u_0)$ and $\tilde p$ of $W(\ZZ_p)$ lie in $\sW(\ZZ_p)$ and hence give elements of $\sW(R)$ for every $R$.
In $\sW(R)$ we have $x\tilde V(y)=\tilde V(F(x)y)$ and thus $\tilde V F=V(u_0)\cdot\id$, moreover $F\tilde V=\tilde p\cdot\id$.

\bRe
Let $\bar\GG_a=\GG_a/\hat\GG_a$ as an fppf quotient. Then $\bar\GG_a(R)=(R_{\red})\colimperf=(R/p)\colimperf$, the colimit perfection, and $\bar\GG_a$ is also the fpqc or syntomic quotient $\GG_a/\hat\GG_a$. Indeed, $\FFF(R)=(R/p)\colimperf$ defines an fpqc sheaf such that $\GG_a\to\FFF$ is syntomic surjective with kernel $\hat\GG_a$.
\eRe

\bPr
\label{Pr:Qperf-pn-Qperf-WnbGa}
There is an exact sequence of fpqc sheaves
\beqn
\label{Eq:Qperf-pn-Qperf-WnbGa}
0\to Q\limperf\xrightarrow {p^n}Q\limperf\to W_n(\bar\GG_a)\to 0
\eeqn
which is exact on the level of countably syntomic sheaves; see Definition \ref{De:N-syntomic-topology}. The same holds with $\tilde V^n$ in place of $p^n$. If $R_{\red}$ is perfect, then \eqref{Eq:Qperf-pn-Qperf-WnbGa} is exact on $R$-points.
\ePr

\bproof
The action of $V^n$ on the syntomic exact sequence  $0\to\hat W\to W\to Q\to 0$ shows that $0\to Q\xrightarrow {V^n}Q\to W_n(\bar\GG_a)\to 0$ is syntomic exact, and the same with $\tilde V$ in place of $V$. Hence $0\to Q\limperf\xrightarrow{\tilde V^n}Q\limperf\to W_n(\bar\GG_a)\to 0$ is fpqc exact, and isomorphic to \eqref{Eq:Qperf-pn-Qperf-WnbGa} since $F$ is bijective on $Q\limperf$ with $p=F\tilde V=\tilde VF$. 
If $R_{\red}$ is perfect, Proposition \ref{Pr:Qperf-Rred-surj} implies that $Q\limperf(R)\to W_n(R_{\red})$ is surjective. Since every element of $W_n(\bar\GG_a)(R)$ is defined over a $\ZZ/p^r$-algebra $R'\subseteq R$ of finite type, and $R'$ admits a countably syntomic cover with perfect reduction, it follows that $Q\limperf\to W_n(\bar\GG_a)$ is countably syntomic surjective.
\eproof

\bCo
\label{Co:sW-Vn-sW-Wn}
There is an exact sequence of fpqc sheaves
\beqn
\label{Eq:sW-Vn-sW-Wn}
0 \to \sW \xrightarrow{\tilde V^n} \sW \to W_n \to 0
\eeqn
which is exact as countably syntomic sheaves, and exact on $R$-points if $R_{\red}$ is perfect.
\eCo

\bproof
This follows from Proposition \ref{Pr:Qperf-pn-Qperf-WnbGa} using the action of $\tilde V$ on the syntomic exact sequence $0\to\hat W\to\sW\to Q\limperf\to 0$, which is exact on $R$-points if $R_{\red}$ is perfect; see Remark \ref{Re:hatW-sW-Qperf}.
\eproof

\subsection{Semiperfect $\FF_p$-algebras}

\bLe
\label{Le:H1-hW-semiperf}
If $R$ is a semiperfect $\FF_p$-algebra then $H^1_{\fpqc}(R,\hat W[{F^n}])=0$.
\eLe

\bproof
The exact sequence $0\to\hat W[{F^n}]\to \hat W\xrightarrow{F^n}\hat W\to 0$ of fpqc sheaves (even of syntomic sheaves) given by Lemma \ref{Le:F-on-hatW} is exact on $R$-points since $R$ is semiperfect, and $H^1_{\fpqc}(R,\hat W)=0$ by Corollary \ref{Co:H^1-hatW-fpqc-syn} since $R_{\red}$ is perfect.
\eproof

\bPr
\label{Pr:sW-semiperf}
Let $R$ be a semiperfect $\FF_p$-algebra and let $J$ be the kernel of $R\limperf\to R$. Then $\sW(R)=W(R\limperf)/\hat W(J)$ where $\hat W(J)=\lim_n\hat W(J/F^n(J))$, i.e.\ $J$ is treated as an ideal of the topological ring $R\limperf$.
\ePr

\bproof
Since $F$ is surjective on $\hat W$ by Lemma \ref{Le:F-on-hatW}, the exact sequence $0\to\hat W\to W\to Q\to 0$ implies that 
$\sW$ is the limit of $\sW^{(n)}=W/\hat W[F^n]$ with transition maps induced by $F$. Lemma \ref{Le:H1-hW-semiperf} gives exact sequences $0\to\hat W(R)[F^n]\to W(R)\to\sW^{(n)}(R)\to 0$ with surjective transition maps given by $F$ on $W(R)$. The limit over $n$ gives the result.
\eproof

\bCo
\label{Co:sW-p-complete}
For every $p$-nilpotent ring $R$ the ring $\sW(R)$ is derived $(p,\tilde p)$-complete.
\eCo

\bproof
Since $\sW$ is an fpqc sheaf we can assume that $R/p$ is semiperfect. Since $\hat W(pR)$ is annihilated by a power of $p$ and thus of $\tilde p$, Proposition \ref{Pr:WN-sWR-sWRN} reduces the assertion to the case where $R$ is semiperfect. This case follows from Proposition \ref{Pr:sW-semiperf}. 
\eproof

\bRe
If $R$ is a quasiregular semiperfect $\FF_p$-algebra, then  the ring $A_{\crys}(R)$ is $p$-torsion-free. Let $N$ be the kernel of the natural homomorphism $A_{\crys}(R)\to W(R)$ and $N^{\nil}\subseteq N$ the submodule of all elements where $\dot\varphi=\varphi/p$ is topologically nilpotent for the $p$-adic topology of $N$. 
Then $A_{\crys}(R)/N^{\nil}\cong\sW(R)$ by \cite{Lau:Gerbe}. This is not used in the following.
\eRe

\subsection{Divided powers}
\label{Se:pd-sW}

\bLe
\label{Le:pd-sW}
Assume that $R$ is an $\FF_p$-algebra, or $p\ge 3$.
Then the kernel of the homomorphism $\sW(R)\to R$ carries natural divided powers.
\eLe

\bproof
The assumption implies that we can take $u_0=1$ in $W(R)$ and thus $\tilde V=V$. By Corollary \ref{Co:sW-Vn-sW-Wn} the kernel of  $\sW(R)\to R$ is equal to $\tilde V(\sW(R))$. We define divided powers on this ideal by the formula $(p-1)!\gamma_p(\tilde V(x))=p^{p-2}\tilde V(x^p)$. The relevant axioms can be verified using that $\sW=W\times_QQ\limperf$: The formula gives the usual divided powers on $VW(R)$, moreover in $Q(R)\limperf$ we have $p=VF=FV$ where $F$ is bijective, and the formula gives the canonical divided powers on $pQ(R)\limperf=VQ(R)\limperf$.
\eproof

\bLe
\label{Le:pd-sW-rel}
Let $R'\to R$ be a nilpotent pd thickening with kernel $\Fa$. Assume that $R'$ is an $\FF_p$-algebra, or $p\ge 3$.
Then the kernel of the natural homomorphism $\sW(R')\to R$ carries natural divided powers.
\eLe

\bproof
Again, the assumption implies that we can take $u_0=1$ in $W(R')$ and thus $\tilde V=V$. 
Proposition \ref{Pr:WN-sWR-sWRN} and Corollary \ref{Co:sW-Vn-sW-Wn} give a commutative diagram with exact rows and columns where the vertical maps $\tilde V$ are injective and $\hat W(\Fa)\to\Fa$ is surjective. (If $R_{\red}$ is perfect, then all rows and columns are short exact sequences.)
\[
\xymatrix@M+0.2em{
0 \ar[r] & \hat W(\Fa) \ar[r] \ar@{^(->}[d]^{\tilde V} & \sW(R') \ar[r] \ar@{^(->}[d]^{\tilde V} & \sW(R) \ar@{^(->}[d]^{\tilde V} \\
0 \ar[r] & \hat W(\Fa) \ar[r] \ar@{->>}[d] & \sW(R') \ar[r] \ar[d] & \sW(R) \ar[d] \\
0 \ar[r] & \Fa \ar[r] & R' \ar[r] & R \ar[r] & 0
} 
\]
The nilpotent divided powers on $\Fa$ give an isomorphism $\hat w'\colon \hat W(\Fa)\cong\bigoplus_{\NN}\Fa$; see \cite[(149)]{Zink:The-display}. Define $i\colon \Fa\to\hat W(\Fa)$ by $\hat w'(i(a))=(a,0,0,\ldots)$. Then $i$ is an injective homomorphism of $W(R')$-modules, and $\hat W(\Fa)=\tilde V\hat W(\Fa)\oplus i(\Fa)$. The diagram implies that the kernel of $\sW(R')\to R$ is equal to $\tilde V\sW(R')\oplus i(\Fa)$. This ideal carries divided powers induced by those on $\tilde V\sW(R')$ given by Lemma \ref{Le:pd-sW} and those on $i(\Fa)$ induced by the divided powers on $\Fa$.
\eproof

\cclearpage
\section{Descent of finite projective $\sW$-modules}

This section is of preliminary nature. We provide an elementary approach to a weak version of the descent properties of finite projective $\sW(R)$-modules. See \S\ref{Se:Descent-II} for a stronger version.

\subsection{Vanishing of some cohomology groups}

Let $f\colon R\to R'$ be a faithfully flat ring homomorphism. For an abelian presheaf $F$ on the category of $R$-algebras we denote by 
$\check C^*(R'/R,F)^+$ the augmented \v Cech complex starting in degree $-1$,
\[
\check C^*(R'/R,F)^+=[ F(R)\to F(R')\to F(R'\otimes_RR')\to\ldots\ ].
\]

\bLe
\label{Le:H*W-fpqc}
The complex $C^*=\check C^*(R'/R,W)^+$ is acyclic, thus $\check H^i(R'/R,W)=0$ for $i\ge 1$.
\eLe

\bproof
Let $C^*_n=\check C^*(R'/R,W_n)^+$. Then $C_1^*$ is acyclic by fpqc descent, $C_n^*$ is acyclic by induction, and $C^*=\lim_n C^*_n$ is acyclic since the transition maps in the limit are surjective.
\eproof

\bLe
\label{Le:H*WFn-fpqc}
If $R$ and $R'$ are semiperfect $\FF_p$-algebras, then $\check C^*(R'/R,W[{F^n}])^+$ is acyclic.
\eLe

\bproof
For every semiperfect $\FF_p$-algebra $S$ we have $0\to W(S)[{F^n}]\to W(S)\xrightarrow{F^n} W(S)\to 0$. Since all tensor products $R'\otimes_R\ldots\otimes_RR'$ are semiperfect, Lemma \ref{Le:H*W-fpqc} gives the result.
\eproof

\bLe
\label{Le:H*hWN-fpqc}
If $R$ is $p$-nilpotent and $N$ is a locally nilpotent ideal of $R$ (or a locally nilpotent $R$-algebra), then $\check C^*(R'/R,\hat W(\tilde N))^+$ is acyclic; here $\tilde N(S)=N\otimes_RS$ for an $R$-algebra $S$. 
\eLe

\bproof
Since the functor $N\mapsto \check C^*(R'/R,\hat W(\tilde N))^+$ preserves filtered colimits we can assume that $N$ is finitely generated and thus nilpotent. Since the functor preserves exact sequences of nilpotent algebras we can assume that $N$ has square zero and $pN=0$. Then $\hat W(\tilde N)\cong\bigoplus_{\NN}\tilde N$, and the lemma follows since $\check C^*(R'/R,\tilde N)^+$ is acyclic by fpqc descent.
\eproof

\bLe
\label{Le:H*hW-indet}
If $f$ is ind-\'etale and $R$ is $p$-nilpotent, then $\check C^*(R'/R,\hat W)^+$ is acyclic.
\eLe

\bproof
Let $N=\Nil(R)$. For an ind-\'etale $R$-algebra $S$ we have $N\otimes_RS=\Nil(S)$ and hence $C(R'/R,\hat W)^+=C(R'/R,\hat W(\tilde N))^+$, which is acyclic by Lemma \ref{Le:H*hWN-fpqc}.
\eproof

\bLe
\label{Le:H*hWFn-indet}
If $f$ is ind-\'etale and $R$ is semiperfect, then $\check C^*(R'/R,\hat W[{F^n}])^+$ is acyclic.
\eLe

\bproof
For every semiperfect $\FF_p$-algebra $S$ we have $0\to\hat W(S)[{F^n}]\to \hat W(S)\xrightarrow{F^n} \hat W(S)\to 0$.
Since all tensor products $R'\otimes_R\ldots\otimes_RR'$ are semiperfect, Lemma \ref{Le:H*hW-indet} gives the result.
\eproof

\bLe
\label{Le:H^*K-indet}
If $f$ is ind-\'etale and $R$ is semiperfect, then $\check C^*(R'/R,T_FQ)^+$ is acyclic.
\eLe

The fpqc sheaf $T_FQ$ was defined in Remark \ref{Re:TFQ-sW-W}.

\bproof
If $S$ is semiperfect then $0\to\hat W(S)[{F^n}]\to W(S)[{F^n}]\to Q(S)[F^n]\to 0$ is exact by Lemma \ref{Le:H1-hW-semiperf}, and $F\colon Q(S)[F^{n+1}]\to Q(S)[F^n]$ is surjective. Hence $\check C^*(R'/R,Q[F^n])^+$ is acyclic by Lemmas \ref{Le:H*WFn-fpqc} and \ref{Le:H*hWFn-indet}, and the limit $\check C^*(R'/R,T_FQ)^+$ is acyclic as well.
\eproof

\subsection{Descent over semiperfect rings}

\bPr
\label{Pr:Descent-LF(sW)-semiperf}
The fibered category of finite projective $\sW(R)$-modules for semiperfect $\FF_p$-algebras $R$ is an ind-\'etale stack over the category $\Aff_{\FF_p}^{\semiperf}$ of semiperfect affine $\FF_p$-schemes.
\ePr

\bproof
Homomorphisms of finite projective $\sW(R)$-modules form fpqc sheaves since $\sW$ is an fpqc sheaf. For every semiperfect $\FF_p$-algebra $S$, \eqref{Eq:TFQ-sW-W} and Proposition \ref{Pr:sW-semiperf} give an exact sequence
\[
0\to T_FQ(S)\to \sW(S)\to W(S)\to 0
\] 
where $T_FQ(S)$ is an ideal of square zero by Lemma \ref{Le:TFQ-Qperf0}.
A descent datum of finite projective $\sW(-)$-modules for an ind-\'etale covering $R\to R'$ of semiperfect $\FF_p$-algebras induces an fpqc sheaf of $\sW$-modules $\sMMM$ over $\Spec R$, which sits in an exact sequence of fpqc sheaves of $\sW$-modules
\[
0\to\MMM\otimes_{W}T_FQ\to\sMMM\to\MMM\to 0
\]
with $\MMM=\sMMM\otimes_{\sW}W$, moreover $\sMMM(S)\to\MMM(S)$ is surjective for every semiperfect $R$-algebra $S$ which admits a homomorphism of $R$-algebras $R'\to S$. By fpqc descent of finite projective $W(-)$-modules \cite[Corollary 34]{Zink:The-display}, there is a finite projective $W(R)$-module $M$ such that $\MMM=M\otimes_{W(R)}W$. Let $\sM$ be a finite projective $\sW(R)$-module that lifts $M$. 
The homomorphism $\sMMM(R)\to\MMM(R)=M$ is surjective since $H^1(R'/R,\MMM\otimes_{W}T_FQ)$ vanishes by Lemma \ref{Le:H^*K-indet}. 
Hence the isomorphism $M\otimes_{W(R)}W\xrightarrow\sim\MMM$ lifts to a homomorphism $\sM\otimes_{\sW(R)}\sW\to\sMMM$, which is an isomorphism since this holds over $R'$. 
\eproof

\subsection{Descent over essentially semiperfect rings}

\bDe
\label{De:essemiperf}
A $p$-nilpotent ring $R$ will be called essentially semiperfect if for some $n\ge 1$ and $\bar R=R/p$ the ring  $\bar R/\bar R[F^n]$ is semiperfect.
\eDe

\bDe
\label{De:proot}
A homomorphism of $p$-nilpotent rings $R\to R'$ will be called a $p^\infty$-root morphism if it arises by pushout from  $\ZZ_p[\{t_i\}_{i\in I}]\to\ZZ_p[\{t_i^{1/p^\infty}\}_{i\in I}]$. Finite compositions of $p^\infty$-root morphisms form the coverings of a topology on $\Aff_{\Spf\ZZ_p}$ that will be called the $p^\infty$-root topology.
\eDe

\bRe
\label{Re:proot-semiperf}
If $R$ is semiperfect, every $p^\infty$-root morphism $\Spec R'\to\Spec R$ has a section.
\eRe

\bPr
\label{Pr:Descent-LF(sW)-essemiperf}
The fibered category of finite projective $\sW(R)$-modules for essentially semiperfect rings $R$ is an ind-\'etale and $p^\infty$-root stack over the category $\Aff^{\essemiperf}$ of essentially semiperfect  affine schemes.
\ePr

\bproof
Homomorphisms of finite projective $\sW(R)$-modules form fpqc sheaves. Let $R$ be given.
We have to show that descent of finite projective $\sW(-)$-modules is effective for every ind-\'etale or $p^\infty$-root covering $R\to R'$.
If $R$ is semiperfect, this holds by Proposition \ref{Pr:Descent-LF(sW)-semiperf} and Remark \ref{Re:proot-semiperf}. By induction we can assume that the assertion holds for $\bar R=R/N$ where $N$ is an ideal of square zero with $pN=0$. For an ind-\'etale or $p^\infty$-root morphism $R\to S$ (it suffices that $R\to S$ is flat and $S_{\red}$ is perfect), Proposition \ref{Pr:WN-sWR-sWRN} gives an exact sequence
\[
0\to\hat W(N\otimes_RS)\to\sW(S)\to\sW(S/NS)\to 0
\]
where $\hat W(N\otimes_RS)$ is an ideal of square zero. A decent datum of finite projective $\sW(-)$-modules for $R\to R'$ induces an fpqc sheaf of $\sW$-modules $\MMM$ over $\Spec R$, which sits in an exact sequence of fpqc sheaves of $\sW$-modules
\[
0\to\bar\MMM\otimes_{\bsW}\hat W(\u N)\to\MMM\to\bar\MMM\to 0
\]
where $\bsW(S)=\sW(S/NS)$ and $\bar\MMM=\MMM\otimes_{\sW}\bsW$ and $\u N(S)=NS$ for an $R$-algebra $S$, moreover the homomorphism $\MMM(S)\to\bar\MMM(S)$ is surjective when $R\to S$ factors over $R'$.
The assumption implies that $\bar\MMM=\bar M\otimes_{\sW(\bar R)}\bsW$ for a finite projective $\sW(\bar R)$-module $\bar M$. Let $M$ be a finite projective $\sW(R)$-module that lift $\bar M$. The homomorphism $\MMM(R)\to\bar\MMM(R)=\bar M$ is surjective since $H^1(R'/R,\bar\MMM\otimes_{\bsW}\hat W(\u N))$ vanishes by Lemma \ref{Le:H*hWN-fpqc}. Hence the isomorphism $\bar M\otimes_{\sW(\bar R)}\bsW\xrightarrow\sim\bar\MMM$ lifts to a homomorphism $M\otimes_{\sW(R)}\sW\to\MMM$, which is an isomorphism since this holds over $R'$. 
\eproof

\cclearpage
\section{Frames and windows}

This section provides a formalism that allows to treat displays, truncated displays, and sheared displays, starting with the category of arrows of abelian groups with a tensor structure, whose algebras are rings with a quasi-ideal in the sense of \cite{Drinfeld:Ring-groupoid}. This is the minimal framework necessary to study these variations of displays. We will follow Zink's terminology of frames and windows, introduced in \cite{Zink:Windows} under slightly different assumptions. One could also use the formalism of higher displays of \cite{Lau:Higher} to treat the relevant display categories, and this might be recommendable for various reasons. The present approach emphasizes the point that in the higher display expression of the theory, all relevant information is visible in degrees $\{0,1\}$.

\subsection{Filtered abelian groups, rings, and modules}

Let $\Ab^{\leftarrow}$ be the abelian category of homomorphisms of abelian groups $M=(M_0\xleftarrow{\:\tau} M_1)$, which will be called filtered abelian groups for simplicity (better would be $1$-filtered). The category $\Ab^{\leftarrow}$ carries the tensor product
\[
M\otimes N=(M_0\otimes N_0\leftarrow M_1\otimes N_0\oplus_{M_1\otimes N_1}M_0\otimes N_1).
\]

A commutative ring in $\Ab^{\leftarrow}$ will be called a filtered ring; this is equivalent to a commutative ring with a quasi-ideal in the sense of \cite{Drinfeld:Ring-groupoid}.
For a filtered ring $A=(A_0\leftarrow A_1)$ we denote by $\Mod(A)$ the category of $A$-modules in $\Ab^{\leftarrow}$ and call these objects filtered $A$-modules. 
This is an abelian category, the forgetful functor $\Mod(A)\to\Ab^\leftarrow$ is exact, and $\Mod(A)$ carries the tensor product
\[
M\otimes_AN=\coker(M\otimes A\otimes N\xrightarrow\psi M\otimes N)
\]
where $\psi(x\otimes a\otimes y)=ax\otimes y-x\otimes ay$.
For a homomorphism of filtered rings $u\colon A\to B$ there is a base change functor $\Mod(A)\to\Mod(B)$ defined by $M\mapsto M^u=M\otimes_AB$. 

Let $A$ be a filtered ring.
There are objects $A(i)\in\Mod(A)$ with $\Hom_A(A(i),M)=M_i$ for $i=0,1$, namely $A(0)=A$ and $A(1)=(A_0\leftarrow A_0)$. 
The functor $\Mod(A)\to\Mod(A_0)$, $M\mapsto M_i$ has a left adjoint $N\mapsto N\otimes_{A_0}A(i)$, where the tensor product is formed componentwise.
(Finite) direct sums of copies of $A(0)$ and $A(1)$ will be called (finite) free filtered $A$-modules. A filtered $A$-module will be called finite if it is a quotient of a finite free filtered module. 
Finite projective filtered modules are direct summands of  finite free filtered modules. 

A filtered ring $A$ will be called Henselian if $(A_0,\tau A_1)$ is a Henselian pair.

\subsection{Structure of filtered modules}

\bLe
\label{Le:fil-proj-mod}
Let $A$ be a Henselian filtered ring. A filtered $A$-module $M$ is finite projective iff $M$ is isomorphic to $M'=L_0\otimes_{A_0}A(0)\oplus L_1\otimes_{A_0}A(1)$ for finite projective $A_0$-modules $L_i$. The assignment $(L_0,L_1)\mapsto M'$ induces a bijection between pairs of isomorphism classes of finite projective $A_0$-modules and isomorphism classes of finite projective filtered $A$-modules.
\eLe

\bproof
Clearly the filtered $A$-module $M'$ is finite projective. Let $R_0=A_0/\tau A_1$ and let $R=(R_0\leftarrow 0)$ as a filtered ring. The hypothesis implies that $A_0\to R_0$ gives a bijection on isomorphism classes of finite projective modules. 

For a filtered $A$-module $M$ let $M_R=M\otimes_AR$ as a filtered  $R$-module, and for an $R$-module $N$ let $\gr(N)=(\coker(N_1\to N_0),N_1)$ as a $\{0,1\}$-graded $R_0$-module. One verifies that $\gr(M_R)_0=\coker(M_1\to M_0)$. The functors $M\mapsto M_R$ and $N\mapsto\gr(N)$ preserve finite free modules and hence finite projective modules, and we have $\gr(M')=(L_0\otimes_{A_0}R_0,L_1\otimes_{A_0}R_0)$. It follows that $(L_0,L_1)\mapsto M'$ is injective on isomorphism classes. 

To prove surjectivity, let $M$ be a finite projective filtered $A$-module. Let $\bar L_0=\gr(M_R)_0$, a finite projective $R_0$-module, and let $L_0$ be a lift of $\bar L_0$ to a finite projective $A_0$-module. Since $M$ is projective there is a homomorphism $\pi\colon M\to L_0\otimes_{A_0}A(0)$ that lifts the identity of $\bar L_0$, and $\pi$ is necessarily surjective. Hence $M=N\oplus L_0\otimes_{A_0}A(0)$ for a filtered  $A$-module $N$. Here $\gr(N_R)_0=0$, so the natural homomorphism $\psi\colon N_1\otimes_{A_0}A(1)\to N$ is surjective, thus $\psi$ splits since $N$ is projective, which implies that $\psi$ is bijective. Necessarily $N_1$ is projective over $A_0$, thus $M$ is of the type $M'$. 
\eproof

\bRe
\label{Re:fil-proj-mod-ideal}
Let $A$ be a Henselian filtered ring such that $A_1\xrightarrow\tau A_0$ is injective. A filtered $A$-module $M$ is finite projective iff $M_0$ is a finite projective $A_0$-module, $M_1\to M_0$ is injective, and $M_0/M_1$ is projective over $A_0/A_1$. This is immediate from Lemma \ref{Le:fil-proj-mod}.
\eRe

\bRe
\label{Re:graded-hom-deg-1}
Let $A$ be a filtered ring such that $A_1$ is an invertible $A_0$-module. For a filtered $A$-module $M$ and an $A_0$-module $N_0$ the restriction to degree $1$ gives an isomorphism
\[
\Hom_{A}(M,N_0\otimes_{A_0}A)\xrightarrow\sim\Hom_{A_0}(M_1,N_0\otimes_{A_0}A_1).
\]
Indeed, for a given homomorphism $g_1\colon M_1\to N_0\otimes_{A_0}A_1$ the composition with the multiplication $M_0\otimes_{A_0}A_1\to M_1$ gives $g_0\colon M_0\to N_0$ and thus $g\colon M\to N$.
\eRe

\subsection{Frames}

\bDe 
\label{De:frame}
A frame is a triple $\u A=(A,d,\sigma)$ where $A=(A_0\xleftarrow{\:\tau}  A_1)$ is a filtered ring, $d\in A_0$ is an element, which gives a filtered ring $A'=(A_0\xleftarrow{\:\text{\normalsize$_d$}} A_0)$, and $\sigma\colon A\to A'$ is a homomorphism of filtered rings such that $\sigma_0\colon A_0\to A_0$ induces the Frobenius on $A_0/p$. We require that $A_0$ and $A_1$ are derived $(p,d)$-complete. 

For another frame $\u B=(B,d',\sigma)$, a frame homomorphism 
$(g,u)\colon\u A\to\u B$ consists of a homomorphism of filtered rings $g\colon A\to B$ and a unit $u\in B_0^*$ with $g_0(d)=ud'$, which implies that 
$g'=(g_0,ug_0)\colon A'\to B'$ is a homomorphism of filtered rings, 
such that $\sigma\circ g=g'\circ\sigma$.
The frame homomorphism is called strict if $u=1$.
\eDe

\bRe
In Definition \ref{De:frame},
the homomorphism of filtered rings $\sigma:A\to A'$ consists of a Frobenius lift $\sigma_0:A_0\to A_0$ and a $\sigma$-linear map $\sigma_1:A_1\to A_0$ such that $\sigma\circ\tau=d\sigma_1$. For a frame homomorphism $(g,u):\u A\to\u B$ we have $g_0\sigma=\sigma g_0$ and $u\cdot g_0\sigma_1=\sigma_1g_1$.
\eRe

\bLe
\label{Le:frame-Henselian}
For a frame $\u A$ as in Definition \ref{De:frame}, the filtered ring $A$ is Henselian.
\eLe

\bproof
The ring $A_0$ is Henselian along $(p,d)$ by the completeness assumption, and the image of $\tau A_1$ in $A_0/(p,d)$ is annihilated by the Frobenius.
\eproof

\bDe 
A prismatic frame is a frame $\u A$ where $A_0$ is a $\delta$-ring such that $\sigma_0\colon A_0\to A_0$ is the associated Frobenius lift and $d$ is distinguished, i.e.\ every homomorphism of $\delta$-rings $A_0\to W(k)$ for a perfect field $k$ sends $d$ to a unit multiple of $p$.
\eDe

\bDe 
\label{De:crys-frame}
A crystalline frame is a prismatic frame with $d=p$ where $\tau:A_1\to A_0$ is injective and the image of this map carries divided powers $\gamma$ such that 
\beqn
\label{Eq:pd-frame-1}
\sigma_1(a)=\tilde\gamma_p(a)+\delta(a),
\eeqn
\beqn
\label{Eq:pd-frame-2}
\sigma_1(\tilde\gamma_p(a))=p^{p-2}\sigma_1(a)^p
\eeqn 
for $a$ in $A_1$ viewed as an ideal of $A_0$, with $\tilde\gamma_p=(p-1)!\gamma_p$. We do not require $pA_0\subseteq A_1$.
\eDe

\bRe
Crystalline frames in the sense of Definition \ref{De:crys-frame} are related to the $\delta$-Cartier rings of \cite[Definition 3.4.4]{Magidson:Witt}. More precisely, crystalline frames $\u A$ such that $\sigma_1\colon A_1\to A_0$ is bijective are equivalent to derived $p$-complete $\delta$-Cartier rings $(B,F,V,\delta)$ such that $V$ is injective by sending $\u A$ to $(A_0,\sigma_0,\tau\circ\sigma_1^{-1},\delta)$.
\eRe

\bRe
\label{Re:pd-frame}
Assume that $(A_0,A_1,\gamma)$ is a pd pair where $A_0$ is a $\delta$-ring. Let $\sigma_0\colon A_0\to A_0$ be the associated Frobenius lift. Then \eqref{Eq:pd-frame-1} defines a $\sigma_0$-linear map $\sigma_1\colon A_1\to A_0$ with $p\sigma_1=\sigma_0$. The $\delta$-structure is equivalent to a lift $s\colon A_0\to W_2(A_0)$ of the projection via $s(a)=(a,\delta(a))$. The kernel of $W_2(A_0)\to A_0/A_1$ carries natural divided powers $\gamma'$ by \cite[2.3]{Zink:The-display}, and \eqref{Eq:pd-frame-2} is equivalent to $s$ being a pd homomorphism. This is a direct consequence of \eqref{Eq:pd-frame-1} and the formula $\tilde\gamma_p'(a,b)=(\tilde\gamma_p(a),p^{p-2}(\tilde\gamma_p(a)+b)^p-\tilde\gamma_p(\tilde\gamma_p(a)))$,  which is immediate from loc.\,cit.
\eRe

\subsection{Modules and windows}

Let $\u A=(A,d,\sigma)$ be a frame and $A'=(A_0\xleftarrow {\:d} A_0)$
as above.

\bDe
\label{De:uA-mod}
An $\u A$-module is a pair $(M,\phi)$ where $M=(M_0\leftarrow M_1)$ is a filtered $A$-module and $\phi\colon M\to M_0\otimes_{A_0}A'$ is a $\sigma$-linear homomorphism of filtered modules.
\eDe

\bRe
The $\sigma$-linear map $\phi$ is equivalent to an $A'$-linear map $\phi^\sharp\colon M^\sigma\to M_0\otimes_{A_0}A'$, which is equivalent to its degree $1$ component $\phi^\sharp_1\colon(M^\sigma)_1\to M_0$ by Remark \ref{Re:graded-hom-deg-1}.
\eRe

\bDe
\label{De:window}
An $\u A$-window is an $\u A$-module $(M,\phi)$ where the filtered $A$-module $M$ is finite projective and the $A_0$-linear map $\phi^\sharp_1$ is bijective. The height $\height(M,\phi)$ is the rank of $M_0$ as an $A_0$-module. Let $R=A_0/\tau A_1$. The dimension $\dim(M,\phi)$ is the rank of $M_0/\tau M_1$ as an $R$-module.
\eDe

\bRe
\label{Re:tensor-Mod-frame}
The category $\Mod(\u A)$ of $\u A$-modules is abelian and carries a tensor product defined by $(M,\phi_M)\otimes(N,\phi_N)=(M\otimes_AN,\phi_M\otimes\phi_N)$. The category $\Win(\u A)$ of $\u A$-windows is an exact category which is not stable under the tensor product. 
\eRe

\bRe[Functoriality]
A frame homomorphism $g\colon \u A\to\u B$ induces a base change functor
\[
g^*\colon \Mod(\u A)\to\Mod(\u B),
\qquad
g^*(M,\phi)=(M\otimes_AB,\phi\otimes \sigma),
\]
that restricts to a functor $g^*\colon \Win(\u A)\to\Win(\u B)$.
Here we use the identification 
\[
(M\otimes_AB)_0\otimes_{B_0}B'=(M_0\otimes_{A_0}A')\otimes_{A'}B'.
\]
We will also write $g^*\uM=\uM\otimes_{\u A}\u B$. 
\eRe

\bRe[Normal representations]
\label{Re:normal-repres}
The category $\Win(\u A)$ is equivalent to the following category $\Win'(\u A)$. Objects are triples $(L_0,L_1,\Psi)$ where the $L_i$ are finite projective $A_0$-modules and 
\[
\Psi\colon L_0^{\sigma_0}\oplus L_1^{\sigma_0}\xrightarrow\sim
L_0\oplus L_1
\]
is an isomorphism of $A_0$-modules. A homomorphism $(L_0,L_1,\Psi)\to(L_0',L_1',\Psi')$ consists of homomorphisms of $A_0$-modules $a\colon L_0\to L_0'$, $b\colon L_1\to L_0'\otimes_{A_0}A_1$, $c\colon L_0\to L_1'$, $e\colon L_1\to L_1'$ 
such that
\[
\left(\begin{matrix}a&\tau(b)\\c&e\end{matrix}\right)\circ\Psi
=
\Psi'\circ
\left(\begin{matrix}\sigma_0(a)&\sigma_1(b)\\d\sigma_0(c)&\sigma_0(e)\end{matrix}\right).
\]
The equivalence $\Win(\u A)'\to\Win(\u A)$ is given by the assignment $(L_0,L_1,\Psi)\mapsto(M,\phi)$ with $M=L_0\otimes_{A_0}A(0)\oplus L_1\otimes_{A_0}A(1)$ and $\phi_1^\sharp=\Psi$ under the identification $(M^\sigma)_1= L_0^{\sigma_0}\oplus L_1^{\sigma_0}$.
\eRe

\bEx
\label{Ex:win-uA-uA(1)}
We have $\u A$-windows $\u A=(A,\sigma)$ and $\u A(1)=(A(1),\phi)$ with $\phi_0=d\sigma$, $\phi_1=\sigma_0$.
\eEx

\subsection{Section complex}

\bDe
\label{De:Gamma(uM)}
For an $\u A$-module $\uM$ with underlying filtered $A$-module $(M_0\xleftarrow{\:\tau} M_1)$ we define a complex of abelian groups in degrees $[0,1]$
\[
\Gamma(\uM)=[M_1\xrightarrow{\;\gamma\;}M_0],
\qquad
\gamma(x)=\phi_1(x)-\tau(x).
\]
\eDe

\bLe
\label{Le:Hi-Gamma-uM}
For $i=0,1$ we have a natural isomorphism
\[
H^i(\Gamma(\uM))\cong\Ext^i_{\Mod(\u A)}(\u A(1),\uM).
\]
\eLe

\bproof
The case $i=0$ is clear. An extension of $\u A$-modules $0\to\uM\to\uN\to\u A(1)\to 0$ splits on the level of filtered $A$-modules. If $N=M\oplus A(1)$ as a filtered $A$-module, then $\phi_N=\phi_M+\phi_{A(1)}+\psi$ for a $\sigma$-linear map $\psi:A(1)\to M_0\otimes_{A_0}A'$. Such $\psi$ correspond to elements of $M_0$, and the group $\Hom_A(A(1),M)=M_1$ acts on $\psi$ via $\gamma$.
\eproof

\bRe[Functoriality I]
The complex $\Gamma(\uM)$ is functorial in $\uM$, and a strict frame homomorphism $\u A\to\u B$ induces a homomorphism of complexes $\Gamma(\uM)\to\Gamma(\uM\otimes_{\u A}\u B)$.
\eRe

\bRe[Functoriality II]
Let $(g,u):\u A\to \u B$ be a general frame homomorphism, and assume that $\alpha\in B_0^*$ with $\alpha=u\sigma_0(\alpha)$ is given. Then multiplication by $\alpha$ again induces a homomorphism of complexes $\Gamma(\uM)\to\Gamma(\uM\otimes_{\u A}\u B$).
\eRe

\subsection{Duality}

Let $A$ be a filtered ring.

\bDe
For $M,N\in\Mod(A)$, a bilinear form $h:M\times N\to A$ consists of bilinear maps of $A_0$-modules $h_i\colon M_i\times N_i\to A_i$ for $i=0,1$ such that the following diagrams commute,  
where $h_0'(m,n)=h_0(m,\tau (n))$ and $h_0''(m,n)=h_0(\tau (m),n)$.
\[
\xymatrix@M+0.2em{
M_0 \times N_0 \ar[r]^-{h_0} & A_0 \\
M_1\times N_1 \ar[u]^{\tau\times\tau} \ar[r]^-{h_1} & A_1 \ar[u]_\tau
}
\qquad
\xymatrix@M+0.2em{
M_1\times N_1 \ar[r]^-{h_1} & A_1
\\
A_1\times M_0\times N_1 \ar[u]^{m\times\id} \ar[r]^-{{\id}\times h_0'} & A_1 \times A_0 \ar[u]_m
}
\qquad
\xymatrix@M+0.2em{
M_1\times N_1 \ar[r]^-{h_1} & A_1
\\
M_1\times N_0\times A_1 \ar[u]^{{\id}\times m} \ar[r]^-{h_0''\times{\id}} & A_0 \times A_1 \ar[u]_m
}
\]
Let $\Bil(M\times N,A)$ denote the abelian group of bilinear forms.
\eDe

\bLe
\label{Le:dual-filt-module}
For each filtered $A$-module $M$ there is a dual filtered $A$-module $M^t$ with a universal bilinear form $M^t\times M\to A$, i.e.\ the resulting map $\Hom_A(N,M^t)\to\Bil(N\times M,A)$ is bijective for every filtered $A$-module $N$. For $M=L\otimes_{A_0}A(i)$ we have $M^t=L^\vee\otimes_{A_0}A(1-i)$. In particular $M\to M^{tt}$ is an isomorphism when $M$ is finite projective.
\eLe

\bproof
We have $M^t_0=\Hom_{A_0}(M_0,A_0)$ and $M^t_1=\Hom_A(M,A)$ such that $M^t_1\to M^t_0$ is given by restriction, the multiplication $A_1\otimes_{A_0}M^t_0\to M^t_1$ sends $a\otimes u$ to the
homomorphism $M\to A$ defined by $M_0\xrightarrow {au} A_1\to A_0$ in degree $0$ and by $M_1\to M_0\xrightarrow{au}A_1$ in degree $1$, and the bilinear form $M^t\times M\to A$ is given by evaluation. The homomorphism 
$\Hom_A(N,M^t)\to\Bil(N\times M,A)$ is bijective when $M$ is free, hence in general by choosing a free resolution of $M$.
The last assertion follows since finite projective filtered modules are direct summands of finite free filtered modules.
\eproof

\bRe
\label{Re:bil-filt-funct}
Let $u:A\to B$ be a homomorphism of filtered rings. Each bilinear form of $A$-modules $N\times M\to A$ induces a bilinear form of $B$-modules $N^u\times M^u\to B$. We get a natural homomorphism $(M^t)^u\to(M^u)^t$, which is an isomorphism when $M$ is finite projective.
\eRe

\bRe
\label{Re:bil-filt-Cartier}
Let $A$ be a filtered ring such that $A_1$ is an invertible $A$-module. For every finite projective filtered $A$-module $M$ the bilinear form $M^t\times M\to A$ is perfect in both degrees and thus $(M^t)_1=\Hom_{A_0}(M_1,A_1)$. 
\eRe

Now let again $\u A$ be a frame.

\bDe
For $\u A$-modules $\uM,\uN$, a bilinear form $h\colon\uM\times\uN\to\u A$ is a bilinear form of filtered $A$-modules $h=(h_0,h_1)\colon M\times N\to A$ such that the following diagrams commute. 
\[
\xymatrix@M+0.2em{
M_0 \times N_0 \ar[r]^-{h_0} & A_0 \\
M_1\times N_1 \ar[u]^{\phi_1\times\phi_1} \ar[r]^-{h_1} & A_1 \ar[u]_{\sigma_1}
}
\qquad
\qquad
\xymatrix@M+0.2em{
M_0 \times N_0 \ar[r]^-{h_0} & A_0 \\
M_0\times N_0 \ar[u]^{\phi_0\times\phi_0} \ar[r]^-{h_0} & A_0 \ar[u]_{d\sigma_0}
}
\]
Let $\Bil(M\times N,\u A)$ denote the abelian group of bilinear forms.
\eDe

\bPr
For each $\u A$-window $\uM=(M,\phi)$ there is a dual $\u A$-window $\uM^t=(M^t,\phi^t)$ with a universal bilinear form $\uM^t\times\uM\to\u A$, i.e.\ for every $\u A$-module $\uN$ the resulting map $\Hom_{\u A}(\uM^t,\uN)\to\Bil(\uN\times\uM,\u A)$ is bijective.
\ePr

\bproof
The $A$-module $M^t$ is given by Lemma \ref{Le:dual-filt-module}. The homomorphism $\phi$ corresponds to an isomorphism of $A_0$-modules $(M^\sigma)_1\xrightarrow\sim M_0$. The inverse of the dual of this isomorphism defines $\phi^t$, 
using that $((M^\sigma)_1)^\vee=(M^{\sigma t})_1=(M^{t\sigma})_1$ by Remarks \ref{Re:bil-filt-Cartier} and \ref{Re:bil-filt-funct}.
\eproof

\bRe
\label{Re:dual-normal-decomp}
The resulting duality operation in the category $\Win'(\u A)$ of Remark \ref{Re:normal-repres} is given by $(L_0,L_1,\Psi)^t=(L_1^\vee,L_0^\vee,(\Psi^{-1})^\vee)$.
\eRe

\subsection{Deformation lemma}

Let us recall the deformation lemma for windows.

\bDe
Let $\u A=(A,d,\sigma)$ be a frame, so $A=(A_0\xleftarrow{\:\tau} A_1)$ is a filtered ring. An ideal of $\u A$ is a filtered $A$-submodule $K\subseteq A$ such that $\sigma_1\colon A_1\to A_0$ induces $\sigma_{1,K}\colon K_1\to K_0$, and $K_0$ and $K_1$ are derived $(p,d)$-complete. 
We can form the quotient frame $\u A/K=(A/K,\bar d,\bar\sigma)$.
An ideal $K$ of $\u A$ will be called leveled if $\tau_K\colon K_1\to K_0$ is bijective. In this case we have a $\sigma_0$-linear map 
$\dot\sigma_K=\sigma_{1,K}\circ\tau_K^{-1}\colon K_0\to K_0$.
\eDe

\bLe
\label{Le:deformation-equivaelence-windows}
Let $K$ be a leveled ideal of the frame $\u A$ such that $\dot\sigma_K$ is pointwise nilpotent. Then the functor $\Win(\u A)\to\Win(\u A/K)$ is an equivalence of categories.
\eLe

This is standard. For completeness we sketch a proof using the present formalism.

\bproof
Let $\u{\bar A}=\u A/K$.
It suffices to show that the induced functor of groupoid cores is an equivalence.
For a given finite projective $A$-module $M$ let $\Win(\u A)_M$ be the groupoid of $\u A$-windows with underlying $A$-module isomorphic to $M$.
We have $A_0/\tau A_1=\bar A_0/\tau\bar A_1$ since $K$ is leveled, moreover $(A_0,\tau A_1)$ and $(\bar A_0,\tau\bar A_1)$ are Henselian pairs by Lemma \ref{Le:frame-Henselian}.
Hence the homomorphism $A_0\to\bar A_0$ induces a bijective map on the sets of isomorphism classes of finite projective modules, the same holds for the homomorphism $A\to\bar A$ by Lemma \ref{Le:fil-proj-mod}, and it suffices to show that the functor $\Win(\u A)_M\to\Win(\u{\bar A})_{\bar M}$ is an equivalence where $\bar M=M\otimes_A\bar A$. 

This follows from a calculation similar to \cite[Proposition 7.1.5]{Lau:Higher}.
Explicitly, $\Win(\u A)_M$ is equivalent to the quotient groupoid $[X/G]$ for $X=\Isom_{A_0}((M^\sigma)_1,M_0)$ and $G=\Aut_A(M)$ with respect to the conjugation action, similarly $\Win(\u{\bar A})_{\bar M}=[\bar X/\bar G]$. The homomorphism $G\to\bar G$ and the map $X\to\bar X$ are surjective, so we have to show that $G_K=\ker(G\to\bar G)$ acts simply transitively on the fibers of $X\to\bar X$. 
One verifies that $G_K$ is isomorphic to $G_K'=\ker(\Aut(M_0)\to\Aut(\bar M_0))$ and that the action on the fibers is simply transitive iff for every $x\in X$ the map of sets
\beqn
\label{Eq:transitive-bijective}
G_K'\to G_K',
\qquad
h\mapsto h^{-1}x\dot\sigma_G(h)x^{-1}
\eeqn
is bijective, where $\dot\sigma_G$ is described as follows. Assume that  $M=L_0\otimes_{A_0}A(0)\oplus L_1\otimes_{A_0}A(1)$ and write $h=1+(a_{ij})$ with $a_{ij}\in\Hom(L_j,L_i)\otimes K_0$. Then $\dot\sigma(h)=1+(b_{ij})$ with $b_{ij}=d^{1+i-j}\dot\sigma_K(a_{ij})$ as an elment of $\Hom(L_j,L_i)^{\sigma_0}\otimes K$. It follows that the operator $h\mapsto x\dot\sigma_G(h)x^{-1}$ is locally nilpotent, and hence \eqref{Eq:transitive-bijective} is bijective as required.
\eproof

\subsection{Hodge filtration}  
\label{Se:Hodge-fil}

Let $A$ be a filtered ring where $A_1\to A_0$ is injective and let $R=A_0/A_1$. The Hodge filtration of a finite projective $A$-module $M$ is the $R$-module $M_1/A_1M_0\subseteq M\otimes_AR$; this is a direct summand.

\bLe
Let $\u A\to\u B$ be a frame homomorphism with bijective underlying ring homomorphism $A_0\to B_0$. Assume that $A_1\to A_0$ and $B_1\to B_0$ are injective with quotients $R=A_0/A_1$ and $S=B_0/B_1$. 
The category $\Win(\u A)$ is equivalent to the category of pairs $(\uN,L)$ where $\uN\in\Win(\u B)$ and $L\subseteq N_0\otimes_{B_0}R$ is a direct summand that lifts the Hodge filtration of $N$.
\eLe

\bproof
The functor sends an $\u A$-module $\uM$ to the pair $(\uN,L)$ where $\uN=\uM\otimes_{\u A}\u B$ and $L$ is the Hodge filtration of $M$.
\eproof

\cclearpage
\section{Displays and group schemes}

In this section, we recall displays and some of their variants, and the construction of formal groups and group schemes associated with displays and truncated displays, adding a dual version in the truncated case.

\subsection{Witt vector frames and displays}

Windows over frames associated to Witt vectors are called displays of various type.

\bEx
[Witt frame]
\label{Ex:Witt-frame}
For every $p$-complete ring $R$ there is a crystalline frame $\u W(R)$ with underlying filtered ring given by $V\colon F_*W(R)\to W(R)$, which is the kernel of $W(R)\to R$, the homomorphism of filtered rings
\[
\sigma\colon (W(R)\xleftarrow{\:V}F_*W(R))\longrightarrow(W(R)\xleftarrow{\: p} W(R))
\] 
is $F$ in degree $0$ and $\id$ in degree $1$, and the divided powers on $VW(R)$ are determined by $\tilde\gamma_p(V(x))=p^{p-2}V(x^p)$ for $x\in W(R)$. The relations \eqref{Eq:pd-frame-1} and \eqref{Eq:pd-frame-2} hold when $R$ is $p$-torsion free, thus in general by functoriality. 
Windows over $\u W(R)$ are called displays over $R$, notation $\Disp(R)=\Win(\u W(R))$. This category is equivalent to the category of $3n$-displays of \cite{Zink:The-display}.
\eEx

\bEx
[Truncated Witt frame]
\label{Ex:Witt-frame-n}
For every $\FF_p$-algebra $R$ the ring of truncated Witt vectors $W_n(R)$ carries the structure of a frame $\u {W\!}_n(R)$ with $d=p$ and filtration given by the quasi-ideal $V\colon F_*W_n(R)\to W_n(R)$, such that the projection $W(R)\to W_n(R)$ extends to a strict frame homomorphism $\u W(R)\to\u {W\!}_n(R)$. 
Windows over $\u {W\!}_n(R)$ are called $n$-truncated displays over $R$, notation $\Disp_n(R)=\Win(\u{W\!}_n(R))$. This category is equivalent to the category of $n$-truncated displays of \cite{Lau:Smoothness}.
\eEx

\bEx
[Relative Witt frame]
\label{Ex:Witt-frame-rel}
Let $R'\to R$  a pd thickening of $p$-nilpotent rings with kernel $\Fa$.
The divided powers induce an isomorphism of $W(R')$-modules $w'\colon W(\Fa)\xrightarrow\sim\prod_{\NN}\Fa$; see \cite[(48)]{Zink:The-display}. Let $i\colon \Fa\to W(\Fa)$ be defined by $w'(i(a))=(a,0,0,\ldots)$. 
We have a crystalline frame $\u W(R'/R)$ with underlying filtered ring given by $V\oplus i\colon F_*W(R')\oplus\Fa\to W(R')$, 
which is the kernel of $W(R')\to R$, the homomorphism of filtered rings
\[
\sigma\colon (W(R')\xleftarrow{\:V\oplus i}F_*W(R')\oplus \Fa)\longrightarrow(W(R')\xleftarrow{\: p} W(R'))
\] 
is $F$ in degree $0$ and ${\id}\oplus 0$ in degree $1$, and the divided powers on $VW(R')\oplus i(\Fa)$ are determined by $\tilde\gamma_p(V(x))=p^{p-2}V(x^p)$ for $x\in W(R')$ and $\gamma_n(i(a))=i(\gamma_n(a))$ for $a\in\Fa$. The relations \eqref{Eq:pd-frame-1} and \eqref{Eq:pd-frame-2} are proved by reduction to the case where $R'$ is $p$-torsion free. Windows over $\u W(R'/R)$ are called relative displays for $R'\to R$, notation $\Disp(R'/R)=\Win(\u W(R'/R))$. This category is equivalent to the category of triples relative to $R'\to R$ of \cite[\S 2.2]{Zink:The-display}.
\eEx

\subsection{Zink functor}
\label{Se:Recall-Zink-functor}

For every $p$-nilpotent ring $R$, by \cite{Zink:The-display} there is a functor 
\beqn
\label{Eq:FZR-DispR-LGR}
\FZ_R:\Disp(R)\to\FG(R)
\eeqn
from displays over $R$ to commutative formal groups over $R$, denoted $BT$ in loc.\,cit., with $\dim(\FZ(\uM))=\dim(\uM)$. Let us recall the definition in terms of the present frame formalism. 
Let $\hat W(R)$ be the set of Witt vectors $(a_0,a_1,\ldots)\in W(R)$ where all $a_i$ are nilpotent and $a_i=0$ for sufficiently large $i$. Then $\hat W(R)$ is a $W(R)$-module, which is upgraded to a $\u W(R)$-module $\u{\hat W}(R)\subseteq\u W(R)$ using the filtration $V\colon F_*\hat W(R)\to\hat W(R)$. For $\uM\in\Disp(R)$ we form a complex $Z(\uM)$ of abelian fpqc sheaves on $\Aff_R$ in degrees $[-1,0]$ defined by
\beqn
\label{Eq:ZMN}
Z(\uM)(S)=\Gamma(\uM\otimes_{\u W(R)}\u{\hat W}(S))[1]
\eeqn
using the tensor product of Remark \ref{Re:tensor-Mod-frame} and the complex $\Gamma$ of Definition \ref{De:Gamma(uM)}.
Then $H^{-1}(Z(\uM))$ vanishes, and $\FZ_R(\uM)(S)=H^0(Z(M)(S))$. Here $H^0$ is taken as presheaves. The formula \eqref{Eq:ZMN} remains valid when $S$ is a nilpotent $R$-algebra.  We refer to \cite[\S11]{Lau:Gerbe} for further details.

\subsection{Truncated Zink functor}
\label{Se:Truncated-Zink-functor}

If $R$ is an $\FF_p$-algebra, a similar construction (originally observed in \cite{Lau-Zink}) yields a functor
\beqn
\label{Eq:FZnR-DispnR-LGnR}
\FZ_{n,R}\colon \Disp_n(R)\to\FG_n(R)
\eeqn
where $\FG_n(R)$ is the category of commutative $n$-smooth group schemes as defined in \cite{Drinfeld:The-Lau}\footnote{A group scheme $G$ over $R$ is called $n$-smooth if the underlying pointed scheme is isomorphic to $\Spec R[x_1,\ldots,x_r]/(x_1^{p^n},\ldots,x_r^{p^n})$ locally in $\Spec R$. By \cite[Theorem~A \& Proposition~4.3.3]{Kothari-Mundinger}, this holds iff locally $G$ is the $F^n$-torsion of a commutative formal group.} such that $\rk(\Lie\FZ_n(\uM))=\dim(\uM)$. Explicitly, for $\uM\in\Disp_n(R)$ we form a complex $Z_n(\uM)$ of abelian fpqc sheaves on $\Aff_R$ in degrees $[-1,0]$ defined by
\beqn
\label{Eq:ZnM}
Z_n(\uM)(S)=\Gamma(\uM\otimes_{\u{W\!}_n(R)}\u{\hat W}(S)[F^n])[1].
\eeqn
Then $H^{-1}(Z_n(\uM))$ vanishes, and $\FZ_{n,R}(\uM)(S)=H^0(Z_n(\uM)(S))$. Again we refer to \cite[\S11]{Lau:Gerbe} for further details.

\subsection{Dual truncated Zink functor}
\label{Se:Dual-Zink-functor}

Let $R$ be an $\FF_p$-algebra.
Following \cite{Drinfeld:The-Lau}, a finite locally free group scheme over $R$ is called $n$-cosmooth if its Cartier dual is $n$-smooth.

\bDe
\label{De:CnM}
For $\uM\in\Disp_n(R)$ let $C_n(\uM)=\Gamma(\tilde{\uM})=[\tilde M_1\xrightarrow{\gamma}\tilde M_0]$ with $\gamma=\phi_1-\tau$ as a complex of fpqc sheaves of abelian groups over $\Spec R$ in degrees $[0,1]$ where $\tilde{\uM}$ is the sheaf of $n$-truncated displays defined by base change of $\uM$.
\eDe

\bPr
\label{Pr:dual-truncated-Zink-functor}
Let $\uM\in\Disp_n(R)$.
Then $H^1(C_n(\uM))_{\syn}$ vanishes, i.e.\ $\gamma\colon\tilde M_1\to\tilde M_0$ is syntomic surjective, and $H^0(C_n(\uM))$ is representable by an $n$-cosmooth commutative group scheme, which we denote by $\FY_n(\uM)$. We have $\rk(\Lie(\FY_n(\uM)^\vee))=\dim(\uM^t)$.
\ePr

\bRe
This result is related to \cite[\S7]{Drinfeld:The-Lau}, see Remark \ref{Re:Drinfeld-S7} below.
\eRe

\bproof
Initial remark.
The complex $C_n(\uM)$ is a complex of smooth commutative unipotent group schemes. Hence $H^0(C_n(\uM))$ is an affine group scheme. If this group scheme is finite locally free (as will be shown),  
then $\gamma$ is automatically a syntomic cover, and the vanishing of $H^1$ follows. This remark will not be used.

\subsubsection*{The case $n=1$} 

We claim that as a complex of presheaves, $C_1(\uM)$ is quasi-isomorphic to a complex $[\tilde L\xrightarrow u \tilde L^{(p)}]$ where $L$ is a finite projective $R$-module of rank $c=\dim(\uM^t)$ and $u(x)=x\otimes 1-h(x)$ for an $R$-linear map $h$.
This complex has trivial $H^1$, and $H^0$ is a $1$-cosmooth group scheme of degree $p^c$ by \cite[Proposition 2.2]{deJong:Finite-locally-free}, based on \cite[VII$_{\text A}$, Th\'eor\`eme 7.4]{SGA3}.

Let us prove the claim.
For any frame $\u A$ and any $\u A$-module $\uM=(M,\phi)$ such that $\phi^\sharp_1\colon(M^\sigma)_1\to M_0$ is bijective, the complex $\Gamma(\uM)$ is isomorphic to
$\Gamma'(\uM)=[\gamma'\colon M_1\to (M^\sigma)_1]$
with $\gamma'={\id}\otimes\sigma-\tau'$ where $\tau'=(\phi^\sharp_1)^{-1}\circ\tau$ is $A_0$-linear.
For $\u A=\u{W\!}_1(R)$ the homomorphism $\tau\colon A_1\to A_0$ is zero and $\sigma_1\colon A_1\to A_0$ is bijective, moreover $A_0=R$. If $\uM$ is an $\u A$-window, thus $M=L_0\otimes_{A_0}A(0)\oplus L_1\otimes_{A_0}A(1)$ as an $A$-module, it follows that
\[
\Gamma(\uM)\cong {}
\Gamma'(\uM)=[(L_0\otimes A_1)\oplus L_1\xrightarrow{\gamma'}(L_0\otimes\sigma_{0*}A_0)\oplus (L_1\otimes\sigma_{0*}A_0)],
\]
which contains the acyclic complex $L_0\otimes[A_1\xrightarrow{\sigma_1}\sigma_{0*}A_0]$. The quotient takes the desired form.

\subsubsection*{The case $n\ge 2$} 

We use induction on $n$. Let $0<r<n$. For an $R$-algebra $S$, the Frobenius induces an endomorphism $F$ of the frame $\u{W\!}_n(S)$, and we have an exact sequence of $\u{W\!}_n(S)$-modules
\[
0\to F^r_*\u{W\!}_{n-r}(S)\xrightarrow {\;v_r\;} \u{W\!}_n(S)\xrightarrow{\;\pi_r\;} \u{W\!}_r(S)\to 0
\]
where $v_r$ is induced by $V^r$. The sequence remains exact under $\uM\otimes_{\u{W\!}_n(R)}-$ and hence induces an exact sequence of complexes of presheaves of abelian groups
\[
0\to
C_{n-r}(F^{r*}\uM_{n-r})\to C_{n}(\uM)\to C_{r}(\uM_{r})\to 0
\]
where $\uM_m\in\Disp_m(R)$ is the reduction of $\uM$.
By induction it follows that $H^1(C_n(\uM))=0$ and that $G_n=H^0(C_n(\uM))$ is a finite locally free group scheme 
sitting in an exact sequence 
\[
0\to F^{r*}G_{n-r}\xrightarrow{\;v_r\;} G_n\xrightarrow{\;\pi_r\;} G_r\to 0.
\]

It remains to verify that $G_n$ is $n$-cosmooth.
Both the components of the complex $C_n(\uM)$ and its cohomology are (affine) flat commutative group schemes and hence carry a functorial Verschiebung homomorphism by \cite[VII$_{\text A}$, \S4.3]{SGA3}. Since the Verschiebung of the group scheme $W_n$ is the Witt vector Verschiebung $V$ and since $V^r:F^r_*W_n(S)\to W_n(S)$ factors as $V^r=v_r\circ\pi_{n-r}$ it follows that the sequence  
\[
G_n^{(p^r)}\xrightarrow{V_{G_n}^r} G_n\xrightarrow{\pi_r} G_r\to 0
\] 
is exact, which implies that $G_n$ is $n$-cosmooth by degree considerations. 
\eproof

\subsection{A duality result}

For completeness we record the following duality that will not be used later.
For an fpqc sheaf of abelian groups $K$ let $K^\vee=\u\Hom(K,\GG_m)$ be the Cartier dual. 

\bPr
\label{Pr:FZ-FY-dual}
For every $n$-truncated display $\uM\in\Disp_n(R)$ with dual $\uM^t$ there is a natural quasi-isomorphism $C_n(\uM)\cong Z_n(\uM^t)^\vee$ of complexes of presheaves of abelian groups, thus also of fpqc sheaves, and hence a natural isomorphism of group schemes $\FY_n(\uM)\cong\FZ_n(\uM^t)^\vee$.
\ePr 

See Definition \ref{De:CnM} for $C_n(\u M)$, Proposition \ref{Pr:dual-truncated-Zink-functor} for $\FY_n(\uM)$, and \S\ref{Se:Truncated-Zink-functor} for $Z_n(\uM)$ and $\FZ_n(\uM)$.

\bproof
Let $\u A=\u{W\!}_n(R)$. 
Assume that $M=L_0\otimes_{A_0}A(0)\oplus L_1\otimes_{A_0}A(1)$ as a filtered $A$-module (normal decomposition of $M$) and let $\Psi:L_0^{F}\oplus L_1^{F}\to L_0\oplus L_1$ be the isomorphism of $A_0$-modules that corresponds to $\phi=\phi_M$ as in Remark \ref{Re:normal-repres}. Then 
\[
C_n(\uM)=[\tilde L_0^F\oplus \tilde L_1^F\xrightarrow\gamma \tilde L_0\oplus \tilde L_1],\qquad
\gamma=\Psi\circ({\id}\oplus F)-V\oplus 1.
\] 
This complex is quasi-isomorphic to 
\[
C_n'(\uM)=[\tilde L_0\oplus \tilde L_1\xrightarrow{\gamma'}\tilde L_0\oplus \tilde L_1^F],\qquad\gamma'={\id}\oplus F-(V\oplus{\id})\circ\Psi^{-1}
\]
because there is a commutative diagram of presheaves
(two commutative squares)
\beqn
\label{Eq:Cn-C'n}
\xymatrix@M+0.2em@C+3em{
\tilde L_0^F \oplus \tilde L_1 \ar@<0.7ex>[r]^{V\oplus{\id}} \ar@<-0.7ex>[r]_{\Psi\circ({\id}\oplus F)} \ar[d]_{V\oplus{\id}} & \tilde L_0 \oplus \tilde L_1 \ar[d]^{(V\oplus{\id})\circ\Psi^{-1}} \\
\tilde L_0 \oplus \tilde L_1 \ar@<0.7ex>[r]^{(V\oplus{\id})\circ\Psi^{-1}} \ar@<-0.7ex>[r]_{{\id}\oplus F} & \tilde L_0 \oplus \tilde L_1^F
}
\eeqn
where on the vertical kernels and on the vertical presheaf cokernels, one of the two parallel maps is zero and the other one is bijective. Since $\hat W[F^n]=W_n^\vee$, using Remark \ref{Re:dual-normal-decomp} it follows that $C_n'(\uM)^\vee\cong Z_n(\uM^t)$. 

In order to make this construction functorial, we now write the complex $C_n'(\uM)$, the quasi-isomorphism $C_n(\uM)\to C_n'(\uM)$, and the isomorphism $C_n'(\uM)^\vee\cong Z_n(\uM^t)$ without referring to a normal decomposition. 

We use that $0\to W_n[F]\to W_n\xrightarrow F W_n\to 0$ is exact as (syntomic) sheaves.
Let $\bar M_0=M_0\otimes_{W_n(R)}R$ and let $\bar L$ be the image of $M_1\to\bar M_0$ (the Hodge filtration of $M$). If a normal decomposition of $M$ is given, then $\bar L=\bar L_1$.
For $i=0,1$ let $K_i$ be the image of the multiplication map $W_n(R)[F]\otimes M_1\to M_i$. Then $K_i=W_n(R)[F]\otimes_R\bar L$, in particular $K_1\to K_0$ is bijective. 
By base change the groups $K_i$ extend to sheaves $\tilde K_i$ over $\Spec R$.
The canonical homomorphism $\tilde M_1\to(\tilde M^\sigma)_1$ is surjective as syntomic sheaves with kernel $\tilde K_1$. 
It follows that $\tau\colon\tilde M_1\to\tilde M_0$ induces a homomorphism of sheaves $v\colon(\tilde M^\sigma)_1\to\tilde M_0/\tilde K_0$ (sheaf quotient).
Now the coordinate free version of \eqref{Eq:Cn-C'n} looks as follows.
\[
\xymatrix@M+0.2em@C+3em{
\tilde M_1 \ar@<0.7ex>[r]^{\tau} \ar@<-0.7ex>[r]_{\phi} \ar[d]_{\tau} & \tilde M_0 \ar[d]^{v\circ(\phi_1^\sharp)^{-1}} \\
\tilde M_0 \ar@<0.7ex>[r]^{v\circ(\phi_1^\sharp)^{-1}} \ar@<-0.7ex>[r]_{\can} & \tilde M_0/\tilde K_0
}
\]
It remains to identify the Cartier dual of $\tilde M_0$ and $\tilde M_0/\tilde K_0$ with the components of $Z_n(\uM^t)$ compatibly with the previous identification derived from a normal decomposition. But the  Cartier dual of $\tilde K_0\to\tilde M_0$ identifies with $M_0^t\otimes_{W_n(R)}\hat W[F^n]\to (M_0^t/\tau M_1^t)\otimes_R\GG_a[F^n]$, and the kernel of this map is $(M^t\otimes\u{\hat W}[F^n])_1$ as required.
\eproof

\bRe
\label{Re:Drinfeld-S7}
Proposition \ref{Pr:FZ-FY-dual} exhibits $\FY_n(\uM)$ as the Cartier dual of $\FZ_n(\uM)$ and thus gives an alternative proof of Proposition \ref{Pr:dual-truncated-Zink-functor}. Similar results are proved in \cite[Proposition 7.2.3, Lemma 7.3.2]{Drinfeld:The-Lau} in the context of (weak) semidisplays.
\eRe

\cclearpage
\section{Sheared displays}

This section introduces plain sheared displays as the analogue of displays with $\sW$ in place of $W$, and sheared displays as the associated stack for the $p^\infty$-root topology. With hindsight (by the equivalence with $p$-divisible groups) this will coincide with the associated fpqc stack.

\subsection{Sheared Witt frames}

In \S\ref{Se:modif-V} we defined $u_0\in W$, $\tilde p\in\sW$, and $\tilde V\colon\sW\to\sW$.

\bDe
For every $p$-nilpotent ring $R$ we define a prismatic frame $\usW(R)$ with underlying  
filtered ring given by $\tilde V\colon F_*\sW(R)\to\sW(R)$, which is the kernel of $\sW(R)\to R$, equipped with the distinguished element $\tilde p$,
and the homomorphism of filtered rings
\[
\sigma\colon (\sW(R)\xleftarrow{\:\tilde V}F_*\sW(R))\longrightarrow(\sW(R)\xleftarrow{\: \tilde p}\sW(R))
\]
is defined to be $F$ in degree $0$ and $\id$ in degree $1$. 

The ring $\sW(R)$ is derived $(p,\tilde p)$-complete by Corollary \ref{Co:sW-p-complete}.
\eDe

\bRe
\label{Re:usW-pd-frame}
If $R$ is an $\FF_p$-algebra or $p\ge 3$ then $\usW(R)$ 
is naturally a crystalline frame. Indeed, Lemma \ref{Le:pd-sW} gives divided powers on the ideal $\tilde V\sW(R)$. The relations \eqref{Eq:pd-frame-1} and \eqref{Eq:pd-frame-2} hold in $Q\limperf(R)$ since that ring is $p$-torsion free, and they hold in $W(R)$ by Example \ref{Ex:Witt-frame}.
\eRe

\bRe
\label{Re:Hom-usW-uW}
The natural ring homomorphism $c\colon\sW(R)\to W(R)$ extends to homomorphism of prismatic frames $(c,u_0)\colon\usW(R)\to\u W(R)$. Moreover, there is a unique unit $\alpha\in W(\ZZ_p)$ that maps to $1$ in $W(\FF_p)$ with $F(\alpha) u_0=\alpha$, namely $\alpha=\prod_{n\ge 0}F^n(u_0)$.
\eRe

\subsection{Relative sheared Witt frames}

Let $R'\to R$ by a nilpotent pd thickening of $p$-nilpotent rings with kernel $\Fa$. The divided powers induce an isomorphism $w':\hat W(\Fa)\xrightarrow\sim\bigoplus_{\NN}\Fa$; see \cite[(149)]{Zink:The-display}. As  earlier let $i:\Fa\to\hat W(\Fa)$ be defined by $w'(i(a))=(a,0,0,\ldots)$. 
We consider $\hat W(\Fa)$ as an ideal of $\sW(R')$ by \eqref{Eq:WN-sWR-sWRN}.

\bDe
We define a prismatic frame $\usW(R'/R)$ with underlying 
filtered ring given by $\tilde V\oplus i\colon F_*\sW(R')\oplus\Fa\to\sW(R')$, which is the kernel of $\sW(R')\to R$, equipped with the distinguished element $\tilde p$, 
and the homomorphism of filtered rings
\[
\sigma\colon (\sW(R')\xleftarrow{\:\tilde V\oplus i}F_*\sW(R')\oplus \Fa)\longrightarrow(\sW(R')\xleftarrow{\: \tilde p} \sW(R'))
\]
is defined to be $F$ in degree zero and ${\id}\oplus 0$ in degree $1$.
\eDe

\bRe
\label{Re:usW-rel-pd-frame}
If $R'$ is an $\FF_p$-algebra or $p\ge 3$ then $\usW(R'/R)$ 
is naturally a crystalline frame. Indeed, Lemma \ref{Le:pd-sW-rel} gives divided powers on the filtration ideal. The relations \eqref{Eq:pd-frame-1} and \eqref{Eq:pd-frame-2} hold in $Q\limperf(R')$ since that ring is $p$-torsion free, and they hold in $W(R')$ by example \ref{Ex:Witt-frame-rel}.
\eRe

\bRe
\label{Re:Hom-usW-uW-rel}
As in Remark \ref{Re:Hom-usW-uW} there is a natural frame homomorphism $\usW(R'/R)\to\u W(R'/R)$.
\eRe

\subsection{Sheared displays}

\bDe
\label{De:Sh-disp}
A plain sheared display over a $p$-nilpotent ring $R$ is a window over $\usW(R)$. 
These objects form an exact category 
$\sDisp^{\plain}(R)=\Win(\usW(R))$,
and we obtain a fibered category
$\sDisp^{\plain}\to\Aff_{\Spf\ZZ_p}$
since a ring homomorphism $R\to R'$ induces a frame homomorphism $\usW(R)\to\usW(R')$ and thus a functor between the window categories. Let
\[
\sDisp\to\Aff_{\Spf\ZZ_p}
\]
be the associated $p^\infty$-root stack (Definition \ref{De:proot}).
Its objects will be called sheared displays.
\eDe

\bPr
\label{Pr:Descent-Sh-disp-essp}
The fibered category of plain sheared displays over essentially semiperfect rings (Definition \ref{De:essemiperf}) is an ind-\'etale and $p^{\infty}$-root stack over over $\Aff^{\essemiperf}$.
\ePr

\bproof
For a plain sheared display $\uM=(M,\phi)$ over an essentially semiperfect ring $R$, the filtered module $M=(M_0\leftarrow M_1)$ corresponds to a pair $(M_0,Q)$ where $M_0$ is a finite projective $\sW(R)$-module and $Q$ is an $R$-projective quotient of $M_0\otimes_{\sW(R)}R$ (using Corollary \ref{Co:sWR-R-surj} and Remark \ref{Re:fil-proj-mod-ideal}). Such pairs satisfy ind-\'etale and $p^\infty$-root descent by Proposition \ref{Pr:Descent-LF(sW)-essemiperf}. If $M$ is given, the extensions to a sheared display $(M,\phi)$ form an fpqc sheaf of sets since $\sW$ is an fqpc sheaf.
\eproof

\bCo
\label{Co:sDissp-ind-etale}
The fibered category $\sDisp\to\Aff_{\Spf\ZZ_p}$ is an ind-\'etale stack.
\eCo

\bCo
\label{Co:essemiperf-plain}
The inclusion $\sDisp^{\plain}\to\sDisp$ is an equivalence over $\Aff^{\essemiperf}$, i.e.\ over essentially semiperfect rings, every sheared display is plain.
\eCo

\bproof[Proof of Corollaries \ref{Co:sDissp-ind-etale} and \ref{Co:essemiperf-plain}]
Use Proposition \ref{Pr:Descent-Sh-disp-essp} and the fact that for each $p$-nilpotent ring $R$ there is an $p^\infty$-root morphism $R\to R'$ where $R'$ is essentially semiperfect.
\eproof

\bRe
\label{Re:sDisp-fpqc-stack}
The fibered category $\sDisp\to\Aff_{\Spf\ZZ_p}$ is also an fpqc stack.
This follows from results of \cite{BMVZ} as explained in Corollary \ref{Co:sheared-displays-fpqc-stack}, but it is also a consequence of the equivalence with $p$-divisible groups (Corollary \ref{Co:sFZ-equiv}).
\eRe

\bLe
\label{Le:sDisp-form-smooth-essp}
The fibered category of sheared displays over $\Aff^{\essemiperf}$ is formally smooth.
\eLe

\bproof
Let $R'\to R$ be a surjective homomorphism of essentially semiperfect rings with nilpotent kernel. Then $\sW(R')\to\sW(R)$ is surjective with nilpotent kernel by Proposition \ref{Pr:WN-sWR-sWRN}. By Corollary \ref{Co:essemiperf-plain} and Remark \ref{Re:normal-repres} it follows that every sheared display over $R$ lifts to $R'$.
\eproof

\bRe
\label{Re:sDisp-Disp}
The frame homomorphism $\usW(R)\to\u W(R)$ of Remark \ref{Re:Hom-usW-uW} gives a base change functor $\sDisp(R)\to\Disp(R)$. On the underlying filtered modules, this functor is given by classical $p$-completion because $W$ is the classical $p$-completion of $\sW$, using that $W(R)$ is the classical $p$-completion of $\sW(R)$ when $R$ is essentially semiperfect. 
\eRe

\subsection{Relative sheared displays}
\label{Se:rel-sh-disp}

Let $\nCrys_{\ZZ_p}$ be the category of nilpotent pd thickenings $Y\to X$ of $p$-nilpotent affine schemes compatible with the canonical divided powers of $p$. This category becomes a site for the ind-\'etale and $p^\infty$-root topology, where a family of morphisms $(Y_i\to X_i)\to(Y\to X)$ is a covering if $Y_i\xrightarrow\sim Y\times_XX_i$ for all $i$ and the family $X_i\to X$ is a covering.
The following is parallel to Definition \ref{De:Sh-disp}.

\bDe
\label{De:Sh-disp-rel}
A plain relative sheared display over $\Spec R\to\Spec R'$ in $\nCrys_{\ZZ_p}$ is a window over the frame $\usW(R'/R)$. 
These objects form a category $\sDisp^{\plain}(R'/R)=\Win(\usW(R'/R))$, and we obtain a fibered category $\sDisp_{\rel}^{\plain}\to\nCrys_{\ZZ_p}$ by base change. 
Let 
\[
\sDisp_{\rel}\to\nCrys_{\ZZ_p}
\]
be the associated $p^\infty$-root stack. Its objects will be called relative sheared displays.
\eDe

Let $\nCrys_{\ZZ_p}^{\essemiperf}\subseteq\nCrys_{\ZZ_p}$ be the full subcategory of all $Y\to X$ where $X$ and $Y$ are essentially semiperfect. The discussion after Definition \ref{De:Sh-disp} carries over to the relative case:

\bPr
\label{Pr:Descent-Sh-disp-rel-essp}
The fibered category of plain relative sheared displays over nilpotent pd thickenings of essentially semiperfect rings 
is an ind-\'etale and $p^{\infty}$-root stack over over $\nCrys_{\ZZ_p}^{\essemiperf}$.
\ePr

\bproof
Analogous to Proposition \ref{Pr:Descent-Sh-disp-essp} with $\sW(R')\to R$ in place of $\sW(R)\to R$.
\eproof

\bCo
The fibered category $\sDisp_{\rel}\to\nCrys_{\ZZ_p}$ is an ind-\'etale stack. \qed
\eCo

\bCo
\label{Co:essemiperf-rel-plain}
The inclusion $\sDisp_{\rel}^{\plain}\to\sDisp_{\rel}$ is an equivalence over $\nCrys_{\ZZ_p}^{\essemiperf}$. \qed
\eCo

\bPr
\label{Pr:sDisp(R'/R)-sDisp(R)}
For every nilpotent pd thickening of essentially semiperfect rings $R'\to R=R'/\Fa$ the functor $\sDisp(R'/R)\to\sDisp(R)$ is an equivalence.
\ePr

\bproof
By Corollaries \ref{Co:essemiperf-plain} and \ref{Co:essemiperf-rel-plain} the functor is equivalent to the functor of window categories induced by  the frame homomorphism $\usW(R'/R)\to\usW(R)$. This frame homomorphism satisfies the hypotheses of Lemma \ref{Le:deformation-equivaelence-windows} because $\hat W(\Fa)=\ker(\sW(R')\to\sW(R))$ is isomorphic to $\bigoplus_{\NN}\Fa$ such that the $F$-linear endomorphism $\dot\sigma$ of $\hat W(\Fa)$ corresponds to the shift on $\bigoplus_{\NN}\Fa$ given by $(a_0,a_1,\ldots)\mapsto(a_1,a_2,\ldots)$, which is locally nilpotent.
\eproof

\bRe
The frame homomorphism $\usW(R'/R)\to\u W(R'/R)$ of Remark \ref{Re:Hom-usW-uW-rel} gives a functor $\sDisp(R'/R)\to\Disp(R'/R)$ similar to Remark \ref{Re:sDisp-Disp}.
\eRe

\subsection{Rim-Schlessinger condition}

A Cartesian diagram of rings
\beqn
\label{Eq:sq-zero-fib-prod}
\xymatrix@M+0.2em{
B' \ar[r] & B \\
A' \ar[r] \ar[u] & A \ar[u]
}
\eeqn
where $B'\to B$ is surjective with kernel of square zero will be called a square zero fiber product. For a base ring $R$ let $\Aff'_R\subseteq\Aff_R$ be a full subcategory which is stable under square zero fiber products. A fibered category $\XXX\to\Aff'_R$ satisfies the strong Rim--Schlessinger condition (RS*) if $\XXX(A')\to\XXX(A)\times_{\XXX(B)}\XXX(B')$ is an equivalence for every square zero fiber product in $\Aff'_R$.

\bEx
\label{Ex:RS*-algebraic}
Every algebraic stack $\XXX\to\Aff_R$ satisfies (RS*) by \cite[\href{https://stacks.math.columbia.edu/tag/0CXP}{Lemma 0CXP}]{Stacks}. 
\eEx

\bEx
\label{Ex:RS*-LF}
The stack $\LF\to\Aff$ of finite projective modules satisies (RS*) by direct verification, or by the more general \cite[Theorem 2.2]{Milnor:Algebraic-K-Theory}, or using that the groupoid core of $\LF$ is $\bigsqcup_{n\ge 0}B{\GL_n}$, which allows to apply Example \ref{Ex:RS*-algebraic}.
\eEx

\bEx
\label{Ex:RS*-BT}
The stack $\BT$ of $p$-divisible groups satisfies (RS*) as a consequence of Example \ref{Ex:RS*-LF} or of Example \ref{Ex:RS*-algebraic} since $\BT=\lim_n\BT_n$ where the groupoid core of $\BT_n$ is algebraic.
\eEx

\bLe
\label{Le:RS*-LFsW}
The fibered category $\LF(\sW(-))\to\Aff^{\essemiperf}$ of finite projective $\sW(-)$-modules over essentially semiperfect rings satisfies (RS*). 
\eLe

\bproof
For a square zero fiber product \eqref{Eq:sq-zero-fib-prod} of essentially semiperfect rings let $I=\ker(B'\to B)=\ker(A'\to A)$. By Proposition \ref{Pr:WN-sWR-sWRN} there is a Cartesian diagram of rings
\[
\xymatrix@M+0.2em{
\sW(B') \ar[r] & \sW(B) \\
\sW(A') \ar[r] \ar[u] & \sW(A) \ar[u]
}
\]
where the horizontal arrows are surjective with kernel $\hat W(I)$, which is an ideal os square zero again. The lemma follows since the fibered category $\LF$ satisfies (RS*).
\eproof

\bLe
\label{Le:RS*-sDisp}
The fibered category $\sDisp^{\essemiperf}\to\Aff^{\essemiperf}$ of sheared displays over essentially semiperfect rings satisfies (RS*).
\eLe

\bproof
For an essentially semiperfect ring $R$ we have $\sDisp(R)=\sDisp^{\plain}(R)=\Win(\usW(R))$ by Corollary \ref{Co:essemiperf-plain}. Since the fibered categories $\LF$ and $\LF(\sW(-))$ satisfy (RS*) by Example \ref{Ex:RS*-LF} and Lemma \ref{Le:RS*-LFsW}, the same follows for finite projective filtered modules over $\usW(-)$ (using Corollary \ref{Co:sWR-R-surj} and Remark \ref{Re:fil-proj-mod-ideal}).
For a given finite projective filtered $\sW(A')$-module $M$ the sets of display structures on $M$ over the four rings of \eqref{Eq:sq-zero-fib-prod} form a Cartesian square because a display structure is an isomorphism of certain finite projective $\sW(-)$-modules.
\eproof

\cclearpage
\section{The Tate module functor}
\label{Se:Tate}

In this section, we construct the functor $\sFZ$ from sheared displays to Tate module schemes, the projective limit version of $p$-divisible groups. The definition is explicit (\S\ref{Se:Tate:constr}), but the proof that it gives a Tate module scheme (Theorem \ref{Th:sC-gives-BT}) occupies most of the section. In characteristic $p$, this follows from a comparison with the truncated Zink functor and its dual version in \S\ref{Se:Tate:char-p}. The general case is established by a study of deformations in \S\S\ref{Se:Tate:deform-seq}--\ref{Se:Tate:first-order-nbhd}, using that $p$-inverted Tate module schemes deform uniquely.

\subsection{Tate module schemes}

Let $R$ be a ring.

\bDe
A Tate module scheme over $R$ is a flat affine commutative group scheme $T$ over $R$ such that $p\colon T\to T$ is injective, the fpqc quotient sheaf $T/pT$ is representable by a finite locally free group scheme $T_1$, which implies that $T/p^nT$ is representable by a finite locally free group scheme $T_n$, and such that $T\to\lim_n T_n$ is an isomorphism.
\eDe

\bRe
The definition can be formulated without fpqc sheafification, thus avoiding set-theoretic questions. One requires the existence of a finite locally free group scheme $T_1$ and a homomorphism $T\to T_1$ such that $0\to T\xrightarrow pT\to T_1\to 0$ is fpqc exact.
\eRe

\bRe
If $T=\Spec A$ is a Tate module scheme, then $A$ is an $\NN$-colimit of finite syntomic $R$-algebras with syntomic transition maps. In particular, $T$ is quasisyntomic, and the quotient $T/p^nT$ can be formed in the category of sheaves for the countably syntomic topology (Definition \ref{De:N-syntomic-topology}). One could use this topology in the definition.
\eRe

\bLe
The category $\BT(R)$ is equivalent to the category of Tate module schemes by sending $G$ to $T=\lim_n G[p^n]$ with transition maps $p\colon G[p^{n+1}]\to G[p^n]$.
\eLe

\bproof The inverse functor sends $T$ to $\colim_n T_n$ with transition maps $p\colon T_n\to T_{n+1}$.
\eproof

\bLe 
If $R$ is an $\FF_p$-algebra, then every Tate module scheme $T$ over $R$ is relatively semiperfect, i.e.\ the relative Frobenius $F\colon T\to T^{(p)}$ is a closed immersion.
\eLe

\bproof
The composition $T\xrightarrow FT^{(p)}\xrightarrow V T$ is multiplication by $p$, which is a closed immersion.
\eproof

\subsection{Deformations of Tate module schemes}

\bLe
\label{Le:deform-quotient}
Let $R'=R/N$ for a nilpotent ideal $N$ and let $H\to G$ be a homomorphism of flat affine commutative group schemes over $R$ which is a closed immersion, with reduction $H'\to G'$ over $R'$. If the fqpc quotient $G'/H'$ is representable by a flat affine group scheme over $R'$ then $G/H$ is representable by a flat affine group scheme over $R$.
\eLe

\bproof
Assume that $N^2=0$.
Let $G=\Spec A$ and $G'=\Spec A'$. 
Let $B=A^H$ and $B'=A'^H$ be the rings of $H$-invariants in $A$ and in $A'$.
Then $X'=\Spec B'$ is a flat group scheme over $R'$ that represents $G'/H'$, and $X=\Spec B$ is (at least) an affine scheme over $R$ with an $H$-invariant morphism $G\to X$.
The exact sequence $0\to A'\otimes_{R'}N\to A\to A'\to 0$ gives
\[
0\to H^0(H,A'\otimes_{R'}N)\to B\to B'\to H^1(H,A'\otimes_{R'}N).
\]
 Here 
the algebraic group cohomology $H^i(H,A'\otimes_{R'}N)$ is the cohomology of the \v Cech complex $\check C^*(G'/X',\OOO_{X'}\otimes_{R'}N)$, thus $H^0=B'\otimes_{R'}N$ and $H^1=0$ by faithfully flat descent. This gives the lower line in the following commutative diagram with exact rows; 
the upper line is clear.
\[
\xymatrix@M+0.2em{
& B\otimes_{R}N \ar[r] \ar[d]^u & B \ar[r] \ar@{=}[d] & B\otimes_RR' \ar[r] \ar[d] & 0 \\
0 \ar[r] & B'\otimes_{R'}N \ar[r] & B \ar[r] & B' \ar[r] & 0 
}
\]
Here $u$ is surjective, thus all vertical arrows are bijective, 
in particular $B\otimes_RN=NB$. Hence $B$ is flat over $R$ since $B'$ is flat over $R'$, using \cite[Theorem 22.3]{Matsumura:Commutative}. Similarly, $B\to A$ is faithfully flat because $B'\to A'=A\otimes_BB'$ is faithfully flat and $\Tor_1^B(B',A)=\Tor_1^R(R',A)$ is zero. The morphism $G\times_{\Spec R} H\to G\times_XG$ of flat affine $R$-schemes defined by projection and multiplication is an isomorphism because this holds over $R'$. Hence $X$ represents $G/H$.
\eproof

\bLe
\label{Le:deform-Tate-module-scheme}
Let $R'=R/N$ for a nilpotent ideal $N$, and let $T$ be a flat affine commutative group scheme over $R$. If $T'=T_{R'}$ is a Tate module scheme, then so is $T$.
\eLe

\bproof
Since $p\colon T'\to T'$ is a closed immersion, the same holds for $p\colon T\to T$. Since $T'/pT'$ is representable by a finite locally free group scheme $T'_1$, Lemma \ref{Le:deform-quotient} implies that $T/pT$ is representable by an affine flat group scheme $T_1$, necessarily finite locally free.
Then $T/p^nT$ is representable by a finite locally free group scheme $T_n$. The homomorphism $T\to\lim T_n$ of flat affine group schemes over $R$ is an isomorphism because this holds over $R'$.
\eproof

\bLe[Deformation sequence]
\label{Le:deform-sequence-TpG}
Let $R$ be a $p$-nilpotent ring and $R'=R/N$ for a nilpotent ideal $N$, 
and let $i\colon \Spec R'\to \Spec R$ be the natural morphism.
For $G\in\BT(R)$ with image $G'\in\BT(R')$ there is an exact sequence
\beqn
\label{Eq:deform-sequence-TpG}
0\to T_pG\to i_*(T_pG')\to \hat G(\u N)\to 0
\eeqn
of fpqc sheaves on $\Spec R$ where $\hat G$ is the formal completion of $G$, and $\u N(S)=NS$ for an $R$-algebra $S$. The sequence is exact on the level of countably syntomic sheaves.
\eLe

\bproof
Since  $G$ is formally smooth by \cite[Chap.~I, Theorem 3.3.13]{Messing:Crystals}, there is an exact sequence of presheaves $0\to\hat G(\u N)\to G\to i_*G'\to 0$, and moreover, $\hat G(\u N)$ is annihilated by a power of $p$. Let $V_pG=\colim(T_pG,p)=\lim(G,p)$ as an fpqc sheaf. It follows that $V_pG\to i_*V_pG'$ is an isomorphism.
We have a commutative diagram of fpqc sheaves whose rows are exact in the category of countably syntomic sheaves, which gives the exact sequence \eqref{Eq:deform-sequence-TpG}.
\beqn
\label{Eq:deform-sequence-TpG-diag}
\xymatrix@M+0.2em{
0 \ar[r] & T_pG \ar[d] \ar[r] & V_pG \ar[d]^\cong \ar[r]^-\pi & G \ar@{->>}[d] \ar[r] & 0 \\
0 \ar[r] & i_*T_pG' \ar[r] & i_*V_pG' \ar[r] & i_*G' \ar[r] & 0 
}
\eeqn
Indeed, $\pi$ is surjective since $p\colon G\to G$ is a surjective syntomic morphism, and the rest is automatic (one could also use that $i_*$ is exact by Lemma \ref{Le:N-syntomic-cocont}). 
\eproof

\bRe
In Lemma \ref{Le:deform-sequence-TpG}, if the ideal $N$ has square zero then $\hat G(\u N)=\Lie(G')\otimes_{R'}\u N$.
\eRe

\subsection{Construction of the Tate module functor}
\label{Se:Tate:constr}

Let $R$ be a $p$-nilpotent ring.

\bDe
\label{De:sC-T}
For a plain sheared display $\uM$ over $R$ (Definition \ref{De:Sh-disp}) let $\tilde{\uM}$ be the associated fpqc sheaf of plain sheared displays over $\Spec R$ defined by base change and let
\[
\sC(\uM)=\Gamma(\tilde{\uM})=[\tilde M_1\xrightarrow{\phi_1-\tau}\tilde M_0]
\] 
as a complex of fpqc sheaves of abelian groups over $\Spec R$ sitting in degrees $[0,1]$, using Definition \ref{De:Gamma(uM)}.
By $p^\infty$-root descent this construction extends to an additive functor
\beqn
\sC\colon \sDisp(R)\to\Ch^{[0,1]}(\Ab(R_{\fpqc}))
\eeqn
where $\Ab(R_{\fpqc})$ denotes the category of fpqc sheaves of abelian groups on $\Aff_R$. We define
\beqn
T\colon \sDisp(R)\to\Ab(R_{\fpqc}),\qquad
T(\uM)=H^0(\sC(\uM)).
\eeqn
We will also consider the truncated version $\sC_n(\uM)=\cone(p^n:\sC(\uM)\to\sC(\uM))$.
\eDe

For a complex $C^*$ of fpqc sheaves of abelian groups over $\Spec R$ we will denote by $H^i(C^*)_\tau$ the cohomology of $C^*$ in the category of $\tau$-sheaves of abelian groups.

\bLe
\label{Le:sC-gives-BT}
Let either $\uM\in\sDisp(R)$ and $\tau=\fpqc$, or $\uM\in\sDisp^{\plain}(R)$ and $\tau=\csyn$;
see Definition \ref{De:N-syntomic-topology}.
 The following are equivalent.
\begin{enumerate}
\item
\label{It:T-is-Tate}
$T(M)$ is representable by a Tate module scheme over $R$, and $H^1(\sC(\uM))_{\tau}$ is zero.
\item
\label{It:T-is-Tate-n}
For every $n$, the sheaf $H^0(\sC_n(\uM))_{\tau}$ is representable by an $n$-truncated $BT$ group over $R$, and $H^i(\sC_n(\uM))_{\tau}$ is zero for $i\ne 0$.
\item
\label{It:T-is-Tate-1}
For $n=1$, the sheaf
$H^0(\sC_1(\uM))_{\tau}$ is representable by a finite locally free group scheme over $R$, and $H^i(\sC_1(\uM))_{\tau}$ is zero for $i\ne 0$.
\end{enumerate}
\eLe

\bproof
\eqref{It:T-is-Tate} $\Rightarrow$ \eqref{It:T-is-Tate-n}
follows from the exact triangle $\sC(\uM)\xrightarrow{p^n} \sC(\uM)\to\sC_n(\uM)$ of presheaves. \eqref{It:T-is-Tate-n}  $\Rightarrow$ \eqref{It:T-is-Tate-1} is clear. Assume \eqref{It:T-is-Tate-1} holds. The exact triangles $\sC_n(\uM)\to\sC_{n+m}(\uM)\to\sC_m(\uM)$ of presheaves show that \eqref{It:T-is-Tate-n} holds, so $T_n=H^0(\sC_n(\uM))_{\tau}$ is an $n$-truncated BT group, and they show that the natural homomorphism $T_{n+1}\to T_n$ is surjective as $\tau$-sheaves and identifies $T_{n+1}[p^n]$ with $T_n$. Hence $T=\lim_nT_n$ is a Tate module scheme. The components of $\sC(\uM)$ are sectionwise derived $p$-complete by Corollary \ref{Co:sW-p-complete}. Thus $\sC(\uM)$ is derived $p$-complete, using that $\tau$-sheaves on $\Aff$ form a replete topos. 
Hence $\sC_n\cong R\lim_nT_n=T$ because $R^1\lim_nT_n=0$ since the $T_n$ form a surjective system.  See \cite[\S B.2]{{Drinfeld:Ring-stacks-conjecturally}} for derived $p$-complete sheaves.
\eproof

The goal of this section is to prove that the equivalent conditions of Lemma \ref{Le:sC-gives-BT} are always true; see Theorem \ref{Th:sC-gives-BT}. If this is known over a certain ring $R$, we obtain an exact functor
\beqn
\label{Eq:sFZ}
\sFZ_R\colon \sDisp(R)\to\BT(R),\qquad T_p(\sFZ_R(\uM))=T(\uM).
\eeqn

\subsection{The Tate module functor in characteristic $p$} 
\label{Se:Tate:char-p}

Let $R$ be an $\FF_p$-algebra. 

\bRe[$n$-smooth-cosmooth decomposition]
For every $n$-truncated BT group over $R$ there is a canonical exact sequence
\[
0\to G[F^n]\to G\to (G^\vee[F^n])^\vee\to 0
\]
where $(G^\vee[F^n])^\vee=G^{(p^n)}[V^n]$, which is uniquely determined by the property that the first group is $n$-smooth (\S\ref{Se:Truncated-Zink-functor}) and the third group is $n$-cosmooth (\S\ref{Se:Dual-Zink-functor}).
\eRe

\bLe
\label{Le:fiber-sequence-sCn-Cn}
For $\uM\in\sDisp(R)$ with image $\uM_n\in\Disp_n(R)$ there is an exact triangle in the derived category of fpqc sheaves of abelian groups (see \eqref{Eq:ZnM} and Definition \ref{De:CnM} for notation)
\[
Z_n(\uM_n)\to \sC_n(\uM)\to C_n(\uM_n)\to Z_n(\uM_n)[1]
\]
If $\uM$ is plain, this exact triangle exists in the derived category of countably syntomic sheaves.
\eLe

\bproof
We assume that $\uM$ is plain, i.e.\ a window over $\usW(R)$.
The general case follows by descent of the construction.
There is an exact sequence of complexes of sheaves of $\sW$-modules in degrees $[-1,0]$
\[
0\longrightarrow
[\sW\xrightarrow{F^n}F^n_*(\sW)]\xrightarrow{(\id,V^n)}[\sW\xrightarrow{\;p^n\;}\sW]\longrightarrow[0\to W_n]
\longrightarrow0
\]
which is exact as countably syntomic sheaves by Corollary \ref{Co:sW-Vn-sW-Wn}, and the first complex is quasi-isomorphic to $\hat W[F^n][1]$ as syntomic sheaves by Lemma \ref{Le:hatWFm-sW-sW}.
There is a similar exact sequence of filtered modules over $\usW$ with $\usW$ and $\u{W\!}_n$ in place of $\sW$ and $W_n$. The functor $\uM\otimes_{\usW(R)}-$ preserves this exact sequence and the quasi-isomorphisms, and the lemma follows from the definitions of $Z_n$ and $C_n$.
\eproof

\bPr
\label{Pr:sC-gives-BT-FFp}
If $R$ is an $\FF_p$-algebra, the equivalent conditions of Lemma \ref{Le:sC-gives-BT} hold, and hence the functor $\sFZ_R$ of \eqref{Eq:sFZ} is defined. This functor preserves height and dimension.
\ePr

\bproof
We can assume that $\uM$ is plain.
Let $\uM_n\in\Disp_n(R)$ be the image of $\uM$.
We work with the countably syntomic topology.
By Lemma \ref{Le:fiber-sequence-sCn-Cn}, \S \ref{Se:Truncated-Zink-functor}, and Proposition \ref{Pr:dual-truncated-Zink-functor}, the complex of sheaves $\sC_n(\uM)$ is acyclic outside degree zero, and $G_n=H^0(\sC_n(\uM))$ sits in an exact sequence of sheaves
\beqn
\label{Eq:FZn-Gn-FYn}
0\to\FZ_n(\uM_n)\to G_n\to\FY_n(\uM_n)\to 0
\eeqn
where the outer term are commutatative finite locally free group schemes, so the same holds for $G_n$. 
In particular, Condition \eqref{It:T-is-Tate-1} of Lemma \ref{Le:sC-gives-BT} is satisfied. 
Hence $G_n$ is an $n$-truncated BT groups by Condition \eqref{It:T-is-Tate-n} of Lemma \ref{Le:sC-gives-BT}.
The exact sequence \eqref{Eq:FZn-Gn-FYn} is the $n$-smooth-cosmooth decomposition of $G_n$ by the uniqueness of that decomposition. Hence \S\ref{Se:Truncated-Zink-functor} and Proposition \ref{Pr:dual-truncated-Zink-functor} give $\dim(G)=\rk(\Lie\FZ_n(\uM_n))=\dim(\uM_n)$ and 
\[
\height(G)=\rk(\Lie\FZ_n(\uM_n))+\rk(\Lie\FY_n(\uM_n))=\dim(\uM_n)+\dim(\uM_n^t)=\height(\uM_n).
\qedhere
\]
\eproof

\subsection{Deformation sequence} 
\label{Se:Tate:deform-seq}

Let $R$ be a $p$-nilpotent ring. The composition of functors
\[
\sDisp(R)\to\Disp(R)\xrightarrow{\FZ_R}\FG(R)
\] 
will be denoted by $\FZ'_R$. The Zink functor $\FZ_R$ is recalled in \S\ref{Se:Recall-Zink-functor}.

\bLe[Deformation sequence]
\label{Le:deform-sequence-TM}
Let $R$ be a $p$-nilpotent ring and $R'=R/N$ for a nilpotent ideal $N$. Let $\uM\in\sDisp^{\plain}(R)$ with image $\uM'\in\sDisp^{\plain}(R')$. Let $H=\FZ'_R(\uM)$. Then there is an exact sequence of countably syntomic sheaves on $\Spec R$ 
\beqn
\label{Eq:deform-sequence-T(M)}
0\to T(\uM)\to i_*T(\uM')\to H(\u N)\to H^1(\sC(\uM))_{\csyn}\to i_*H^1(\sC(\uM))_{\csyn}\to 0
\eeqn
where $i\colon \Spec R'\to \Spec R$ is the natural morphism and $\u N(S)=NS$ for an $R$-algebra $S$. 
\eLe

\bproof
Let $\uM^\wedge\in\Disp(R)$ be the image of $\uM$, thus $\FZ'_R(\uM)=\FZ_R(\uM^\wedge)$. For a nilpotent $R$-algebra $N'$ we have $H(N')=H^0(Z(\uM^\wedge)(N'))$ where 
\[
Z(\uM^\wedge)(N')=\Gamma(\uM^\wedge\otimes_{\u W(R)}\u{\hat W}(N')) [1] 
\cong\Gamma(\uM\otimes_{\usW(R)}\u{\hat W}(N')) [1];
\]
the equality is \eqref{Eq:ZMN} and the isomorphism involves multiplication by $\alpha^{-1}$ with $\alpha$ as in Remark \ref{Re:Hom-usW-uW}. Hence by Proposition \ref{Pr:WN-sWR-sWRN} and Corollary \ref{Co:sWR-sWR/N-Nsynt} we get an exact sequence of complexes of fpqc sheaves on $\Spec R$ concentrated in degrees $[0,1]$ 
\beqn
\label{Eq:ZMN-sC-isC}
0\to Z(\uM^\wedge)(\u N)[-1]\to \sC(\uM)\to i_*(\sC(\uM'))\to 0
\eeqn
which is exact as countably syntomic sheaves.
Now \eqref{Eq:deform-sequence-T(M)} is the cohomology sequence of \eqref{Eq:ZMN-sC-isC}, using that $i_*$ is exact for the countably syntomic topology by Lemma \ref{Le:N-syntomic-cocont}.
\eproof

\bRe
The deformation sequence \eqref{Eq:deform-sequence-T(M)} is used in the proof of Theorem \ref{Th:sC-gives-BT}, which then implies that \eqref{Eq:deform-sequence-T(M)} is actually a short exact sequence $0\to T(\uM)\to i_*T(\uM')\to H(\u N)\to 0$. 
\eRe

\bLe
[Comparison of deformation sequences]
\label{Le:Compare-T(M)-TpG}
Let $\uM\in\sDisp^{\plain}(R)$ and assume that $T(\uM)\cong T_pG$ with $G\in\BT(R)$.  Let $H=\FZ'_R(\uM)$. 
Then there is a homomorphism of formal groups $u\colon \hat G\to H$ with the property that for every nilpotent thickening of $R$-algebras $S\to S'=S/N$ the following diagram of fpqc sheaves on $\Spec S$ commutes 
\beqn
\label{Eq:Compare-T(M)-TpG}
\xymatrix@M+0.2em{
0 \ar[r] & T_pG_S \ar[r] \ar[d]^\cong & i_*(T_pG_{S'}) \ar[r] \ar[d]^\cong & \hat G(\u N) \ar[r] \ar[d]^u & 0 \\
0 \ar[r] & T(\uM_S) \ar[r] & i_*(T(\uM_{S'})) \ar[r] & H(\u N) 
}
\eeqn
where $i\colon \Spec S'\to\Spec S$ is the natural morphism, we set $\u N(B)=NB$ for an $S$-algebra $B$ as earlier, 
the upper row is \eqref{Eq:deform-sequence-TpG}, and the lower row is \eqref{Eq:deform-sequence-T(M)}.
Moreover, the homomorphism $u$ is an isomorphism.
\eLe

\bproof
For a given pair $(S,N)$ where $S$ is an $R$-algebra and $N\subseteq S$ is a nilpotent ideal, the 
rows of \eqref{Eq:Compare-T(M)-TpG} are exact by 
Lemmas \ref{Le:deform-sequence-TpG} and \ref{Le:deform-sequence-TM}, 
and the isomorphism $T_pG\xrightarrow\sim T(\uM)$ gives an injective homomorphism $u_{S,N}\colon \hat G(\u N)\to H(\u N)$ of fpqc sheaves over $\Spec S$, which is functorial in the pair $(S,N)$. For every nilpotent $R$-algebra $N$, using the pair $(R\oplus N,N)$ and taking global sections gives an injective homomorphism of abelian groups $u_N\colon \hat G(N)\to H(N)$ which is functorial in $N$, thus an injective homomorphism of formal groups $u\colon \hat G\to H$. The homomorphism $u$ gives back the collection of homomorphisms $u_{S,N}$ since for each pair $(S,N)$ there is a homomorphism of pairs $(R\oplus N,N)\to(S,N)$ which is the identity on $N$.
Now $\hat G$ and $H$ are formal groups of the same dimension because this can be verified over the reduction to $R/p$, 
where Proposition \ref{Pr:sC-gives-BT-FFp} applies.
Hence the monomorphism $u$ is an isomorphism. Indeed, $u$ induces a monomorphism of the Lie algebras, and a homomorphism of projective modules over a ring $C$ which gives a monomorphism of functors on $C$-algebras is a split injection.
\eproof

\subsection{First order neighborhoods of Tate module schemes}
\label{Se:Tate:first-order-nbhd}

\bCon[Absolute first order neighborhood of $T_pG$]
\label{Con:inf-nbhd-TpG-abs}
Let $R$ be a $p$-nilpotent ring and $G\in\BT(R)$. We denote by $(T_pG)_1$ the first order (square zero) infinitesimal neighborhood of $T_pG$ in $V_pG=\colim(G,p)$.
This can be made explicit as follows. The first order infinitesimal neighborhood of zero in $G$ is isomorphic to $\Spec(R\oplus\Lie(G)^\vee)$, and we have a commutative diagram of fpqc sheaves with Cartesian squares and affine faithfully flat vertical morphisms
\beqn
\xymatrix@M+0.2em{
T_pG \ar[r] \ar[d] & (T_pG)_1 \ar[r]^-j \ar[d] & V_pG \ar[d] \\
\Spec R \ar[r] & \Spec(R\oplus\Lie(G)^\vee) \ar[r]^-{j_0} & G.
}
\eeqn
Hence $(T_pG)_1$ is an affine scheme. If $T_pG=\Spec A$ and $(T_pG)_1=\Spec B$ then $A=B/\Fa$ with $\Fa=\Lie(G)^\vee\otimes_RA$ as an $A$-module, in particular $A,B,\Fa$ are flat over $R$. 
The deformation sequence \eqref{Eq:deform-sequence-TpG} for $B\to A$ in place of $R\to R'$ gives an exact sequence of abelian groups
\beqn
\label{Eq:def-sequence-inf-nbhd-abs}
0\to T_pG(B)\to T_pG(A)\xrightarrow{\eta}\hat G(\Fa)\cong\Lie G\otimes_R\Fa\cong\End_A(\Fa).
\eeqn
\eCon

\bLe
\label{Le:inf-nbhd-image-idA}
In \eqref{Eq:def-sequence-inf-nbhd-abs} we have $\eta(\id_A)=\id_{\Fa}$.
\eLe

\bproof
By the proof of Lemma \ref{Le:deform-sequence-TpG} the homomorphism $\eta$ is the composition
\[
T_pG(A)\to V_pG(A)\xleftarrow\sim V_pG(B)\to G(B).
\]
The image of $\id_A\in T_pG(A)$ in $V_pG(A)$ lifts to the canonical element $j$ of $V_pG(B)$, whose image in $G(B)$ comes from the canonical element $j_0$ of $G(\Spec R\oplus\Lie(G)^\vee)$.
\eproof

\bCon
[Relative first order neighborhood of $T_pG$]
\label{Con:inf-nbhd-TpG-rel}
Let $R\to\bar R=R/N$ where $N$ has square zero and let $\bar G\in\BT(\bar R)$. 
Let $V=V_pG$ for some lift $G\in\BT(R)$ of $\bar G$. Then $V$ does not depend on the choice of $G$. 
We denote by $(T_p\bar G)_{1,R}$ the square zero neighborhood of $T_p\bar G$ in $V$. This is also the square zero neighborhood of $T_p\bar G$ in $(T_pG)_{1}$ and hence an affine scheme.
\eCon

\bLe
\label{Le:inf-nbhd-kernels}
In the situation of Construction \ref{Con:inf-nbhd-TpG-rel} let $T_p\bar G=\Spec\bar A$, $(T_p\bar G)_1=\Spec\bar B$, $(T_p\bar G)_{1,R}=\Spec C$, and $\bar A=\bar B/\bar\Fa=C/\Fb$. Then there is an exact sequence of $C$-modules
\[
0\to N\otimes_R\bar A\to\Fb\to\bar\Fa\to 0.
\]
\eLe

\bproof
Clearly $NC\to\Fb\to\bar\Fa\to 0$ is exact.
Let $T_pG=\Spec A$ and $(T_pG)_1=\Spec B$ and $A=B/\Fa$. Then $C=B/(NB+\Fa)^2=B/N\Fa$ and thus $NC=NB/N\Fa=N\otimes_{\bar R}\bar A$, using that $B$ is flat over $R$ for the last equality.
\eproof

\bPr
\label{Pr:TM-is-Tate-deforms}
Let $R\to\bar R=R/N$ where $N$ has square zero and let $\uM\in\sDisp^{\plain}(R)$ with image $\u{\bar M}\in\sDisp^{\plain}(\bar R)$. If $T(\u{\bar M})$ is a Tate module scheme, then so is $T(\uM)$.
\ePr

\bproof
Let $T(\u{\bar M})\cong T_p\bar G=\Spec\bar A$ with $G\in\BT(\bar R)$, let $(T_p\bar G)_1=\Spec\bar B$ be the first infinitesimal neighborhood of $T_p{\bar G}$ as in Construction \ref{Con:inf-nbhd-TpG-abs}, and let $(T_p\bar G)_{1,R}=\Spec C$ be the relative first infinitesimal neighborhood of $T_p\bar G$ as in Construction \ref{Con:inf-nbhd-TpG-rel}. The deformation sequences of Lemma \ref{Le:deform-sequence-TM} for $C\to\bar A$ and for $\bar B\to\bar A$ together with the comparison isomorphism of Lemma \ref{Le:Compare-T(M)-TpG} give the following commutative diagram with exact rows.
\[
\xymatrix@M+0.2em{
0 \ar[r] & T(\uM)(C) \ar[r] \ar[d] & T_p\bar G(\bar A) \ar[r] \ar@{=}[d] & \Lie(\bar G)\otimes_{\bar R}\Fb \ar[d] \cong\Hom_{\bar A}(\bar\Fa,\Fb) \hspace{-6.3em} \\
0 \ar[r] & T(\uM)(\bar B) \ar[r] & T_p\bar G(\bar A) \ar[r] & \Lie(\bar G)\otimes_{\bar R}\bar\Fa 
\cong\Hom_{\bar A}(\bar\Fa,\bar\Fa)\hspace{-6.3em}
}\hspace{6.3em}
\]
Let $s\in\Hom_{\bar A}(\bar\Fa,\Fb)$ be the image of $\id_{\bar A}\in T_p\bar G(\bar A)$ in the upper line, and let $A=C/s(\bar\Fa)$. Then $s$ is a section of the projection $\Fb\to\bar\Fa$ by Lemma \ref{Le:inf-nbhd-image-idA}. Using Lemma \ref{Le:inf-nbhd-kernels} we get an exact sequence $0\to N\otimes_{\bar R}\bar A\to A\to\bar A\to 0$, which implies that $A$ is flat over $R$ with $A\otimes_R\bar R=\bar A$. Now the deformation sequence \eqref{Eq:deform-sequence-T(M)} for $A\to\bar A$ reads
\[
0\longrightarrow T(\uM)(A) \longrightarrow T_p\bar G(\bar A)\longrightarrow \Hom_{\bar A}(\bar\Fa,\bar\Fb/s(\bar\Fa)).
\]
Here $\id_{\bar A}$ maps to zero by the construction and hence lifts to an element $g\in T(\uM)(\bar A$). 

We claim that $g\colon \Spec A\to T(\uM)$ is an isomorphism. 
If this holds, $\Spec A$ is a flat affine group scheme over $R$ whose reduction over $\bar R$ is a Tate module scheme, hence $\Spec A$ is a Tate module scheme by Lemma \ref{Le:deform-Tate-module-scheme}, which finishes the proof.

To prove the claim let $i\colon \Spec\bar R\to\Spec R$.
The morphism $g$ is injective because the composition $\Spec A\to T(\uM)\to i_*T_p\bar G\to i_*V_p\bar G$ factors into $\Spec A\to V_pG\xrightarrow\sim i_*V_p{\bar G}$ where the first arrow is injective.
To prove surjectivity, let $S$ be an $R$-algebra and $\bar S=S/NS$. We have to show that for each $h\in T(\uM)(S)$ the image $\bar h\in T_p\bar G(\bar S)=(\Spec\bar A)(\bar S)$ lifts to an element of $(\Spec A)(S)$. As earlier let $G\in\BT(R)$ be a lift of $\bar G$. Clearly $\bar h$ lifts to an element $h'\in V_pG(S)$, which necessarily lies in $(\Spec C)(S)$. The resulting ring homomorphism $u\colon C\to S$ induces a commutative diagram of deformation sequences
\[
\xymatrix@M+0.2em{
0 \ar[r] & T(\uM)(C) \ar[r] \ar[d] & T_p\bar G(\bar A) \ar[r] \ar[d] & \Hom_{\bar A}(\bar\Fa,\Fb) \ar[d] \\
0 \ar[r] & T(\uM)(S) \ar[r] & T_p\bar G(\bar S) \ar[r] & \Hom_{\bar A}(\bar\Fa,NS).
}
\]
Here $\id_{\bar A}\in T_p\bar G(\bar A)$ maps to $\bar h\in T_p\bar G(\bar S)$, which comes from $h\in T(\uM)(S)$. It follows that the composition $\bar\Fa\xrightarrow s\Fb\xrightarrow u NS$ is zero, so $u$ factors as $C\to A\to S$, in other words $h'$ lies in $(\Spec A)(S)$ as desired.
\eproof

\subsection{Conclusions}
\label{Se:Tate-conclusions}

\bTh
\label{Th:sC-gives-BT}
The equivalent conditions of Lemma \ref{Le:sC-gives-BT} are always true, and hence the functor $\sFZ_R$ of \eqref{Eq:sFZ} exists for every $R$. This functor preserves height and dimension.
\eTh

\bproof 
The case of $\FF_p$-algebras holds by Proposition \ref{Pr:sC-gives-BT-FFp}. In general, by induction we can assume that there is an ideal of square zero $N\subseteq R$ such that the theorem holds over $R/N$, and we can assume that the given sheared display $\uM$ over $R$ is plain. Then Proposition \ref{Pr:TM-is-Tate-deforms} implies that $T(\uM)$ is a Tate module scheme. The deformation sequence of Lemma \ref{Le:deform-sequence-TM} reduces to $0\to T(\uM)\to i_*T(\uM')\xrightarrow v H(\u N)\to H^1(\sC(\uM))_{\csyn}\to 0$. Here the homomorphism $v$ is surjective by Lemma \ref{Le:Compare-T(M)-TpG}, hence $H^1(\sC(\uM))_{\csyn}$ vanishes. Height and dimension over $R$ are determined over $R/p$.
\eproof

\bPr
\label{Pr:FZ-sFZ}
There is a commutative diagram of fibered categories over $\Aff_{\Spf\ZZ_p}$:
\[
\xymatrix@M+0.2em{
\sDisp \ar[r]^{\sFZ} \ar[d] & \BT \ar[d]^{G\mapsto\hat G} \\
\Disp \ar[r]^{\FZ} & \FG
}
\]
\ePr

\bproof
Immediate from Lemma \ref{Le:Compare-T(M)-TpG}.
\eproof

We recall the windows $\u A$ and $\u A(1)$ over a frame $\u A$ as in Example \ref{Ex:win-uA-uA(1)}.

\bCo
\label{Co:sFZ-usW-Gm}
We have $\sFZ_R(\usW(R))=\mu_{p^\infty}$.
\eCo

\bproof
We know that $G=\sFZ_R(\usW(R))$ is a $p$-divisible group of height $1$ and dimension $1$, in particular $G=\hat G$, and $\hat G=\FZ_R(\u W(R))=\mu_{p^{\infty}}$ by Proposition \ref{Pr:FZ-sFZ} and \cite[(211)]{Zink:The-display}. 
\eproof

\bRe
Corollary \ref{Co:sFZ-usW-Gm} is also proved in \cite[\S 3.6]{Drinfeld:Ring-stacks-conjecturally}.
\eRe

\bPr
\label{Pr:sFZ-usW1-QpZp}
We have $\sFZ_R(\usW(R)(1))=\QQ_p/\ZZ_p$.
\ePr

\bproof
We have $T(\usW(R)(1))=\sW^{F=1}$, the $F$-invariants in $\sW$. Let $T=T_p(\QQ_p/\ZZ_p)$ as a Tate module scheme. Then $T$ is the derived or classical $p$-completion of the constant sheaf $\ZZ$, so there is a natural homomorphism $T\to\sW^{F=1}$. This homomorphism is an isomorphism because the composition $T\to\sW^{F=1}\to W^{F=1}$ is an isomorphism, and we have an exact sequence of sheaves $0\to T_FQ\to\sW\to W\to 0$ where $(T_FQ)^{F=1}$ is zero.
\eproof

\bPr
\label{Pr:sFZ-duality}
The functor $\sFZ_R$ is compatible with duality, i.e.\ there is a natural isomorphism of $p$-divisible groups $\sFZ_R(\uM)^\vee\cong\sFZ_R(\uM^t)$ for $\uM\in\sDisp(R)$.
\ePr

\bproof
Let $\sFZ'_R(\uM)=\sFZ_R(\uM^t)^\vee$.
The natural bilinear form $\uM\times\uM^t\to\usW(R)$ gives a bilinear map $T(\uM)\times T(\uM^t)\to T(\usW(R))=T_p(\mu_{p^\infty})$ using Corollary \ref{Co:sFZ-usW-Gm}, or equivalently a natural homomorphism $c:\sFZ_R(\uM)\to\sFZ'_R(\uM)$ of $p$-divisible groups. We claim that $c$ is an isomorphism. It suffices to prove this over algebraically closed fields. The homomorphism $c$ is an isomorphism when $\uM=\usW$ or $\uM=\usW(1)$ by Proposition \ref{Pr:sFZ-usW1-QpZp}. Hence $c$ is an isomorphism when $\uM$ is ordinary.
Over perfect rings, displays and sheared displays are equivalent.
For a given (sheared) display $\uM$ over a perfect field $k$ there is a deformation to a (sheared) display $\uN$ over $R=k[[t]]\limperf$ with ordinary generic fiber. The homomorphism $\sFZ_R(\uN)\to\sFZ'_R(\uN)$ of $p$-divisible groups over $R$ is an isomorphism in the generic fiber, hence an isomorphism by Lemma \ref{Le:Hom-BT-isom-clopen} below.
\eproof

\bRe
\label{Re:Pr-sFZ-duality-alt-proof}
It should be possible to given alternative proof of Proposition \ref{Pr:sFZ-duality} using Proposition \ref{Pr:FZ-FY-dual} instead of the deformation to the ordinary case. 
\eRe

\bLe
\label{Le:Hom-BT-isom-clopen}
Let $f:G\to H$ be a homomorphism of $p$-divisible groups over a ring $R$. The locus where $f$ is an isomorphism is open and closed.
\eLe

\bproof
Since $f$ is an isomorphism iff $f\colon G[p]\to H[p]$ is an isomorphism, the locus is open, so it suffices to show that it is closed as a set. Thus we can assume that $R$ is perfect. Let $(M,F,V)$ and $(N,F,V)$ be the Dieudonn\'e modules of $G$ and $H$ over $W(R)$ and let $h$ be the common rank. Let $M'=\det(M)$ and $N'=\det(N)$. These are locally free $W(R)$-modules of rank $1$ equipped with $F$ and $V$ such that $FV=VF=p^h$. Necessarily we have $F=p^dF_0$ and $V=p^{h-d}V_0$ with $F_0V_0=\id$ on $M'$ and on $N'$. The pairs $(M',F_0)$ and $(N',F_0)$ are equivalent to \'etale $\ZZ_p$-local systems. Hence the isomorphism locus is an open and closed set.
\eproof

For later use we note the following fact.

\bLe
\label{Le:sFZ-eq-Fp-gen}
If the functor\/ $\sFZ_R$ of \eqref{Eq:sFZ} is an equivalence over all\/ $\FF_p$-algebras then it is an equivalence over all $p$-nilpotent rings.
\eLe

\bproof
By $p^\infty$-root descent it suffices to prove the equivalence over essentially semiperfect rings. This is formal, using that the stack of sheared displays is formally smooth by Lemma \ref{Le:sDisp-form-smooth-essp} and that both sides of the functor satisfy (RS*) over essentially semiperfect rings by Example \ref{Ex:RS*-BT} and Lemma \ref{Le:RS*-sDisp}. More precisely, let $R'=R/N$ for an ideal $N$ of square zero and let $S=R\times_{R'}R$. Then $S\cong R\oplus N$ and hence $S\cong R\times_{R'}(R'\oplus N)$. The condition (RS*) on both sides and formal smoothness of the source implies that for $\uM\in\sDisp(R')$ the functor $\sFZ$ induces an equivalence of the groupoids of deformations to $R$ of $\uM$ and of $\sFZ(\uM)$. It follows that $\sFZ$ is an equivalence of the groupoid cores, hence an equivalence.
\eproof

\subsection{A dual version of the Zink functor}

Let us record the following complement.

\bDe
A $V$-Tate module scheme over an $\FF_p$-algebra $R$ is a flat affine commutative group scheme $Y$ over $R$ such that $V\colon Y^{(p)}\to Y$ is injective with cokernel representable by a finite locally free group scheme $Y_1$, thus the cokernel of $V^n=Y^{(p^n)}\to Y$ is representable by a finite locally free group scheme $Y_n$, and such that $Y\to\lim Y_n$ is an isomorphism.
\eDe

\bRe
Over $\FF_p$-algebras, 
Cartier duality gives a contravariant equivalence between $V$-Tate module schemes and commutative formal groups.
\eRe

We can repeat Definition \ref{De:sC-T} for displays. Let again $R$ be a $p$-nilpotent ring.

\bDe
For a display $\uM$ over $R$ let $C(M)=\Gamma(\tilde{\uM})$ where $\tilde{\uM}$ is the sheaf of displays over $\Spec R$ defined by base change, and let $\YYY(\uM)=H^0(C(\uM))$ as an fpqc sheaf.
\eDe

\bPr
For $\uM\in\Disp(R)$ we have $H^1(C(\uM))_{\fpqc}=0$, and $\YYY(\uM)$ is representable by a flat affine group scheme $Y$ such that $Y_{R/p}$ is a $V$-Tate module scheme.
\ePr

\bproof
The complex $C(\u M)$ is given by a homomorphism $\gamma\colon\tilde M_1\to\tilde M_0$ of pro-unipotent affine group schemes. 
If $R$ is an $\FF_p$-algebra, the proposition is a consequence of Proposition \ref{Pr:dual-truncated-Zink-functor}; in this case $\gamma$ is a countably syntomic cover by the proof of Proposition \ref{Pr:dual-truncated-Zink-functor}. In general it follows that $\gamma$ is flat, hence $T(\uM)$ is a flat affine group scheme, and $H^1(C(\uM))_{\fpqc}=0$.
\eproof

\bRe
For $\uM\in\sDisp(R)$ with image $\uM^\wedge\in\Disp(R)$, there is a natural homomorphism $\sC(\uM)\to C(\uM^\wedge)$, which induces a homomorphism $T(\uM)\to\YYY(\uM^\wedge)$.
If $\uM$ is $F$-nilpotent, this homomorphism is an isomorphism.\footnote{Sketch of proof of this fact. As a consequence of \cite[Corollary 82]{Zink:The-display}, for an $F$-nilpotent sheared display $\uM$, the sheaf $\FY(\uM)$ sits in a deformation sequence similar to \eqref{Eq:deform-sequence-T(M)}. Hence we can assume that $R$ is an $\FF_p$-algebra. In that case we use that for a unipotent $p$-divisible group $G$ we have $T_pG=\lim_n\coker(V^n\colon G^{(p^n)}\to G)$.}
If $R$ is an $\FF_p$-algebra, this homomorphism is the dual of $\widehat{G^\vee}\to G^\vee$ where $G=\sFZ(\uM)$, or equivalently the limit of $G[p^n]\to G[p^n]/V^n$, and there is an exact sequence
\[
0\to \lim G[F^n]\to T_pG\to\YYY(\uM)\to 0
\]
where the limit is taken with respect to $p\colon G[F^{n+1}]\to G[F^n]$.
\eRe

\bRe
Let us call an fpqc sheaf of abelian groups $F$ dualisable if $F\to F^{\vee\vee}$ is an isomorphism, where $F^\vee=\u\Hom(F,\GG_m)$. Then $p$-divisible groups and formal groups are dualisable. For $\uM\in\sDisp(R)$ and $G=\sFZ(\uM)$ we have $G^\vee\cong T(\uM^t)$. Is it true that under this isomorphism, the dual of $\hat G\to G$ is $T(\uM^t)\to\YYY(\uM^t)$?
\eRe

\cclearpage
\section{Sheared display functors}
\label{Se:sh-disp-functor}

The objective of this section is to prove the following result (see \S\ref{Se:Sheared-display-functor:conclusions}).

\bTh
\label{Th:rel-sh-disp-functor}
Let $R'\to R$ be a nilpotent pd thickening of $p$-nilpotent rings. Assume that $R'$ is an $\FF_p$-algebra, or $p\ge 3$. There is an exact additive functor
\beqn
\label{Eq:sDR'R}
\sD_{R'/R}\colon \BT(R)\to\sDisp(R'/R)
\eeqn
which is functorial in $R'/R$,
together with a natural isomorphism of finite projective $R'$-modules 
\beqn
\label{Eq:DDG-sDGrel}
\DD(G)_{R'/R}\cong\sD_{R'/R}(G)\otimes_{\sW(R')} R'
\eeqn
that preserves the Hodge filtration,  
where $\DD$ is the covariant crystalline Dieudonn\'e functor.
\eTh

The category of relative sheared displays $\sDisp(R'/R)$ appears in Definition \ref{De:Sh-disp-rel}.

The proof of Theorem \ref{Th:rel-sh-disp-functor} for $\FF_p$-algebras proceeds as follows. If $R$ is qrsp, we obtain the functor $\sD_{R'/R}$ from the crystalline Dieudonn\'e functor by a natural frame homomorphism $\u A_{\crys}(R)\to\usW(R'/R)$. This extends first to the case where $R$ is quasisyntomic by $p^\infty$-root descent, and then to a functor of the groupoid cores for general $R'\to R$ by a sheafified Kan extension. Finally one verifies that the functor of groupoid cores extends to an additive functor. The case of $p$-nilpotent rings for $p\ge 3$ is similar.

n the case $R'=R$ the functor \eqref{Eq:sDR'R} will be denoted by
\beqn
\label{Eq:sDR}
\sD_R\colon \BT(R)\to\sDisp(R).
\eeqn

\subsection{A crystalline construction in characteristic $p$}

Let $R$ be a qrsp $\FF_p$-algebra. Then the ring $A_{\crys}(R)$ is torsion free, hence the natural Frobenius lift on $A_{\crys}(R)$ gives it the structure of $\delta$-ring, and we get a crystalline frame $\u A_{\crys}(R)$ in the sense of Definition \ref{De:crys-frame} with filtration ideal $I_{\crys}(R)$, the kernel of $A_{\crys}(R)\to R$. 

\bLe
\label{Le:Acrys-universal-derived}
The pd thickening $A_{\crys}(R)\to R$ is initial in the category of derived $p$-complete pd thickenings $B\to S$ with a homomorphism of $\FF_p$-algebras $R\to S$.
\eLe

\bproof
It is well known that $A_{\crys}(R) \to R$ has the above universal property
for classically $p$-complete test objects.
The derived $p$-complete derived pd thickening $\mathrm{R}
\Gamma_{\crys}^{\wedge}(R/\ZZ_p) \to R$ constructed in \cite[Definition
4.3.1 and Example 4.3.3]{Magidson:Divided} is static and $p$-torsionfree, and thus also classically $p$-complete, by the comparison with derived de Rham cohomology
\cite[Corollary 4.3.5.1]{Magidson:Divided}.
It satisfies the above universal property (even a stronger one allowing
derived test objects) and consequently agrees with $A_{\crys}(R) \to R$.
See also \cite[Theorem 8.8]{Mathew:Some-recent}.
\eproof

\bLe
\label{Le:Acris-initial-pd-frame}
For every qrsp $\FF_p$-algebra $R$ the frame $\u A_{\crys}(R)$ is initial in the category of crystalline frames $\u A$ with a homomorphism $R\to A_0/A_1$.
\eLe

\bproof
Let $\u A$ with $R\to A_0/A_1$ be given. The universal property of $A_{\crys}(R)$ of Lemma \ref{Le:Acrys-universal-derived} gives a pd homomorphism $\rho\colon A_{\crys}(R)\to A_0$. We have to show that $\rho$ is a homomorphism of $\delta$-rings, or equivalently that ($\star$) $s_A\circ\rho=W_2(\rho)\circ s_{A_{\crys}}$ where $s_A\colon A_0\to W_2(A_0)$ corresponds to the $\delta$-structure, and similarly for $A_{\crys}(R)$. The relation ($\star$) follows from the universal property of $A_{\crys}(R)$ since all four maps in ($\star$) are pd homomorphisms; see Remark \ref{Re:pd-frame}.
\eproof

\bEx
Lemma \ref{Le:Acris-initial-pd-frame} implies that for every nilpotent pd thickening of $\FF_p$-algebras $R'\to R$ where $R$ is qrsp there is a natural homomorphism of crystalline frames%
\footnote{\label{Fn:srho}This homomorphism can be constructed in a more elementary way (without reference to the notion of animated pd rings, which is implicit in \cite{Magidson:Divided}) as follows. 
There is a natural ring homomorphism $R^\flat\to R'$, hence a homomorphism of $\delta$-rings $W(R^\flat)=\sW(R^\flat)\to\sW(R')$.
Let $D$ be the non-completed pd envelope of $W(R^\flat)\to R$. Then $D$ is $p$-torsion free by \cite[Remark 8.3]{Mathew:Some-recent}, and $A_{\crys}(R)$ is the classical or derived $p$-completion of $D$. The universal property of $D$ gives a pd homomorphism $\srho'\colon D\to\sW(R')$, which is a homomorphism of $\delta$-rings by the uniqueness of the pd homomorphism $D\to W_2(\sW(R))$. Hence $\srho'$ also preserves the divided Frobenius, using the relation \eqref{Eq:pd-frame-1} on both sides. By derived $p$-completion we get a ring homomorphism $\srho\colon A_{\crys}(R)\to\sW(R')$ which preserves the divided Frobenius and the $\delta$-structure, using that $W_2(A_{\crys}(R))$ is the derived or classical $p$-completion of $W_2(D)$. Hence $\srho$ also preserves the divided powers, using  \eqref{Eq:pd-frame-1} again.}
\beqn
\label{Eq:srhoRR'}
\srho_{R'/R}\colon\u A_{\crys}(R)\to\usW(R'/R).
\eeqn
\eEx

\bDe
\label{De:sDR'/R-qrsp}
Let $R'\to R$ be a nilpotent pd thickening of $\FF_p$-algebras where $R$ is qrsp. In this case, we define the relative sheared display functor $\sD_{R'/R}$ as the composition
\beqn
\label{Eq:De:sDR'/R-qrsp}
\BT(R)\xrightarrow{\;D^{\crys}_R\;}\Win(\u A_{\crys}(R))\xrightarrow{\;\srho_{R'/R}\;}\sDisp(R'/R)
\eeqn
where $D^{\crys}_R$ is given by the covariant crystalline Dieudonn\'e functor as in \cite[Lemma A.2]{Kisin:Crystalline} (contravariant) or \cite[Proposition 3.17]{Lau:Relations}, and the second functor is induced by \eqref{Eq:srhoRR'}.
\eDe

\bLe
\label{Le:sDR'/R-qrsp-comp-crys}
The functor $\sD_{R'/R}$ of Definition \ref{De:sDR'/R-qrsp} is exact, and there is a natural isomorphism $\DD(G)_{R'/R}\cong\sD_{R'/R}(G)\otimes R'$ that preserves the Hodge filtration for $G\in\BT(R)$.
\eLe

\bproof
Both functors in \eqref{Eq:De:sDR'/R-qrsp} are exact. The comparison isomorphism follows from the construction of $D^{\crys}_R$.
\eproof

\subsection{A $p$-nilpotent version of the crystalline construction}
\label{Se:Acris-p-nilp}

Let $R$ be a $p$-nilpotent ring such that $R/p$ is a qrsp $\FF_p$-algebra. The universal property of $A_{\crys}(R/p)$ gives a homomorphism  $A_{\crys}(R/p)\to R$ of pd thickenings of $R/p$; let $I_{\crys}(R)$ be its kernel. We define a crystalline frame $\u A_{\crys}(R)\subseteq\u A_{\crys}(R/p)$ with the same underlying ring and filtration ideal $I_{\crys}(R)$.

\bLe
\label{Le:Acris-initial-pd-frame-p-nilp}
The frame $\u A_{\crys}(R)$ is initial in the category of crystalline frames $\u A$ with a ring homomorphism $R\to A_0/A_1$.
\eLe

\bproof
For given $\u A$ and $R\to A_0/A_1$ let $\tilde A_1=A_1+pA_0$. The divided powers on $A_1$ extend to $\tilde A_1$. This gives a crystalline frame $\u{\tilde A}$ with underlying ring $A_0$ and filtration ideal $\tilde A_1$. Indeed, \eqref{Eq:pd-frame-1} defines $\sigma_1\colon\tilde A_1\to A_0$, and \eqref{Eq:pd-frame-2} holds on $\tilde A_1$ because this relation defines an ideal in $\tilde A_1$ that contains $A_1$ and $p$. By the universal property of $\u A_{\crys}(R/p)$ there is a unique homomorphism of crystalline frames $\tilde\rho\colon\u A_{\crys}(R/p)\to\u{\tilde A}$. The underlying ring homomorphism preserves the augmentations to $A_0/A_1$ by the universal property of $A_{\crys}(R/p)$. Hence $\tilde\rho$ restricts to a homomorphism of crystalline frames $\rho\colon\u A_{\crys}(R)\to\u A$. Conversely, any homomorphism $\u A_{\crys}(R)\to\u A$ extends to a homomorphism $\u A_{\crys}(R/p)\to\u{\tilde A}$, hence $\rho$ is unique.
\eproof

\bEx
Lemma \ref{Le:Acris-initial-pd-frame-p-nilp} implies that for a nilpotent pd thickening $R'\to R$ of $p$-nilpotent rings where $\bar R=R/p$ is qrsp and $p\ge 3$ there is a homomorphism of crystalline frames%
\footnote{This homomorphism admits an elementary construction similar to footnote \ref{Fn:srho}. We can assume that $R=\bar R$. The Cartier homomorphism $\Delta\colon W(R^\flat)\to W(W(R^\flat))$ has image in $\sW(W(R^\flat))$ by Example \ref{Ex:Artin}. There is a natural ring homomorphism $W(R^\flat)\to R'$, hence a homomorphism of $\delta$-rings $W(R^\flat)\to\sW(W(R^\flat))\to\sW(R')$. Now proceed as before.}
\beqn
\label{Eq:srhoRR'-p=3}
\srho_{R'/R}\colon\u A_{\crys}(R)\to\usW(R'/R).
\eeqn
\eEx

\bDe
\label{De:sDR'/R-qrsp-p3}
Assume that $p\ge 3$. Let $R'\to R$ be a nilpotent pd thickening of $p$-nilpotent rings such that $\bar R=R/p$ is qrsp. In this case, we define the relative sheared display functor $\sD_{R'/R}$ by \eqref{Eq:De:sDR'/R-qrsp} with a slightly different meaning of the symbols,
where $D^{\crys}_R$ is given by the crystalline Dieudonn\'e functor again,\footnote{For the frame $\u A_{\crys}(R)$ the homomorphism $\sigma_1$ does not necessarily generate the unit ideal, which is assumed in the cited literature. One can either remark that this condition is not essential, or one can use the functor $D^{\crys}_{R/p}$ and note that windows over $\u A_{\crys}(R)$ are equivalent to windows over $\u A_{\crys}(R/p)$ together with a lift of the Hodge filtration to $R$.} 
and the second functor is induced by \eqref{Eq:srhoRR'-p=3}.
\eDe

\bRe
\label{Re:sDR'/R-qrsp-comp-crys-3}
Lemma \ref{Le:sDR'/R-qrsp-comp-crys} extends to the functors $\sD_{R'/R}$ of Definition \ref{De:sDR'/R-qrsp-p3} literally.
\eRe

\subsection{Groupoid cores of additive categories} 

For an additive category $\AAA$ let $\AAA^\simeq$ be the groupoid core (forget all non-isomorphisms) endowed with the symmetric monoidal structure $\oplus$. Any additive functor $\tilde F\colon \AAA\to\BBB$ induces a symmetric monoidal functor $F\colon \AAA^\simeq\to\BBB^\simeq$. In this situation we say that $\tilde F$ extends $F$.
We will need the following general observation.

\bLe
\label{Le:groupoid-core-extend}
Let $\AAA$, $\BBB$ be additive categories. A symmetric monoidal functor $F\colon \AAA^\simeq\to\BBB^\simeq$ extends to an additive functor $\tilde F\colon \AAA\to\BBB$ iff the following diagram commutes for all $X\in\AAA$. 
\[
\xymatrix@M+0.2em{
F(X)\oplus F(X) \ar[r]_-\sim^-{\mu_{X,X}} \ar[d]_{\left(\begin{smallmatrix}1&1\\0&1\end{smallmatrix}\right)} & F(X\oplus X) \ar[d]^{F\left(\begin{smallmatrix}1&1\\0&1\end{smallmatrix}\right)} \\
F(X)\oplus F(X) \ar[r]_-\sim^-{\mu_{X,X}} & F(X\oplus X) 
}
\]
The extension $\tilde F$ is unique if it exists. Moreover, if $\tilde F$ extends $F$ and $\tilde G$ extends $G$, each isomorphism $F\cong G$ of symmetric monoidal functors uniquely lifts to an isomorphism $\tilde F\cong\tilde G$.
\eLe

\bRe
The condition involving $\left(\begin{smallmatrix}1&1\\0&1\end{smallmatrix}\right)$ is necessary, as the following example shows. Let $\MMM$ be the category of matrices over a field $k$, with object set $\NN$ and $\Hom(n,m)$ the set of $(m\times n)$ matrices, where composition is given by matrix multiplication. The autoequivalence of the groupoid core $\MMM^\simeq$ defined by the identity on objects and by $A\mapsto (A^t)^{-1}$ on morphisms does not extend to an additive functor; it sends $\left(\begin{smallmatrix}1&1\\0&1\end{smallmatrix}\right)$ to $\left(\begin{smallmatrix}1&0\\-1&1\end{smallmatrix}\right)$.
\eRe

\bproof[Proof of Lemma \ref{Le:groupoid-core-extend}]
Assume that $F$ is a symmetric monoidal functor such that the diagram commutes for all $X$. For $X,Y\in\AAA$ we define an injective group homomorphism
\[
\psi_{X,Y}\colon \Hom(X,Y)\to \Aut(X\oplus Y),\qquad u\mapsto\left(\begin{matrix}1&u\\0&1\end{matrix}\right).
\]
If $\tilde F$ extends $F$ then there is a commutative diagram 
\[
\xymatrix@M+0.2em{
\Hom(X,Y) \ar[d]_{\psi_{X,Y}} \ar[rr]^{\tilde F} && \Hom(FX,FY) \ar[d]^{\psi_{FX,FY}} \\
\Aut(X\oplus Y) \ar[r]^-F & \Aut(F(X\oplus Y)) \ar[r]^-{\mu'} & \Aut(FX\oplus FY)
}
\]
with $\mu'(u)=\mu_{X,Y}^{-1}u\mu_{X,Y}$. Hence $\tilde F$ is unique. 
Let $\Hom(X,Y)^a$ be the set of all $h\in\Hom(X,Y)$ such that  $(\mu'\circ F\circ\psi_{X,Y})(h)$ lies in the image of $\psi_{FX,FY}$  and define 
\[
\tilde F\colon \Hom(X,Y)^a\to\Hom(FX,FY)
\] 
by $\psi_{FX,FY}\circ\tilde F=\mu'\circ F\circ\psi_{X,Y}$. Then $\Hom(X,Y)^a$ is a subgroup of $\Hom(X,Y)$, and $\tilde F$ is a  group  homomorphism. 

By functoriality in $\AAA^\simeq$, the abelian groups $\Hom(X,Y)^a$ for varying $X,Y$ are stable under multiplication by isomorphisms on both sides, and $\tilde F(uhv)=F(u)\tilde F(h)F(v)$ when $u$ and $v$ are isomorphisms. The hypothesis on $\left(\begin{smallmatrix}1&1\\0&1\end{smallmatrix}\right)$ implies that $\id_X\in\Hom(X,X)^a$ with $\tilde F(\id_X)=\id_{F(X)}$. Hence every isomorphism $u\colon X\xrightarrow\sim Y$ lies in $\Hom(X,Y)^a$,  and  $\tilde F(u)=F(u)$.

Assume that $X''=X\oplus X'$ and $Y''=Y\oplus Y'$ and let $\xi\colon \Hom(X,Y)\to\Hom(X'',Y'')$ be the homomorphism defined by $\xi(h)=\left(\begin{smallmatrix}h&0\\0&0\end{smallmatrix}\right)$. Then $\Hom(X,Y)^a$ is the inverse image under $\xi$ of $\Hom(X'',Y'')^a$, and $\xi\tilde F=\tilde F\xi$, using that $F$ preserves the natural decomposition $X''\oplus Y''=(X\oplus Y)\oplus Z$.
Now, for given $h\colon X\to Y$ we take $X'=Y$ and $Y'=X$. Then $\xi(h)=u-v$ where $u=\left(\begin{smallmatrix}h&1\\1&0\end{smallmatrix}\right)$ and $v=\left(\begin{smallmatrix}0&1\\1&0\end{smallmatrix}\right)$ are isomorphisms. If follows that $h\in\Hom(X,Y)^a$.

A similar argument show that $\tilde F$ preserves composition, using that this holds when one of the ingredients is an isomorphism.
More precisely, for $h\colon X\to Y$ and $g\colon Y\to Z$ let $\delta(g,h)=\tilde F(gh)-\tilde F(g)\tilde F(h)$. For given decompositions $X''=X\oplus X'$ and $Y''=Y\oplus Y'$ and $Z''=Z\oplus Z'$ the operator $\xi$ preserves composition and hence $\xi(\delta(g,h))=\delta(\xi (g),\xi(h))$. Now let $X'=Y$ and $Y'=X$ and $Z'=0$ so that $\xi(h)=u-v$ as before. It follows that $\delta(\xi(g),\xi(h))=\delta(u-v,\xi(h))=\delta(u,\xi(h))-\delta(v-\xi(h))=0$.

Hence $\tilde F$ is an additive functor which coincides with $F$ on the groupoid core, and $\tilde F$ induces the given symmetric monoidal structure on $F$ by the above formula $\xi\tilde F=\tilde F\xi$.
\eproof

\subsection{Reduction to quasisyntomic rings}
\label{Se:reduction-qsyn-rings}

Morphisms from a smooth algebraic stack to an etale stack can be defined over smooth rings. In this section we record a similar result for the stack $\BT^\simeq$ of $p$-divisible groups, using the smoothness properties of the stacks of truncated BT groups.
A ring homomorphism $A\to B$ will be called a countably smooth covering if $B$ is the colimit of a sequence of faithfully flat smooth ring homomorphisms $A=B_0\to B_1\to B_2\to\ldots$~. The morphisms $\BT^\simeq_{n+1}\to\BT^\simeq_n$ of algebraic stacks are smooth by a theorem of Grothendieck, \cite[Theorem 4.4]{Illusie:BT}. This implies the following.

\bLe[{\cite[Lemma 1.1]{Lau:Smoothness}}]
\label{Le:BT-presentation}
There is a countably smooth morphism $\Spec A\to\BT^\simeq$ such that $\Spec A\to\Spec\ZZ$ is a countably smooth covering as well.
\eLe

For a $p$-nilpotent ring $R$ let $\Aff_R^{\qsyn}$ be the category of quasisyntomic affine $R$-schemes and 
\[
u\colon \Aff_R^{\qsyn}\to\Aff_R
\] 
the inclusion. 

\bCo
\label{Co:BT-groupoid}
There is a groupoid $X_\bullet$ in $\Aff_R^{\qsyn}$ and a morphism $X_0\to\BT^\simeq_R$ of fibered categories over $\Aff_R$ which is ind-\'etale surjective such that $X_\bullet$ is the associated \v Cech nerve.
\eCo

\bproof
Take $X_0=\Spec(A\otimes R)$ with $A$ as in Lemma \ref{Le:BT-presentation}.
\eproof

\bLe
\label{Le:Hom-BTR-Y}
For an ind-\'etale stack $\YYY$ over $\Aff_R$ the restriction functor 
\beqn
\label{Eq:Hom-BTR-Y}
u^*\colon \Hom_{\Aff_R}(\BT^\simeq_R,\YYY)\to\Hom_{\Aff_R^{\qsyn}}(u^*\BT^\simeq_R,u^*\YYY)
\eeqn
is an equivalence. 
\eLe

\bproof
Corollary \ref{Co:BT-groupoid} provides a groupoid $X_\bullet$ in $\Aff_R^{\qsyn}$ such that $\BT^\simeq_R$ is the ind-\'etale stack associated to $X_\bullet$ in $\Aff_R$. The functor \eqref{Eq:Hom-BTR-Y} with $X_n$ in place of $\BT^\simeq_R$ is an equivalence by the Yoneda lemma. Hence \eqref{Eq:Hom-BTR-Y} is a limit of equivalences, thus an equivalence.
\eproof

For a $p$-nilpotent ring $R$ let $\nCrys_R$ be the category of nilpotent pd thickenings $Y\to X$ of affine $R$-schemes, endowed with the ind-\'etale topology as in \S\ref{Se:rel-sh-disp}, and let $\nCrys_R^{\qsyn}$ be the full subcategory where $Y$ is quasisyntomic over $R$. We have a Cartesian diagram of categories
\[
\xymatrix@M+0.2em{
\nCrys_R^{\qsyn} \ar[r]^-v \ar[d]_{\pi'} & \nCrys_R \ar[d]^\pi \\
\Aff_R^{\qsyn} \ar[r]^-u & \Aff_R
}
\]
where $\pi$ sends $Y\to X$ to $Y$.

\bRe
\label{Re:pi-cocart}
The functor $\pi$ is an opfibration, i.e.\ the functor of opposite categories is a fibration in the sense of \cite[\href{https://stacks.math.columbia.edu/tag/02XM}{Definition 02XM}]{Stacks}, because for an object $Y\to X$ of $\nCrys_R$ and a morphism $Y\to Y'$ in $\Aff_R$ we can form the object $Y'\to X\sqcup_YY'$ of $\nCrys_R$ using the fiber product of rings.
\eRe

\bLe
\label{Le:Hom-BTR-Y-nCrys}
For a ind-\'etale stack $\YYY$ over $\nCrys_R$ the restriction functor 
\beqn
\label{Eq:Hom-BTR-Y-nCrys}
v^*\colon \Hom_{\nCrys_R}(\pi^*\BT^\simeq_R,\YYY)\to\Hom_{\nCrys_R^{\qsyn}}(v^*\pi^*\BT^\simeq_R,v^*\YYY)
\eeqn
is an equivalence. 
\eLe

\bproof
Let the groupoid $X_\bullet$ in $\Aff_R^{\qsyn}$ and the morphism $X_0\to\BT^\simeq_R$ be given by Corollary \ref{Co:BT-groupoid}. Since the functor $\pi$ induces equivalences of the ind-\'etale coverings of a given object and its image, $\pi^*X_0\to\pi^*\BT^\simeq_R$ is ind-\'etale surjective with \v Cech nerve $\pi^*X_\bullet$. In particular $\pi^*\BT^\simeq$ is the ind-\'etale stack associated to the groupoid $\pi^*X_\bullet$, and it suffices to show that 
\beqn
\label{Eq:Hom-uXn-Y-nCrys}
v^*\colon \Hom_{\nCrys_R}(\pi^*X_n,\YYY)\to\Hom_{\nCrys_R^{\qsyn}}(v^*\pi^*X_n,v^*\YYY)
\eeqn
is an equivalence. The right hand side can be identified with $\Hom_{\nCrys_R^{\qsyn}}(\pi'^*X_n,v^*\YYY)$. Let us construct an inverse of \eqref{Eq:Hom-uXn-Y-nCrys}, using that $\pi$ is an opfibration by Remark \ref{Re:pi-cocart}. For a given functor $\alpha\colon \pi'^*X_n\to v^*\YYY$ over $\nCrys_R^{\qsyn}$ we have to construct an extension of $\alpha$ to a functor $\tilde\alpha\colon \pi^*X_n\to\YYY$ over $\nCrys_R$. Let an object $Z=(Y\to X)$ of $\nCrys_R$ and an element $g\in\pi^*X_n(Z)$ be given, thus $g\colon Y\to X_n$. Then $Z'=(X_n\to X\sqcup_YX_n)$ lies in $\nCrys_R^{\qsyn}$, and we define $\tilde\alpha(g)$ to be the image of $\id_{X_n}$ under 
\[
\pi'^*X_n(Z')\xrightarrow\alpha v^*\YYY(Z')=\YYY(Z')\to\YYY(Z).
\qedhere
\]
\eproof

\subsection{The absolute sheared display functor in characteristic $p$}

We consider the category $\Aff_{\FF_p}$, its subcategories $\Aff_{\FF_p}^{\semiperf}$ (semiperfect) and $\Aff_{\FF_p}^{\qsyn}$ (quasisyntomic), and 
their intersection $\Aff_{\FF_p}^{\qrsp}$ (quasiregular semiperfect). 
The fibered categories of $p$-divisible groups and sheared displays will be denoted by $\BT\to\Aff_{\FF_p}$ and $\sDisp\to\Aff_{\FF_p}$ by a temporary change of notation. Here $\BT$ is an fpqc stack, and $\sDisp$ is (at least) an ind-\'etale and $p^\infty$-root stack. For $c\in\{\qsyn,\qrsp\}$ let $\BT^c\to\Aff_{\FF_p}^c$ and $\sDisp^c\to\Aff_{\FF_p}^c$ be the restrictions.

\bCon 
\label{Con:sh-disp-fun-p}

Definition \ref{De:sDR'/R-qrsp} for $R'=R$ gives a fiberwise additive functor 
\[
\sD^{\qrsp}\colon \BT^{\qrsp}\to\sDisp^{\qrsp}
\] 
over $\Aff_{\FF_p}^{\qrsp}$, 
which extends uniquely to a fiberwise additive functor 
\[
\sD^{\qsyn}\colon \BT^{\qsyn}\to\sDisp^{\qsyn}
\] 
over $\Aff_{\FF_p}^{\qsyn}$ by $p^\infty$-root descent. The restriction of $\sD^{\qsyn}$ to the groupoid cores extends uniquely to a morphism
\[
\sD^\simeq\colon \BT^\simeq\to \sDisp^\simeq
\]
of stacks over $\Aff_{\FF_p}$
by Lemma \ref{Le:Hom-BTR-Y}.
The functor $\sD^\simeq$ has a natural extension to a fiberwise symmetric monoidal functor with respect to the direct sum. Indeed, we have two functors $\BT^\simeq\times\BT^\simeq\to\sDisp^\simeq$ over $\Aff_{\FF_p}$ given by $\sD^\simeq(G\oplus H)$ and $\sD^\simeq(G)\oplus\sD^\simeq(H)$, they are isomorphic over $\Aff_{\FF_p}^{\qsyn}$ by the additive structure, hence they are isomorphic 
by $p^\infty$-root descent and Lemma \ref{Le:Hom-BTR-Y} again. 
The remaining axioms of a symmetric monoidal structure follow similarly. 
The symmetric monoidal functor $\sD^\simeq$ satisfies the $\left(\begin{smallmatrix}1&1\\0&1\end{smallmatrix}\right)$-condition of Lemma \ref{Le:groupoid-core-extend} because this holds over $\Aff_{\FF_p}^{\qsyn}$, and hence $\sD^\simeq$ extends to a fiberwise additive functor
\[
\sD\colon \BT\to \sDisp
\]
over $\Aff_{\FF_p}$.
\eCon

\subsection{The relative sheared display functor in characteristic $p$}
\label{Se:rel-sheared-disp-fun-p}

We consider the category $\nCrys_{\FF_p}$ of nilpotent thickenings of affine $\FF_p$-schemes and the functor $\pi\colon \nCrys_{\FF_p}\to\Aff_{\FF_p}$ as in \S\ref{Se:reduction-qsyn-rings}.
We have fibered categories $\sDisp_{\rel}\to\nCrys_{\FF_p}$ and $\pi^*\BT\to\nCrys_{\FF_p}$. The latter is an fpqc stack, and $\sDisp_{\rel}$ is (at least) an ind-\'etale and $p^\infty$-root stack.
For $c\in\{\qsyn,\qrsp\}$ let $\nCrys_{\FF_p}^c\subseteq\nCrys_{\FF_p}$ be the preimage of $\Aff_{\FF_p}^c\subseteq\Aff_{\FF_p}$ under $\pi$,
and let $\pi^*\BT^c$ and $\sDisp_{\rel}^c$ be the restrictions of $\pi^*\BT$ and $\sDisp_{\rel}$ to $\nCrys_{\FF_p}^c$.
Now Construction \ref{Con:sh-disp-fun-p} carries over to the relative case as follows.

\bCon 
\label{Con:sh-disp-fun-rel}

Definition \ref{De:sDR'/R-qrsp} gives a fiberwise additive functor 
\[
\sD_{\rel}^{\qrsp}\colon \pi^*\BT^{\qrsp}\to\sDisp_{\rel}^{\qrsp}
\]
over $\nCrys_{\FF_p}^{\qrsp}$,
which extends uniquely to a fiberwise additive functor 
\[
\sD_{\rel}^{\qsyn}\colon \pi^*\BT^{\qsyn}\to\sDisp_{\rel}^{\qsyn}
\]
over $\nCrys_{\FF_p}^{\qsyn}$ by $p^\infty$-root descent. The restriction of $\sD_{\rel}^{\qsyn}$ to the groupoid cores extends uniquely to a morphism
\[
\sD^\simeq_{\rel}\colon \pi^*\BT^\simeq\to \sDisp^\simeq_{\rel}
\]
of stacks over $\nCrys_{\FF_p}$ by Lemma \ref{Le:Hom-BTR-Y-nCrys}.
As in the absolute case (Construction \ref{Con:sh-disp-fun-p}), the functor $\sD^\simeq_{\rel}$ has a natural extension to a fiberwise symmetric monoidal functor with respect to the direct sum that satisfies the $\left(\begin{smallmatrix}1&1\\0&1\end{smallmatrix}\right)$-condition of Lemma \ref{Le:groupoid-core-extend}, hence it extends to a fiberwise additive functor
\[
\sD_{\rel}\colon \pi^*\BT\to\sDisp_{\rel}
\]
over $\nCrys_{\FF_p}$. 
\eCon

\bCon
\label{Con:sh-disp-fun-rel-3}
 For  $p\ge 3$, Construction \ref{Con:sh-disp-fun-rel} can be extended to give an additive functor $\sD_{\rel}\colon \pi^*\BT\to\sDisp_{\rel}$ over $\nCrys_{\ZZ_p}$ as follows. In the first paragraph of \S\ref{Se:rel-sheared-disp-fun-p} we replace $\FF_p$ by $R=\ZZ/p^r$ and consider the categories $\Aff_{R}^{\qsyn}$ of quasisyntomic affine $R$-schemes, $\Aff_{R}^{\semiperf}$ of affine $R$-schemes which are semiperfect modulo $p$, and their intersection $\Aff_{R}^{\qrsp}$. Now the construction works literally, starting with Definition \ref{De:sDR'/R-qrsp-p3} instead of Definition \ref{De:sDR'/R-qrsp}.
\eCon

\subsection{Conclusions}
\label{Se:Sheared-display-functor:conclusions}

\bproof[Proof of Theorem \ref{Th:rel-sh-disp-functor}]
The functor $\sD_{\rel}$ of Constructions \ref{Con:sh-disp-fun-rel} and \ref{Con:sh-disp-fun-rel-3} gives the functors $\sD_{R'/R}$ of \eqref{Eq:sDR'R}. The functor $\sD_{R'/R}$ is exact because this can be verified over perfect fields, and it holds when $R$ is qrsp by Lemma \ref{Le:sDR'/R-qrsp-comp-crys}. This lemma and Remark \ref{Re:sDR'/R-qrsp-comp-crys-3} also give a canonical isomorphism \eqref{Eq:DDG-sDGrel} when $R$ is qrsp, which extends to the case where $R$ is quasisyntomic by $p$-root descent, and then to the general case by Kan extension.
\eproof

\bPr
\label{Pr:sD-QpZp}
We have $\sD_R(\QQ_p/\ZZ_p)=\usW(R)(1)$;
see Example \ref{Ex:win-uA-uA(1)}.
\ePr

\bproof
We can assume that $R=\ZZ/p^r$ with $r=1$ when $p=2$. Then the assertion is classical; note that $\sD_R$ is given by Definitions \ref{De:sDR'/R-qrsp} and \ref{De:sDR'/R-qrsp-p3}, and $A_{\crys}(R/p)=\sW(\FF_p)=W(\FF_p)$.
\eproof

\bPr
The functor $\sD_{R'/R}$ preserves duality, i.e.\ for every nilpotent pd thickening $R'\to R$ as in Theorem \ref{Th:rel-sh-disp-functor} and every $G\in\BT(R)$ there is a natural isomorphism
\beqn
\label{Eq:sD-duality}
\sD_{R'/R}(G^\vee)\cong\sD_{R'/R}(G)^t.
\eeqn
\ePr

\bproof 
The crystalline duality theorem \cite[\S5.3]{BBM} implies that the functor $D_R^{\crys}$ in \eqref{Eq:De:sDR'/R-qrsp} preserves duality,
including the case where $R$ is $p$-nilpotent used in Definition \ref{De:sDR'/R-qrsp-p3}. This gives a duality isomorphism \eqref{Eq:sD-duality} over $\nCrys_{\FF_p}^{\qrsp}$, or over $\nCrys_{\ZZ/p^r}^{\qrsp}$ when $p\ge 3$, which carries over to all of $\nCrys_{\FF_p}$ or $\nCrys_{\ZZ/p^r}$ by Constructions \ref{Con:sh-disp-fun-rel} and \ref{Con:sh-disp-fun-rel-3}.
\eproof

\subsection{Relation with the display functor}

By \cite{Lau:Smoothness}, for every $p$-nilpotent ring $R$ there is an additive functor $D_R\colon\BT(R)\to\Disp(R)$ which is functorial in $R$, such that the underlying filtered $F$-module is given by the  Dieudonn\'e crystal, and this determines $D_R$ uniquely. Similarly, for every pd thickening of $p$-nilpotent rings $R'\to R$ there is an additive functor $D_{R'/R}\colon\BT(R)\to\Disp(R'/R)$ with the same property and uniqueness.

\bPr
\label{Pr:compare-sD-D-rel}
For every nilpotent pd thickening $R'\to R$ as in Theorem \ref{Th:rel-sh-disp-functor}, the composition
\[
\BT(R)\xrightarrow{\sD_{R'/R}}\sDisp(R'/R)\to\Disp(R'/R)
\]
is isomorphic to $D_{R'/R}$, and the isomorphism is functorial in $(R'\to R)$.
\ePr

\bproof
Let $E_{R'/R}\colon\BT(R)\to\Disp(R'/R)$ be the composition and let 
\[
\bar D_{R'/R},\,\bar E_{R'/R}\colon\BT(R)\to(\text{filtered $F$-modules over }\u W(R'/R))
\]
be the functors associated to $D_{R'/R}$ and $E_{R'/R}$ by forgetting the divided Frobenius. It suffices to construct an isomorphism $\bar E_{R'/R}\cong\bar D_{R'/R}$ when $R$ is quasisyntomic over $\ZZ/p^r$ and $R/p$ is semiperfect (with $r=1$ when $p=2$). In this case we have homomorphisms of crystalline frames $\u A_{\crys}(R)\to\usW(R'/R)\to\u W(R'/R)$, and by the construction of $\sD_{R'/R}$ in Definitions \ref{De:sDR'/R-qrsp} and \ref{De:sDR'/R-qrsp-p3} it follows that $\bar E_{R'/R}(G)$ is given by the Dieudonn\'e crystal as required.
\eproof

\cclearpage
\section{Dieudonn\'e equivalence and consequences}

By descent, the functor $\sD_R$ is an equivalence for all $\FF_p$-algebras if this holds for semiperfect $\FF_p$-algebras. This case can be deduced from the perfect case using deformations over Frobenius by a method borrowed from \cite[\S5.2]{Lau:Semiperf}. As a preparation we need that the categories of $p$-divisible groups and sheared displays preserve certain limits and colimits. 

\subsection{Limits}
\label{Se:Dieud-equiv:limits}

Let $R$ be a semiperfect $\FF_p$-algebra and $R^\flat=R\limperf$ its limit perfection.

\bLe
\label{Le:BTRflat-limBTRF}
The natural functor $\BT(R^\flat)\to\lim(\BT(R),F)$ is an equivalence.
\eLe

\bproof
We have $R^\flat=\lim(R,F)$ where the transition maps are surjective with locally nilpotent kernel. The exact categories of locally free modules and of finite locally free group schemes preserve such limits, and hence the same holds for $p$-divisible groups.
\eproof

\bLe
\label{Le:WRflat-limsWRF}
The natural homomorphism $\pi\colon W(R^\flat)\to\lim(\sW(R),F)$ is an isomorphism.
\eLe

\bproof
The homomorphism $\pi$ exists since $W(R^\flat)=\sW(R^\flat)$.
We have an exact sequence 
\beqn
\label{Eq:KR-sWR-WR}
0\to K(R)\to\sW(R)\to W(R)\to 0
\eeqn
where $\sW(R)=W(R^\flat)/\hat W(J)$ as in Proposition \ref{Pr:sW-semiperf} and $K(R)=W(J)/\hat W(J)$. The sequence induces a homomorphism
$u\colon \lim(\sW(R),F)\to\lim(W(R),F)=W(R^\flat)$ with $u\pi=\id$.
Now $\lim(K(R),F)=0$ since $K(R)=\lim_nK_n(R)$ with $K_n(R)=W(R[F^n])/\hat W(R[F^n])$ where $F^n$ annihilates $K_n(R)$. Hence $u$ is injective, thus $u$ and $\pi$ are bijective.
\eproof

For a ring $A$ let $\LF(A)$ be the category of finite projective $A$-modules.

\bLe
The natural functor $\LF(W(R^\flat))\to\lim(\LF(\sW(R)),F)$ is an equivalence.
\eLe

\bproof
The homomorphism $F\colon \sW(R)\to\sW(R)$ is surjective with kernel $\hat W(R[F^n])$, which is a locally nilpotent ideal, and the limit is $W(R^\flat)$ by Lemma \ref{Le:WRflat-limsWRF}. Hence \cite[\href{https://stacks.math.columbia.edu/tag/0D4B}{Lemma 0D4B}]{Stacks} applies.
\eproof

Recall that $\sDisp(R)=\sDisp^{\plain}(R)=\Win(\usW(R))$ by Corollary \ref{Co:essemiperf-plain}.

\bLe
\label{Le:DispRflat-limsDispRF}
The natural functor $\Disp(R^\flat)\to\lim(\sDisp(R),F)$ is an equivalence.
\eLe

\bproof
Use Lemma \ref{Le:WRflat-limsWRF} and the analogous equivalence $\LF(R^\flat)\xrightarrow\sim\lim(\LF(R),F)$.
\eproof

\subsection{Partial finite presentation} 

Let $R$ be a semiperfect $\FF_p$-algebra.

\bLe
\label{Le:finite-pres-BT}
Let $\Fa\subseteq R$ be an ideal contained in $R[F^n]$ for some $n$ and let $E$ be the set of all finitely generated ideals $\Fa'$ contained in $\Fa$. Then
$\BT(R/\Fa)=\colim_{\Fa'\in E}\BT(R/\Fa')$.
\eLe

\bproof
Assume that $\pi\colon S\to R$ is surjective with kernel contained in $S[F^n]$ and $S$ is semiperfect. There is a homomorphism $\psi\colon R\to S$ with $\psi\circ \pi=F^n$. For a $p$-divisible group $G$ over $S$ we have an exact sequence $0\to K\to G^{(p^n)}\to G\to 0$ where $K=(G^\vee[F^n])^\vee$. Here $K$ is a finite locally free group scheme contained in $G^{(^pn)}[p^n]=\psi^*(G_R[p^n])$.
This construction gives an equivalence between $\BT(S)$ and the category of pairs $(H,K)$ where $H\in\BT(R)$ and $K\subseteq\psi^*(H[p^n])$ is $n$-cosmooth with $n$-smooth quotient.
For $\Fa'\in E$, the homomorphism $F^n\colon R\to R$ factors as $R\to R/\Fa'\to R/\Fa\to R$. We apply the preceding equivalence to $S\to R$ where $S=R/\Fa'$ or $S=R/\Fa$. It follows that $\BT(R/\Fa)$ is equivalent to the category of pairs $(H,K)$ where $H\in\BT(R)$ and where $K\subseteq H_{R/\Fa}[p^n]$ is $n$-cosmooth with $n$-smooth quotient. For given $H$, such $K$ are parametrised by a projective scheme over $R$ of finite presentation,  and hence the value over $R/\Fa$ is the colimit over $E$ of the values over $R/\Fa'$. 
\eproof

\bLe
\label{Le:finite-pres-sDisp}
In the situation of Lemma \ref{Le:finite-pres-BT} we have $\sDisp(R/\Fa)=\colim_{\Fa'\in E}\sDisp(R/\Fa')$.
\eLe

\bproof
For $\Fb\subseteq R[F^n]$ we have $0\to\hat W(\Fb)\to\sW(R)\to\sW(R/\Fb)\to 0$ by Proposition \ref{Pr:WN-sWR-sWRN}.
Since $\hat W(\Fa)$ is the colimit of $\hat W(\Fa')$, it follows that $\usW(R/\Fa)$ is the colimit of $\usW(R/\Fa')$, and hence $\sDisp(R/\Fa)=\Win(\usW(R/\Fa))$ is the colimit of $\sDisp(R/\Fa')$ by Remark \ref{Re:normal-repres}.
\eproof

\subsection{Deformations over semiperfect $\FF_p$-algebras}
\label{Se:Dieud.equiv:deform}

Let $R$ be a semiperfect $\FF_p$-algebra.

\bPr
\label{Pr:Deform-semiperf}
For an ideal $\Fa$ of $R$ with $\Fa\subseteq R[F^n]$ for some $n$, the following square of categories is $2$-Cartesian, and the horizontal arrows are essentially surjective.
\[
\xymatrix@M+0.2em{
\BT(R) \ar[r] \ar[d]_{\sD_R} & \BT(R/\Fa) \ar[d]^{\sD_{R/\Fa}} \\
\sDisp(R) \ar[r] & \sDisp(R/\Fa)
}
\]
\ePr

\bproof
By Lemmas \ref{Le:finite-pres-BT} and \ref{Le:finite-pres-sDisp} we can assume that $\Fa$ is finitely generated. Then we can assume that $\Fa$ has square zero, in particular $\Fa$ carries nilpotent divided powers. In that case the classical deformation argument applies: Writing $\bar R=R/\Fa$, the diagram extends to
\[
\xymatrix@M+0.2em{
\BT(R) \ar[r] \ar[d]_{\sD_R} \ar@{}[dr]|\square & \BT(\bar R) \ar[r]^-= \ar[d]^{\sD_{R/\bar R}} & \BT(\bar R) \ar[d]^{\sD_{\bar R}} \\
\sDisp(R) \ar[r] & \sDisp(R/\bar R) \ar[r]^-\sim & \sDisp(\bar R).
}
\]
where $\sDisp(-)$  coincides with $\sDisp^{\plain}(-)$ in all three cases  by Corollaries \ref{Co:essemiperf-plain} and \ref{Co:essemiperf-rel-plain}.
The lower right functor is an equivalence by Proposition \ref{Pr:sDisp(R'/R)-sDisp(R)}. The left square is $2$-Cartesian because the horizontal fibers correspond to lifts of the Hodge filtration in the same way by the Grothendieck--Messing theorem for $p$-divisible groups, its trivial window analogue as in \S\ref{Se:Hodge-fil}, and the comparison isomorphism \eqref{Eq:DDG-sDGrel}.
The functor  $\sDisp(R)\to\sDisp(R/\bar R)$  is essentially surjective since lifts of the Hodge filtration exist. 
\eproof

\bCo
\label{Co:Deform-Rflat-R}
The following   square of categories  is $2$-Cartesian, and the horizontal arrows are essentially surjective.
\[
\xymatrix@M+0.2em{
\BT(R^\flat) \ar[r] \ar[d]^{D_{R^\flat}} & \BT(R) \ar[d]^{\sD_R} \\
\Disp(R^\flat) \ar[r] & \sDisp(R)
}
\]
\eCo

\bproof
By Lemmas \ref{Le:BTRflat-limBTRF} and \ref{Le:DispRflat-limsDispRF}, the square is the limit of similar squares with $R$ in place of $R^\flat$ and horizontal arrows $F^n$. These are $2$-Cartesian with essentially surjective horizontal arrows by Proposition \ref{Pr:Deform-semiperf}.
\eproof

\bTh
\label{Th:sD-equiv-Fp}
For every $\FF_p$-algebra $R$, the functor $\sD_R$ is an equivalence.
\eTh

\bproof
The functors $\sD_R$ form a morphism of fibered categories over $\Aff_{\FF_p}$, which are $p^\infty$-root stacks, and semiperfect $\FF_p$-algebras form a basis of this topology. Hence we may assume that $R$ is semiperfect.
In the diagram of Corollary \ref{Co:Deform-Rflat-R}, the functor $D_{R^\flat}$ is an equivalence by a theorem of Gabber; see \cite[Theorem D]{Lau:Smoothness}. The result follows; see \cite[Lemma 5.9]{Lau:Semiperf}.
\eproof

\bRe
In the arguments leading to Theorem \ref{Th:sD-equiv-Fp} we use that over semiperfect $\FF_p$-algebras, all sheared displays are plain. Using the present descent results for finite projective $\sW(-)$-modules, which are limited to ind-\'etale and $p^{\infty}$-root coverings, this is ensured by the definition of sheared displays over general $p$-nilpotent rings by $p^{\infty}$-root descent. Even when restricting to semiperfect $\FF_p$-algebras, such a definition is required for the construction of the functor $\sD_R$, which proceeds as $\qrsp\to\qsyn\to$ general $\to$ semiperfect.
\eRe

\subsection{Conclusions}
\label{Se:Dieud-equiv:conclusions}

\bCo
\label{Co:sD-equiv-3}
If $p\ge 3$, the functor $\sD_R$ is an equivalence for every $p$-nilpotent ring.
\eCo

\bproof
Similar to the proof of Lemma \ref{Le:sFZ-eq-Fp-gen}, using that the stack $\BT$ is formally smooth.
\eproof

\bCo
\label{Co:sD-sFZ}
The composition of functors $\BT\xrightarrow{\sD}\sDisp\xrightarrow{\sFZ}\BT$ between categories fibered over $\Aff_{\FF_p}$, or over $\Aff_{\Spf\ZZ_p}$ when $p\ge 3$, is isomorphic to the identity.
\eCo

\bproof
Let $G\in\BT(R)$ and $\uM=\sD_R(G)$. 
For an $R$-algebra $S$ we have
\[
T_pG(S)=\Hom_S(\QQ_p/\ZZ_p,G_S)=\Hom_{\sDisp(S)}(\usW(S)(1),\uM_S)=T(\u M)(S)
\] 
using that $\sD_R$ is fully faithful, Proposition \ref{Pr:sD-QpZp}, and Lemma \ref{Le:Hi-Gamma-uM}.
\eproof

\bCo
\label{Co:sFZ-equiv}
The functor $\sFZ_R\colon\sDisp(R)\to\BT(R)$ provided by Therorem \ref{Th:sC-gives-BT} is an equivalence for every $p$-nilpotent ring $R$.
\eCo

\bproof
For $p\ge 3$ this follows from Theorem \ref{Th:sD-equiv-Fp} and Corollary \ref{Co:sD-sFZ}. For $p=2$ this argument gives that $\sFZ_R$ is an equivalence when $R$ is an $\FF_p$-algebra. Then use Lemma \ref{Le:sFZ-eq-Fp-gen}.
\eproof

\bDe
\label{De:sD-general-p=2}
For $p=2$ we define $\sD\colon\BT\to\sDisp$ as the inverse of the functor $\sFZ$.
\eDe

\bRe
\label{Re:sD-exact-duality}
By Corollary \ref{Co:sD-sFZ}, the restriction of $\sD$ to $\Aff_{\FF_p}$ coincides with the functor provided by Theorem \ref{Th:rel-sh-disp-functor}.
The functor $\sD$ is exact since this can be verified over $\FF_p$-algebras, and compatible with duality by Proposition \ref{Pr:sFZ-duality}. 
\eRe

\bproof[Proof of Theorem \ref{Th:Intro-sFZR}]
The functor $\sFZ_R$ is well-defined and exact by Theorem \ref{Th:sC-gives-BT}, it preserves duality by Proposition \ref{Pr:sFZ-duality}, the relation with the Zink functor $\FZ_R$ is given by Proposition \ref{Pr:FZ-sFZ}, it is an equivalence by Corollary \ref{Co:sFZ-equiv}, and the inverse is exact by Remark \ref{Re:sD-exact-duality}.
\eproof

\subsection{Rigidity and uniqueness}
\label{Se:Dieud-equiv:rigidity}

Using rigidity properties of the fibered category of $p$-divisible groups, the functor $\sD$ can be characterised without reference to the functor $\sFZ$ as follows.
Let $\SSS=\Aff_{\FF_p}\!$ or $\SSS=\Aff_{\Spf\ZZ_p}\!$ and let $\BT_{\SSS}\to\SSS$ denote the fibered category of $p$-divisible groups.

\bLe
\label{Le:End-BT-Fp-Zp}
Let $\Psi$ be a fiberwise additive endomorphism of\/ $\BT_{\SSS}$. Every homomorphism $h_0\colon\QQ_p/\ZZ_p\to\Psi(\QQ_p/\ZZ_p)$ of $p$-divisible groups over $\FF_p$ uniquely extends to a natural homomorphism $h_G\colon G\to\Psi(G)$ for all $G\in\BT_{\SSS}$. If $h_0$ is an isomorphism and $\Psi$ is fully faithful for unipotent $p$-divisible groups over $\FF_p$-algebras, then $h_G$ is an isomorphism for all $G$.
\eLe

\bproof
Assume first that $\SSS=\Aff_{\FF_p}$.
The homomorphisms $h_0$ and $h_G$ correspond to homomorphisms of the respective Tate module schemes, and $h_G$ exists uniquely for all $G$ since $T_pG=\u\Hom(\QQ_p/\ZZ_p,G)$. Under the additional hypotheses, $h_G$ is an isomorphism when $G$ is unipotent, and by functoriality it remains to show that $h_D$ is an isomorphism for $D=\mu_{p^\infty}$ over some field.
Let $G=E[p^{\infty}]$ where $E$ is an elliptic curve over $R=k[[t]]$ with supersingular special fiber and ordinary generic fiber, for $k=\bar\FF_p$. Then $h_G$ is an isomorphism in the special fiber, hence an isomorphism over $R$ by Lemma \ref{Le:Hom-BT-isom-clopen}. It follows that $h_D$ is an isomorphism over the algebraic closure of $k((t))$.
Assume now that $\SSS=\Aff_{\Spf\ZZ_p}$. By rigidity of \'etale $p$-divisible groups, $h_0$ extends to an isomorphism $\QQ_p/\ZZ_p\cong\Psi(\QQ_p/\ZZ_p)$ over $\ZZ_p$, and $h_G$ exists uniquely as before. Under the additional hypotheses, $h_G$ is an isomorphism because this can be verified over the reduction modulo $p$.
\eproof

\bCo
\label{Co:End-BT-Fp-Zp-dual}
Let $\Psi$ be a fiberwise additive endomorphism of\/ $\BT_{\SSS}$. Every homomorphism $h_0\colon\Psi(\mu_{p^{\infty}})\to\mu_{p^{\infty}}$ of $p$-divisible groups over $\FF_p$ uniquely extends to a natural homomorphism $h_G\colon\Psi(G)\to G$ for all $G\in\BT_{\SSS}$. If $h_0$ is an isomorphism and $\Psi$ is fully faithful for infinitesimal $p$-divisible groups over $\FF_p$-algebras, then $h_G$ is an isomorphism for all $G$.
\eCo

\bproof
Apply Lemma \ref{Le:End-BT-Fp-Zp} to $\Psi'=\vee\circ\Psi\circ\vee$.
\eproof

\bRe
\label{Re:End-BT-Fp-Zp}
The proof of Lemma \ref{Le:End-BT-Fp-Zp} shows that the lemma holds for the fibered category $\BT_{\SSS,\uni}$ of unipotent $p$-divisible groups in place of $\BT_\SSS$, and Corollary \ref{Co:End-BT-Fp-Zp-dual} holds for the fibered category $\BT_{\SSS,\inf}$ of infinitesimal $p$-divisible groups in place of $\BT_\SSS$.
\eRe

\bRe
\label{Re:FZ-D-equiv}
Recall the functors $\FZ\colon\Disp\to\FG$ defined in \cite{Zink:The-display} and $D\colon\BT\to\Disp$ defined in \cite{Lau:Smoothness} between fibered categories over $\Aff_{\Spf\ZZ_p}$. The composition $\FZ\circ D$ is isomorphic to the formal completion functor, and the functors restrict to mutually inverse equivalences between infinitesimal $p$-divisible groups and $V$-nilpotent displays,
$D_{\nil}\colon \BT_{\inf}\xrightarrow\sim\Disp_{\Vnil}$ and $\FZ_{\inf}\colon\Disp_{\Vnil}\xrightarrow\sim\BT_{\inf}$, by loc.\,cit.
By Remark \ref{Re:End-BT-Fp-Zp}, the equivalence $D_{\inf}$ is uniquely determined by the property of being fully faithful and the value $D_{\inf}(\mu_{p^\infty})=\u W(\FF_p)$.
\eRe

\bPr
\label{Pr:sD-unique}
The functor $\sD\colon \BT\to\sDisp$ of fibered categories over $\Aff_{\Spf\ZZ_p}$ given by Theorem \ref{Th:rel-sh-disp-functor} for $p\ge 3$ and by Definition \ref{De:sD-general-p=2} for $p=2$ is uniquely determined by either of the following properties. Let $E_{\inf}$ be the composition $\BT_{\inf}\to\BT\xrightarrow{\sD}\sDisp\to\Disp$.
\begin{enumerate}
\item
\label{It:sD-unique-1}
The functor $E_{\inf}$ is isomorphic to the functor $D_{\inf}$ of Remark \ref{Re:FZ-D-equiv}.
\item
\label{It:sD-unique-2}
The restrictions of $E_{\inf}$ and $D_{\inf}$ over $\Aff_{\FF_p}$ are isomorphic.
\item
\label{It:sD-unique-3}
The functor $E_{\inf}$ restricted to $\Aff_{\FF_p}$ is fully faithful, it has image in $\Disp_{\Vnil}$, and takes the value $E_{\inf}(\mu_{p^\infty})=\u W(\FF_p)$.
\end{enumerate}
\ePr

\bproof
For any $p$, the restriction of $\sD$ to $\Aff_{\FF_p}$ is given by Theorem \ref{Th:rel-sh-disp-functor}, and Proposition \ref{Pr:compare-sD-D-rel} implies that $\sD$ has property \eqref{It:sD-unique-2}, thus also \eqref{It:sD-unique-3}. The Rim--Schlessinger condition implies that $E_{\inf}$ is fully faithful since this holds in characteristic $p$, hence \eqref{It:sD-unique-1} holds for $\sD$. Assume that $\sD'\colon \BT\to\sDisp$ is another functor that satisfies \eqref{It:sD-unique-3} and let $\Psi=\sFZ\circ\sD'$. Then $\Psi$ preserves infinitesimal $p$-divisible groups, $\Psi(\mu_{p^\infty})=\mu_{p^\infty}$, and the composition
\[
\BT_{\inf}\longrightarrow\BT\xrightarrow{\;\sD'\;}\sDisp\xrightarrow{\;\sFZ\;}\BT\xrightarrow {\text{compl.}}\FG
\]
is isomorphic to $\BT_{\inf}\to\BT\xrightarrow{\sD'}\sDisp\to\Disp\xrightarrow{\FZ}\FG$ by Proposition \ref{Pr:FZ-sFZ}. Hence $\Psi$ induces a fully faithful endomorphism of $\BT_{\inf,\FF_p}$, and Corollary \ref{Co:End-BT-Fp-Zp-dual} gives the result.
\eproof

\bproof[Proof of Theorem \ref{Th:Intro-sDR}]
For $pR=0$ or $p\ge 3$, the exact functor $\sD_R$ is given by Theorem \ref{Th:sD-equiv-Fp} and is an equivalence by Theorem \ref{Th:sD-equiv-Fp}. This construction extends to an equivalence $\sD_R$ in general by Corollaries \ref{Co:sD-sFZ} and \ref{Co:sFZ-equiv} (see Definition \ref{De:sD-general-p=2}), and $\sD_R$ preserves the exact structure and duality by Remark \ref{Re:sD-exact-duality}. The functor $\sD\colon\BT\to\sDisp$ of fibered categories is compatible with $D_{\nil}\colon\BT_{\inf}\to\Disp_{\Vnil}$, and this property determines $\sD$, by Proposition  \ref{Pr:sD-unique}.
\eproof

\cclearpage
\section{Descent of finite projective $\sW$-modules, II}
\label{Se:Descent-II}

A $p$-nilpotent ring $R$ will be called good if $R_{\red}$ is perfect.
According to \cite{BMVZ}, the higher flat cohomology of $\hat W$ and $\sW$ vanishes over good rings;
this goes beyond the elementary theory of $\sW$ as presented in \cite{Drinfeld:Ring-stacks-conjecturally} or \S\ref{Se:Sheard-Witt-vect}.
It follows that finite projective modules over $\sW(R)$ satisfy fpqc descent with respect to $R$ when $R$ is good, by an inductive argument (due to Drinfeld) for descent of finite projective modules over the animated rings $\sW(R)/^Lp^n$.
In this section, we present a classical argument  based on the exact sequences $0\to T_FQ\to\sW\to W\to 0$.

\subsection{Vanishing of more cohomology}

The following vanishing result is proved in \cite{BMVZ} as a consequence of \cite[Proposition 7.4.8]{Bhatt-Lurie:Absolute} and \cite[Theorem 1.18]{Bragg-Olsson:Representability}.

\bPr
\label{Pr:Hi-fpqc-hatW-sW}
If $R$ is good, then $H^i_{\fpqc}(R,\hat W)=0$ and $H^i_{\fpqc}(R,\sW)=0$ for $i>0$.
\ePr

\bRe
A priori, the definition of fpqc cohomology requires the choice of a suitable cardinality bound as for example in \cite[\S1.4]{Cesnavicius-Scholze}. This is not relevant here because for a given fpqc sheaf $\FFF$, we have $H^i_{\fpqc}(R,\FFF)=0$ for every good ring $R$ and $i>0$ iff for every faithfully flat homomorphism $R\to R'$ between good rings, we have $\check H^i(R'/R,\FFF)=0$ for $i>0$.
\eRe

\bCo
\label{Co:Hi-fpqc-TFQ}
If $R$ is a semiperfect $\FF_p$-algebra, then $H^i_{\fpqc}(R,T_FQ)=0$ for $i>0$.
\eCo

\bproof
Since $R$ is semiperfect, $F$ is surjective on $Q(R)$, hence $0\to T_FQ(R)\to\sW(R)\to W(R)\to 0$ is exact. For a faithfully flat homomorphism of semiperfect $\FF_p$-algebras $R\to R'$ it follows that $\check H^i(R'/R,T_FQ)=0$ for $i>0$, using Proposition \ref{Pr:Hi-fpqc-hatW-sW} and Lemma \ref{Le:H*W-fpqc}.
\eproof

\bRe
\label{Re:Hi-fpqc-TFQ}
For an $\FF_p$-algebra $R$ and a faithfully flat ring homomorphism $R\to R'$ with semiperfect $R'$ we have $H^i_{\fpqc}(R,T_FQ)=\check H^i(R'/R,T_FQ)$, which depends only on $R$ by Corollary \ref{Co:Hi-fpqc-TFQ}. This avoids possible size issues in the definition of fpqc  cohomology.
\eRe

\subsection{Sheaves of infinitesimal elements}

Let $\bar Q=W(\bar\GG_a)$.
For a sheaf of rings $\AAA$ equipped with a homomorphism $\AAA\to\bar Q$ let $\AAA^0$ be the kernel of this map. This applies to $Q,Q\limperf,W,\sW$.

\bLe
\label{Le:Hi-fpqc-barQ}
For every $p$-nilpotent $R$ we have $H^i_{\fpqc}(R,\bar Q)=0$ for $i>0$.
\eLe

\bproof
Similar to Lemma \ref{Le:H*W-fpqc}.
\eproof

\bLe
\label{Le:sWR-barQr-surj}
If $R$ is good, then $\sW(R)\to\bar Q(R)=W(R_{\red})$ is surjective.
\eLe

\bproof
Source and target are derived $p$-complete by Corollary \ref{Co:sW-p-complete},
and the homomorphism $\sW(R)\to R_{\red}$ is surjective by Corollary \ref{Co:sWR-R-surj}.
\eproof

\bPr
\label{Pr:Hi-Qperf-Qperf0}
If $R$ is good, then $H^i_{\fpqc}(R,Q^{\perf,0})=H^i_{\fpqc}(R,Q\limperf)=0$ for $i>0$.
\ePr

\bproof
Let $R\to R'$ be a faithfully flat homomorphism of good rings.
We have exact sequences $0\to\hat W\to\sW\to Q\limperf\to 0$ (syntomic) and $0\to Q^{\perf,0}\to Q\limperf\to\bar Q\to 0$ (countably syntomic), both exact as presheaves on good rings, by Remark \ref{Re:hatW-sW-Qperf} and Lemma \ref{Le:sWR-barQr-surj}. Hence Proposition \ref{Pr:Hi-fpqc-hatW-sW} and Lemma \ref{Le:Hi-fpqc-barQ} give the result.
\eproof

\subsection{Some non-commutative cohomology groups and descent}

For a sheaf of rings $\AAA$ equipped with a homomorphism $\AAA\to\bar Q$ let $\GL_d(\AAA)^0$ be the kernel of $\GL_d(\AAA)\to\GL_d(\bar Q)$. 
We have an ind-fppf exact sequence of sheaves of groups 
\beqn
\label{Eq:MdTFQ-GLdsW-GLdW}
1\to 1+M_d(T_FQ)\to\GL_d(\sW)^0\to\GL_d(W)^0\to 1.
\eeqn

\bLe
\label{Le:conn-surj-Fp}
For a good $\FF_p$-algebra $R$, \eqref{Eq:MdTFQ-GLdsW-GLdW}
induces a surjective connecting homomorphism
\beqn
\label{Eq:conn-first}
\delta\colon 1+M_d(W^0(R))\to H^1_{\fpqc}(R,1+M_d(T_FQ)).
\eeqn
\eLe

\bproof
We can work with \v Cech cohomology for a chosen faithfully flat ring homomorphism $R\to R'$ where $R'$ is semiperfect; see Remark \ref{Re:Hi-fpqc-TFQ}. We have a commutative diagram of sets
\[
\xymatrix@M+0.2em{
M_d(Q^0(R)) \ar[r]^-\delta \ar[d]_f &
\check H^1(R'/R,M_d(T_FQ)) \ar[d]^g \\
1+M_d(Q^0(R)) \ar[r]^-\delta &
\check H^1(R'/R,1+M_d(T_FQ))
}
\]
where the upper line is the connecting homomorphism associated to the additive sequence 
\[
0\to M_d(T_FQ)\to M_d(Q^{\perf,0})\to M_d(Q^0)\to 0
\]
and the lower line is the connecting homomorphism associated to the multiplicative sequence 
\[
1\to 1+M_d(T_FQ)\to\GL_d(Q\limperf)^0\to\GL_d(Q)^0\to 1,
\]
$f$ is the bijective map $x\mapsto 1+x$, and $g$ is induced by the isomorphism $M_d(T_FQ)\cong 1+M_d(T_FQ)$, $y\mapsto 1+y$.
Both sequences are exact as presheaves over semiperfect $\FF_p$-algebras.
The diagram commutes by a simple calculation: 
For $a\in M_d(Q^0(R))$ with inverse image $b\in M_d(Q^{\perf,0}(R'))$, the element $c=b_1-b_0\in M_d(T_FQ(R'\otimes_RR'))$ represents $\delta(a)$, and $\delta(1+a)$ is represented by $(1+b_1)(1+b_0)^{-1}=(1+c)(1+b_0)(1+b_0)^{-1}=1+c$, using Lemma \ref{Le:TFQ-Qperf0}.
Now the upper line is surjective by Proposition \ref{Pr:Hi-Qperf-Qperf0}.
The proposition follows since $W(R)\to Q(R)$ is surjective by Remark \ref{Re:hatW-sW-Qperf}. 
\eproof

\bPr
\label{Pr:sW-descent-fpqc}
The category of finite projective $\sW(-)$-modules over good rings satisfies fpqc descent. 
\ePr

\bproof
Finite projective $\sW(-)$-modules form a prestack since $\sW$ is an fpqc sheaf. We have to show that every descent datum $(M,\alpha)$ of finite projective $\sW(-)$-modules for a faithfully flat ring homomorphism $R\to S$ with good $R$ is effective. Here $M$ is a finite projective module over $\sW(S)$, and $\alpha$ is a descent isomorphism over $\sW(S\otimes_RS)$. We can assume that $S/pS$ is semiperfect.

\subsubsection*{Step 0}
The descent datum of finite projective modules for $R\to S$ induced by $(M,\alpha)$ via $\sW\to\GG_a$ is effective and gives a finite projective $R$-module $\bar M_0$.

\subsubsection*{Step 1}
Reduction to free modules.
Let $\bar N_0$ be a finite projective $R$-module such that $\bar M_0\oplus\bar N_0$ is free. The ring $\sW(R)$ is henselian along the kernel of $\sW(R)\to R$ by Lemma \ref{Le:frame-Henselian}.
Hence $\bar N_0$ lifts to a finite projective $\sW(R)$-module $N_0$. Let $(N,\beta)$ the associated descent datum, thus $N=N_0\otimes_{\sW(R)}\sW(S)$. The modules $M\oplus N$ is free over $\sW(S)$ since its reduction over $S$ is free. If $(M,\alpha)\oplus(N,\beta)$ is effective, the same holds for the summands. Hence we can assume that $M$ is free over $\sW(S)$ and $\bar M_0$ is free over $R$. In this case the assertion is that $(M,\alpha)$ descends to a free $\sW(R)$-module.

\subsubsection*{Step 2}
Reduction to characteristic $p$.
Assume that there is an ideal $N\subseteq R$ of square zero such that the descent datum of finite free $\sW(-)$-modules for $R/N\to S/NS$ induced by $(M,\alpha)$ is effective. Let $M'_0$ be the resulting finite free $\sW(R/N)$-module.
Let $i\colon\Spec(R/N)\to\Spec R$ be the natural morphism. There is an exact sequence 
\[
0\to\hat W(\u N)\to\sW\to i_*(\sW)\to 0
\] 
of sheaves over $\Spec R$, which is exact on $S$-points since $S$ is good; see Proposition \ref{Pr:WN-sWR-sWRN}.
Hence, if we choose a basis of $M_0'$, the resulting basis of $M_0'\otimes_{\sW(R/N)}\sW(S/NS)$ lifts to a basis of $M$ over $\sW(S)$. In terms of this basis we have $\alpha\in 1+M_d(\hat W(\u N(S\otimes_RS)))$. Then $(M,\alpha)$ is effective by Lemma \ref{Le:H*hWN-fpqc}.
By induction it follows that we can assume that $R$ is an $\FF_p$-algebra.
\subsubsection*{Step 3}

The descent datum of finite projective $W(-)$-modules for $R\to S$ induced by $(M,\alpha)$ via $\sW\to W$ is effective and gives a finite projective $W(R)$-module $M''_0$. This module is free because its reduction over $R$ is free. The exact sequence of sheaves
\[
0\to T_FQ\to\sW\to W\to 0
\] 
is exact on $S$-points since $S$ is semiperfect. Hence, if we choose a basis of $M''_0$, the resulting basis of $M''_0\otimes_{W(R)}W(S)$ can be lifted to a basis of $M$ over $\sW(S)$. In terms of this basis we have $\alpha\in 1+M_d(T_FQ(S\otimes_RS))$. 
Lemma \ref{Le:conn-surj-Fp} implies that $(M,\alpha)$ is effective.
\eproof

\bCo
\label{Co:plain-sh-disp-good-fpqc}
Plain sheared displays over good rings satisfy fpqc descent.
\eCo

\bproof
Similar to Proposition \ref{Pr:Descent-Sh-disp-essp}, using Proposition \ref{Pr:sW-descent-fpqc} instead of Proposition \ref{Pr:Descent-LF(sW)-essemiperf}.
\eproof

\bCo
\label{Co:sheared-displays-fpqc-stack}
Sheared displays over $p$-nilpotent rings satisfy fpqc descent.
\eCo

\bproof
Clear from Corollary \ref{Co:plain-sh-disp-good-fpqc} and Definition \ref{De:Sh-disp}.
\eproof

\bCo
\label{Co:sheared-disp-over-good-rings-are-plain}
If $R$ is good, then all sheared displays over $R$ are plain (Definition \ref{De:Sh-disp}).
\eCo

\bproof
Every good ring has an fpqc covering by an essentially semiperfect ring. Now use Corollaries \ref{Co:essemiperf-plain} and \ref{Co:plain-sh-disp-good-fpqc}.
\eproof

\cclearpage

\end{document}